\documentclass[11pt,twoside]{amsart}
\usepackage[
  paperwidth=210mm,
  paperheight=280mm,
  inner=26.5mm,
  outer=22.6mm,
  includehead,
  top=21.8pt,
  headheight=40.8pt,
  headsep=12pt,
  bottom=29.13mm,
  footskip=41.4pt
]{geometry}

\usepackage{etoolbox}
\usepackage{xparse}
\usepackage{mathrsfs}
\let\Sec\S

\usepackage{clrcollabmacros}

\usepackage{amssymb}
\usepackage[all,cmtip]{xy}
\input diagxy
\usepackage{hyperref}
\usepackage{todonotes}

\newtheorem{theorem}{Theorem}[section]
\newtheorem*{theorem*}{Theorem}
\newtheorem{proposition}[theorem]{Proposition}
\newtheorem{corollary}[theorem]{Corollary}
\newtheorem*{corollary*}{Corollary}
\newtheorem*{proposition*}{Proposition}
\newtheorem{lemma}[theorem]{Lemma}

\theoremstyle{definition}
\newtheorem{definition}[theorem]{Definition}
\newtheorem{example}[theorem]{Example}
\newtheorem*{example*}{Example}
\newtheorem{notation}[theorem]{Notation}
\AtBeginEnvironment{example}{%
  \pushQED{\qed}%
}
\AtEndEnvironment{example}{\popQED\endexample}

\numberwithin{equation}{section}

\newtheorem{remark}[theorem]{Remark}

\DeclareMathOperator{\csk}{cosk}
\DeclareMathOperator{\Csk}{Cosk}
\DeclareMathOperator{\Dec}{Dec}

\DeclareMathOperator{\Image}{Im}
\DeclareMathOperator{\spine}{sp}
\DeclareMathOperator{\Tot}{Tot}
\DeclareMathOperator{\str}{str}
\DeclareMathOperator{\String}{String}
\DeclareMathOperator{\Spin}{Spin}

\newcommand{\cString}{\mathcal{S}tring}
\newcommand{\Linfft}{\LnA{\infty}^{\ft}}

\newcommand{\LGp}{\mathsf{Lie}_{\infty}\mathsf{Grp}}

\newcommand{\Nb}{\mathbb{N}}
\newcommand{\Rb}{\mathbb{R}}
\newcommand{\Id}{Id}

\newcommand{\sm}{\mathsf{Smooth}}
\newcommand{\ssm}{\mathsf{sSmooth}}

\renewcommand{\sSet}{\mathsf{sSet}}

\renewcommand{\C}{\mathsf{C}}
\newcommand{\D}{\mathsf{D}}
\newcommand{\sD}{\mathsf{sD}}
\newcommand{\sC}{\mathsf{sC}}

\newcommand{\LW}{\overline{W}}

\newcommand{\into}{\hookrightarrow}
\newcommand{\onto}{\twoheadrightarrow}

\newtheorem{pdef}[theorem]{Proposition-Definition}
\usepackage{epigraph}
\newcommand{\bl}{\bullet}

\usepackage{booktabs} 

\newcommand{\LnG}[1]{\mathsf{Lie}_{#1}\mathsf{Grp}}

\newcommand{\OS}{\Omega(S)}
\newcommand{\dR}{\mathrm{dR}}
\newcommand{\Om}{\Omega}
\newcommand{\cG}{\mathcal{G}}
\newcommand{\df}[1]{{\em #1}}

\newcommand{\afib}{\overset{\sim}{\twoheadrightarrow}}
\newcommand{\LE}[1]{L_{EA(#1)}}
\newcommand{\LB}[1]{L_{BA(#1)}}
\newcommand{\LEO}[2]{L_{EA(#1)/#2}}
\newcommand{\LBO}[2]{L_{BA(#1)/#2}}
\renewcommand{\odot}{\vee}
\newcommand{\jm}{\jmath}

\newcommand{\vph}{\vphi}
\newcommand{\cc}{\circ}

\renewcommand{\horn}{\Lam^m_k}
\newcommand{\hrn}[1]{#1(\horn)}
\newcommand{\jmm}{\jm^\ast}

\newcommand{\Ups}{\Upsilon}
\newcommand{\Ga}{\Gamma}
\newcommand{\ka}{\kappa}
\newcommand{\ph}{\phi}
\newcommand{\Del}{\Delta}

\newcommand{\eps}{\epsilon}

\newcommand{\ba}[1]{\bar{#1}}
\newcommand{\Lam}{\Lambda}
\newcommand{\lam}{\lambda}
\newcommand{\Hen}{\mathrm{Hen}}
\renewcommand{\cl}{\mathrm{closed}}
\newcommand{\et}{\acute{e}t}
\newcommand{\loc}{loc}

\let\oldaddcontentsline\addcontentsline
\newcommand{\stoptocentries}{\renewcommand{\addcontentsline}[3]{}}
\newcommand{\starttocentries}{\let\addcontentsline\oldaddcontentsline}

\tikzset{%
  clr style/.style={fill=white,inner sep=2.5pt},
  baseline=(current  bounding  box.center),
  ampersand replacement=\&, row sep=2em,column sep=2em,
  std/.style={->,font=\scriptsize}
}

\newenvironment{tikzdiag}[2]
{
\begin{tikzpicture}[clr style]
\matrix (m) [matrix of math nodes, row sep=#1em, column sep=#2em]
}
{
\end{tikzpicture}
}

\newcommand{\pbdiag}[1][.4]{\begin{scope}[shift=($(m-1-1)!{#1}!(m-2-2)$)]\draw +(-0.25,0) -- +(0,0)  -- +(0,0.25); \end{scope}}

\makeatletter
\NewDocumentCommand{\extp}{e{^}}{%
  \mathop{\mathpalette\extp@{#1}}\nolimits
}
\NewDocumentCommand{\extp@}{mm}{%
  \bigwedge\nolimits\IfValueT{#2}{^{\extp@@{#1}#2}}%
  \IfValueT{#1}{\kern-2\scriptspace\nonscript\kern2\scriptspace}%
}
\newcommand{\extp@@}[1]{%
  \mkern
    \ifx#1\displaystyle-1.8\else
    \ifx#1\textstyle-1\else
    \ifx#1\scriptstyle-1\else
    -0.5\fi\fi\fi
  \thinmuskip
}
\makeatletter

\newcommand{\Alt}{\textstyle{\extp}}
\newcommand{\ideal}{\trianglelefteq}
\DeclareMathOperator{\coker}{coker}
\DeclareMathOperator{\Coder}{Coder}
\newcommand{\dgco}{\mathsf{dgCoCom}}
\newcommand{\Bun}{\mathrm{Bun}}
\newcommand{\cT}{\mathcal{T}}

\makeatletter
\newcommand{\oset}[3][0ex]{%
  \mathrel{\mathop{#3}\limits^{
    \vbox to#1{\kern-2\ex@
    \hbox{$\scriptstyle#2$}\vss}}}}
\makeatother
\newcommand{\trivfib}{\oset[-.4ex]{\sim}{\twoheadrightarrow}}

\DeclareMathOperator{\GL}{GL}

\DeclareMathOperator{\Fun}{\mathsf{Fun}}
\DeclareMathOperator{\sSh}{\mathsf{sSh}}

\title{Lie's Third Theorem for Lie $\infty$-algebras}
\author{Christopher L.\ Rogers and Jesse Wolfson}
\address{Department of Mathematics \& Statistics, University of Nevada,  Reno.}
\thanks{CLR was supported in part by NSF Grant DMS-2305407 and Simons Foundation grant 585631}
\email{chrisrogers@unr.edu, chris.rogers.math@gmail.com}
\address{ Department of Mathematics, University of California-Irvine}
\email{wolfson@uci.edu}

\thanks{JW was supported in part by NSF Grants DMS-1811846 and DMS-1944862.}
\begin{document}
	
\begin{abstract} 
We introduce the theory of local minimal models for Kan simplicial manifolds, which provide the appropriate generalization of minimal Kan simplicial sets to geometric contexts. We use this to obtain the first proof of Lie’s third theorem for finite-type Lie $\infty$-algebras: Every finite-type, homologically and non-negatively graded $L_\infty$-algebra over $\R$ integrates to a finite-dimensional Lie $\infty$-group. As a corollary, our construction yields a new explicit finite-dimensional model for the string Lie 2-group. 
\end{abstract}
\maketitle

\tableofcontents

\vfill
\epigraph{``It seems likely that no `simple' proof of Lie's Third Theorem exists.'' J-P.\ Serre \cite[p. 155]{Serre:Lie} }

\newpage
\section{Introduction}
Lie $\infty$-algebras are the coherent homotopy analogs of differential chain graded Lie algebras concentrated in non-negative degrees, i.e.\ connective homologically graded $L_\infty$-algebras\footnote{We recall the precise definition in \Sec \ref{s:lien}.}.  We say such an algebra is {\bf finite-type} if it is finite-dimensional in each degree. Any finite-dimensional Lie algebra over a field of characteristic zero is an example. Finite-type Lie $\infty$-algebras satisfying a nilpotency condition are featured implicitly, yet prominently, in Sullivan's approach to rational homotopy theory \cite{Sull}. In particular, the ``Koszul dual'' of a minimal model for a finite-type nilpotent space is a nilpotent Lie $\infty$-algebra. Sullivan's theory provides foundational evidence for considering Lie $\infty$-algebras -- specifically, those beyond the classical nilpotent Lie algebras -- as legitimate infinitesimal avatars of global topological objects. Furthermore, Sullivan's work, along with the work of Gabriel and Zisman \cite{GZ:1967}, and Barratt, Gugenheim, and Moore \cite{BGM} on simplicial fibrations, demonstrates the importance of having a {\it good theory of minimal models}, especially when the goal is to extract geometric information from homotopy-theoretic constructions. 
It is therefore somewhat ironic that the full strength of Sullivan's theory, in its original manifestation, fails to include the examples that were of primary interest to Lie, Killing, and Cartan, e.g., the complex simple Lie algebras and their compact  real forms.         

A key ingredient in Sullivan's theory is an  ``integration'' procedure which leverages simplicial de Rham theory to produce a Kan simplicial set from a nilpotent Lie $\infty$-algebra. Building on this, Henriques \cite{Hen} proved in 2008 the existence of an integration functor 
\begin{equation} \label{eq:intro-int}
\sint\colon \Linfft\to \LGp
\end{equation}
from the category of {\it all} finite type Lie $\infty$-algebras over $\R$ to the category of infinite dimensional {\bf Lie $\infty$-groups}, i.e.\, reduced Kan simplicial Banach manifolds (Def.\ \ref{d:inftygpd}).\footnote{Simultaneously, Getzler \cite{Get} gave a detailed study of the integration of nilpotent $L_\infty$-algebras, which foreshadows some of the present work and which we discuss below.} By later joint work of the first author  \cite{RZ} with Zhu, \eqref{eq:intro-int} extends canonically to a ``partially left-exact'' $\infty$-functor, where the source is presented as an $\infty$-category via its category of fibrant objects structure, as constructed by the first author in \cite{R}, and the target, via its structure of an incomplete category of fibrant objects as in \cite{RZ}.\\ 

\paragraph{{\bf The String Group and Kac-Moody Theory}}
Henriques' work initiates the connection between the theory of classical (non-nilpotent) Lie algebras and Sullivan's formalism. Moreover, it brings into focus a different class of Lie $\infty$-algebras with connections to mainstream representation theory, specifically the theory of Kac-Moody Lie algebras. Indeed, given a compact, simple and simply connected Lie group $G$, Henriques functorially constructs a family of Lie 2-groups $\cString_k(G):=\tau_{\leq 2}\sint \str_k(\g)$, one for each $k \in \Z$, by applying Duskin's 2-truncation \cite{Duskin} to the integration of $\str_k(\g)$. Here, $\str_k(\g)$ is the level $k$ {\bf string Lie 2-algebra}. As we recall in Example \ref{ex:s2g}, this is a natural example of a 
{\it minimal} finite-dimensional Lie $\infty$-algebra: its underlying complex is $\g \dsum \R[1]$ with vanishing differential, and its $L_\infty$ structure is completely determined by $k \in H^3_{\CE}(\g,\R)$. The upshot is that the geometric realization of $\cString_k(G)$ produces an extension of the classifying space $BG$ by $K(\Z,3)$. In particular, for $k=1$ and $G=\Spin(n)$, this provides a Lie 2-group model $\cString(n)$ for the classifying space of the {\bf string group} $\String(n)$, i.e.\ the 3-connected topological group covering $\Spin(n)$. Equivalently, one can understand $\cString(n)$ as an explicit geometric, albeit infinite-dimensional, model for the homotopy fiber of the {\bf  universal first fractional Pontryagin class} $\frac{1}{2} p_1 \maps B\Spin(n) \to K(\Z,4)$.        

As emphasized in \cite{BSCS} and \cite{W:2012} (see also \cite[Rmk.\ 8.7]{Hen}), the infinite-dimensional Lie 2-group $\cString_k(G)$ melds the data of the level-$k$ Kac-Moody central extension of the loop group $LG$ with a presentation of the $k$th power $\mathscr{G}^{\boxtimes k}$ of the canonical $U(1)$-gerbe $\mathscr{G}$ over $G$. The relationship between these objects over $LG$ and $G$, respectively, is the subject of an established line of research spanning several decades. For example, Brylinski posed as an open question thirty years ago  whether one could obtain the positive energy representations of the Kac-Moody central extension geometrically from the canonical gerbe on $G$ \cite[Sec. 6.4]{B:1993}. As he noted, this was to be contrasted with the infinite-dimensional analog of the Borel-Weil Theorem developed in \cite{PS:1986}. Recall that the positive energy representations admit, via a canonical circle action, a decomposition into {finite-dimensional} eigenspaces, with the irreducibles parameterized by a subset of the integral weights for the finite-dimensional Lie algebra $\g$. We therefore understand Brylinski's question as one which asks for the integrable highest-weight modules of the level-$k$ affine Lie algebra $\widehat{L_k\g}$ via finite-dimensional geometric constructions involving only the compact group $G$, in analogy with classical Lie theory for $\g$. Related to this (e.g.\ see the exposition in \cite{F:2001}), and going back even further, Chern and Simons explain in \cite{CS:1974} that the starting point of their work on secondary characteristic classes was an attempt to derive an explicit combinatorial formula for the first Pontryagin number.

Despite a substantial amount of progress on the above circle of ideas, an explicit finite-dimensional model of $\cString_k(G)$ has remained elusive, with the best attempts prior to this present work appearing in the results of Schommer-Pries \cite{SP}, Waldorf \cite{W:2012}, and work of the second author \cite{W}.\footnote{Note that all of these are essentially conditional on an explicit integral cocycle for half the first Pontryagin class, which, as noted below, is currently lacking in the literature.} The theorem of Freed, Hopkins, and Teleman \cite{FHT:2011}, and also Meinrenken's related work \cite{M:2012}, establish the relationship between the twisted $K$-theory of $G$, the Dixmier-Douady class of $\mathscr{G}^{\boxtimes k}$ and the Verlinde ring $R^k(LG)$. Nevertheless, Brylinski's question with its implicit emphasis on finite-dimensional methods remains open. Finally, the work of Brylinski and MacLaughlin \cite{BM:1996} completes the program initiated by Chern and Simons on cocycle models for characteristic classes, specifically, the Chern, Euler, and Pontryagin classes. On the other hand, an analogous explicit construction for $\frac{1}{2}p_1$ has yet to appear.\\

\paragraph{\bf{Main Results (Prelude)}}
The above discussion motivates the need for a true Lie theoretic package for Lie $\infty$-groups and Lie $\infty$-algebras. This is the goal of the present paper, current work in progress, and subsequent work by the authors. Our focus here concerns extending the following cornerstone of classical Lie theory:
\begin{quote}
{\bf Lie's 3rd Theorem:} Every finite-dimensional Lie algebra over $\R$ integrates to a finite-dimensional Lie group. 
\end{quote}
Our main result is:
\begin{theorem}[Lie's 3rd Theorem for Lie $\infty$-algebras]\label{t:lie3}
  Every finite type Lie $\infty$-algebra over $\R$ integrates to a finite-dimensional Lie $\infty$-group. Concretely,
  there exists a finite dimensional Lie $\infty$-group $\cG^{\ft}(L)$ and a weak equivalence
  \begin{equation} \label{eq:intro-finite}
    \cG^{\ft}(L) \xto{\sim} \sint L,
  \end{equation}
  where $\int L$ denotes Henriques' smooth Sullivan integration of $L$ as in Eq.\ \ref{eq:intro-int}.
\end{theorem}
The above theorem has been a long-standing desideratum in higher Lie theory. It can be interpreted as an extension of the results of Getzler \cite{Get} for the nilpotent case, and a globalization of the local integration featured in \cite{SevSir}. We in particular obtain an explicit finite-dimensional Lie 2-group model for $B\String(n)$:
\begin{corollary}[Explicit finite-dimensional String 2-Group] \label{cor:FinString}
Let $\g$ be a simple Lie algebra of compact type and $k \in \Z$. Let $G$ denote the simply connected Lie group of $\g$, and let $U$ be a good cover of $G$. There exists a finite-dimensional Lie 2-group $\cString^{\ft}_k(G, U)$ with 1-simplices given by $\cString^{\ft}_k(G, U)_1=U$, and a weak equivalence 
\[
\cString^{\ft}_k(G, U) \xto{\sim} \cString_k(G).
\]
\end{corollary}
The explicit description of $\cString^{\ft}_k(G, U)$, and its dependence on the cover $U$ and two auxiliary choices is spelled out in   
Example \ref{ex:s2g}. Furthermore, from a careful analysis of our construction for $G=\Spin(n)$ and $k=1$, we see that it provides the sought-after universal cocycle model for the fractional Pontryagin class $\frac{1}{2}p_1$, as discussed above. We spell out the details of this observation in forthcoming work.\\

\paragraph{\bf{Locally Minimal Models}} Our main technical innovation in the above is inspired by the use of minimal models in classical homotopy theory. As already mentioned, they are a core tool in the basic theory of Kan fibrations. Furthermore, on the algebraic side, the minimal model of the de Rham algebra of a nilpotent space represents another key ingredient, along with integration, of Sullivan's theory. What ties these two notions of minimal model together is higher Lie theory: namely, Getzler's work on the integration of nilpotent $L_\infty$-algebras \cite[Prop.\ 5.7]{Get}.  

\newcommand{\relb}{\text{rel $\pa \Del$}}
Recall that a reduced Kan complex $X_\bl \in \sSet_0$ is minimal if, for all $n > 0$, the canonical map $X_n \to X_n/ \simeq_{\relb}$ is an isomorphism. Here, $\simeq_{\relb}$ denotes the equivalence relation corresponding to homotopy relative to the boundary $\pa \Del^n$. Naively, one might try to lift this definition directly to the diffeo-geometric context for Lie $\infty$-groups. Of course, since homotopy rel boundary may not lead to a smooth quotient, the definition needs to be phrased in the category of sheaves $\Sh(\sm)$. Here $\sm$ denotes the site of smooth Banach manifolds with surjective submersions as covers. Note that this category of sheaves has been in the background all along, as it is used in the precise definition  of the Kan horn filling conditions (Def.\ \ref{d:inftygpd}) for a simplicial manifold. 
Other than that, everything goes through and we say $X_\bl \in \LGp$ is {\bf minimal} if for all $n > 0$ the map $X_n \to X_n/ \simeq_{\relb}$ is an isomorphism of sheaves. Furthermore, a {\bf minimal model} for $Y_\bl \in \LGp$ is a minimal Lie $\infty$-group $X_\bl$ equipped with a weak equivalence $X_\bl \xto{\sim} Y_\bl$. 
Unfortunately, yet interestingly, this naive notion of minimal model quickly leads to a dead-end:
\begin{proposition*}[Minimal Model No Go Theorem] 
Let $\g$ be a simple Lie algebra of compact type, $G$ its simply connected Lie group, and $k \in \Z$. There is no minimal model for the Lie 2-group $\cString_k(G) = \tau_{\leq 2}\sint \str_k(\g)$.  
\end{proposition*}   
As we explain in detail in Example \ref{ex:nomin}, such a minimal model $X$ for $\cString_k(G)$ would be an extension of $NG$, i.e.\ the nerve of $G$ as a one object Lie groupoid, by the smooth Eilenberg-Mac Lane space $K(S^1,2)$. Furthermore, the extension would necessarily be classified by a nontrivial element in $H^3_{naive}(G;S^1)$, the naive continuous group cohomology of $G$. But it is well-known that $H^3_{naive}(G;S^1)=0$, and this yields a contradiction.

Although the above definition of minimal model is still useful, this negative result forces us to a consider a weaker notion. We say a morphism $F \to F'$ in $\Sh(\sm)$ is {\bf surjective \'etale} if it is a surjective morphism of sheaves, and for all sections $U \to F'$, the pullback $F \times_{F'} U \to U$ is represented by a (surjective) local diffeomorphism of smooth manifolds. We then define $X_\bl \in \LGp$ to be {\bf locally minimal} if $X_n \to X_n/ \simeq_{\relb}$ is surjective \'etale for all $n > 0$, and we define a {\bf locally minimal model} for $Y_\bl \in \LGp$ in the analogous way. Since we need to work relatively throughout much of the paper, we also introduce the notion of a locally minimal model for a Kan fibration between simplicial manifolds. See Def.\ \ref{d:locmin} the precise definition. Furthermore, we show in Lemma \ref{l:locmin} that any {\bf \'etale hypercover} between simplicial manifolds is a locally minimal fibration. The precise definition of \'etale hypercover is given in Def.\ \ref{def:etale-maps}. For our purposes here, we emphasize that such a hypercover is necessarily a level-wise surjective local diffeomorphism.    

In order to appreciate that this notion of locally minimal model is on the right track, let us first recall that the simplicial homotopy group $\pi_n(X)$ of a Lie $\infty$-group $X_\bl$ is, by definition a certain coequalizer, e.g.\ \cite[Eq.\ 5.2]{RZ}). In general, it is not a manifold. However, it is often the case that $\pi_n(X)$ is isomorphic in $\Sh(\sm)$ to a quotient of a representable group object, i.e.\ a quotient of a Lie group $H$ by a not-necessarily discrete finitely generated subgroup $A \leq H$. In this situation, following \cite[p. 1034]{Hen}, we say $\pi_n(X)$ is a ``finite-dimensional diffeological group''. Here is the main example of such a Lie $\infty$-group: 
\begin{example*}
By \cite[Thm.\ 6.4]{Hen}, if $L$ is a finite-type Lie $\infty$-algebra, then the simplicial homotopy groups of its integration $\sint L$, as in \eqref{eq:intro-int}, are finite-dimensional diffeological groups. 
\end{example*}
In light of this fact, we have the following  encouraging result:
\begin{proposition*}[Finite-dimensional from locally minimal] 
Suppose $Y_\bl$ is a Lie $\infty$-group such that $\pi_n(Y)$ is a finite-dimensional diffeological group for all $n \geq 1$. If $X_\bl \xto{\sim} Y_\bl$ is a locally minimal model for $Y_\bl$, then $X_\bl$ is a finite-dimensional Lie $\infty$-group. 

Moreover, the dimensions $\{\dim X_n\}_{n\ge 0}$ are invariants of $Y_\bullet$ in the sense that for any weakly equivalent Lie $\infty$-group $Z_\bullet\simeq Y_\bullet$ and any locally minimal model $W_\bl\xto{\sim} Z_\bullet$, $\dim X_n=\dim W_n$ for all $n$.
\end{proposition*}
The reason for the finite-dimensionality of $X_\bl$ is conceptually straightforward, as we now explain. First, the proof of the above proposition involves walking inductively up the Moore-Duskin Postnikov tower (Section \ref{s:post}) of $X_\bl$. The local minimality of $X_\bl$ implies that the map connecting the $(n+1)$-st Moore stage with the $n$-th Duskin stage is an \'etale hypercover. As discussed above, this implies that it is a level-wise surjective local diffeomorphism.  On the other hand, the assumption on the simplicial homotopy groups leads one to conclude that the map 
connecting the $n$-th Duskin stage with the $n$-th Moore stage is level-wise a fiber bundle whose fibers are finite-dimensional non-Hausdorff manifolds. By combining these facts along with the inductive hypothesis for the $n$-th Moore stage and the implicit function theorem, we conclude that the $(n+1)$-st Moore stage is level-wise finite-dimensional. The remaining details are given in the proof of Lemma \ref{l:locminfd}. \\

\paragraph{{\bf Main Results Revisited: Constructing $\cG^{\ft}_\bl(L)$}} 
From the above discussion, it is clear that in order to prove Theorem \ref{t:lie3} and exhibit a finite-dimensional integration $\cG^{\ft}_\bl(L)$, it suffices to construct a locally minimal model for Henriques' integration $\sint L$. Our approach exploits the interplay between the Moore-Duskin Postnikov tower for $\sint L$ and the analogous tower (Section \ref{sec:Linf-post}) for $L$. The work of the second author on functorial $k$-invariants \cite{W2}, and the smooth relative Duskin truncation \cite{W} provides the essential tools for carrying this out. 

Our construction consists of three steps, which we package into the following three theorems. For the purposes of exposition, we postpone the list of needed technical ingredients until the end of this introductory section.

Since the zeroth homology of any Lie $\infty$-algebra is a Lie algebra $\g:=H_0(L)$, we begin by exhibiting a locally minimal model for $\sint \g$. 
The smooth Sullivan integration of a Lie algebra is, in a sense, a classical construction i.e.\ \cite[``Thm 8.1'']{Sull}. Hence, it is interesting to note that this first step of our construction turns out to be the most delicate. 

Let us say that $X_\bl \in \ssm $ is a {\bf good simplicial manifold} if for all $n \geq 0$, we have $X_n\simeq \pi_0 X_n$. In other words, $X_n$ is a disjoint union of contractible pieces. As usual, we say $X_\bl$ is {\bf reduced} if $X_0=\ast$. Given a Kan fibration $f \maps X\to Y$, we denote by $\tau_{\le n}(X,f)\to Y$ its $n$-th {\bf Duskin truncation}. See Section~\ref{s:post} for its precise definition.     

\begin{theorem}[Section \ref{s:base}]\label{t:main1} 
    Let $L$ be a finite type Lie $\infty$-algebra and $G$ the simply-connected Lie group integrating $\g:=H_0L$.  Let $NG$ denote the nerve of $G$ as defined above.
    \begin{enumerate}
        \item\label{t:main1.1} There exists a good, reduced, \'etale hypercover $U_\bl \trivfib NG$; in particular $U_\bl$ is locally minimal.
        \item\label{t:main1.2} Given a good, reduced, \'etale hypercover $U_\bl \trivfib NG$, there exists a commuting square
        \begin{equation}\label{e:Lie30}
            \xymatrix{
                \cG^0_\bl \ar[r]^\sim \ar[d] & \sint \g \ar[d]^{\tau_{\le 1}} \\
                U_\bl \ar@{->>}[r]^\sim & NG
            }
        \end{equation}
        in which the left vertical map is a minimal Kan fibration, and $\cG^0_\bl$ is good, reduced, and finite dimensional Lie $\infty$-group. Equivalently, the Kan fibration 
        \[
            \sint\g \to NG
        \]
        admits a good minimal model when pulled back along $U_\bl \trivfib NG$.  
    \end{enumerate}
\end{theorem}

As we recall in Prop.\ \ref{prop:int-exact}, Henriques' functor $\sint(-)$ sends a large class of fibrations in $\Linfft$  to Kan fibrations in $\LGp$. In particular, the canonical surjection $L \to \g=H_0(L)$ integrates to a Kan fibration $\sint \tau_{\le 0} \maps \sint L \to \sint \g$. The next step in our construction is to exhibit a local minimal model for this fibration. In fact, in this case we can produce a minimal model in the ``naive'' sense discussed earlier. This underscores the subtlety of the previous step.
\begin{theorem}[Section \ref{sec:main3-pf}]\label{t:main3}
The Kan fibration $\sint\tau_{\le 0}\maps \sint L\to \sint \g$ admits a minimal model
\begin{equation*}
  \xymatrix{
    \mathcal{L}_\bl \ar[r]^\sim \ar[dr] & \sint L \ar[d]^{\int \tau_{\le 0}}\\
    & \sint \g
  }
\end{equation*}
\end{theorem}
For the last step, we combine the locally minimal models $\cG^{0}_\bl \xto{\sim} \sint \g$ and $\mathcal{L}_\bl \xto{\sim} \sint L$ from the previous two. Again, the key idea is a careful analysis of the principal bundles appearing in the Moore-Duskin Postnikov tower.  
\begin{theorem}[Section \ref{sec:main4-pf}] \label{t:main4} 
Given $\cG^0_\bl \xto{\sim} \sint\g$ and $\mathcal{L}_\bl \xto{\sim} \sint L$
as above, the pullback 
\[
\cG^{\ft}(L)_\bl:=\cG^0\times_{\int \g}\mathcal{L} 
\]
is a good, locally minimal, and finite-dimensional Lie $\infty$-group. In particular, the canonical weak equivalence $\cG^{\ft}(L)_\bl \xto{\sim} \sint L$ 
presents $\cG^{\ft}(L)_\bl$ as a finite dimensional locally minimal model for $\sint L$. 
\end{theorem}
At this point, it is worthwhile to see how the output of construction compares 
to the minimal model for $\sint L$ constructed by Getzler for the special case when $L$ is nilpotent.  
\begin{example}
Let $L$ be a finite type Lie $\infty$-algebra such that $\g:=H_0L$ is a nilpotent Lie algebra. Theorems \ref{t:main1}, \ref{t:main3}, and \ref{t:main4} simplify considerably in this case: first, the nilpotency of $\g$ implies that the exponential map is an isomorphism of manifolds $\g\xto{\cong} G$, and thus the identity map $NG\xto{\id} NG$ is itself a good, reduced \'etale hypercover. Further, by \cite[Example 5.5]{Hen} and \cite[Prop.\ 6.7]{RZ}, the map $\sint\g\trivfib NG$ is a hypercover. Therefore, the left vertical map in the square~\eqref{e:Lie30} is a minimal trivial fibration, and thus an isomorphism (Lemma~\ref{l:minkan}). Theorem~\ref{t:main3} therefore provides a finite dimensional {\em minimal} model (not just locally minimal)
    \[
            \xymatrix{
                \cG\ar[r]^\sim \ar[dr]_{q} & \sint L \ar[d]^{\tau_{\le 1}} \\
                    & NG
            }.
    \]
When $L$ is nilpotent, this gives a minimal model for Getzler's integration $\gamma_\bullet(L)$ of \cite{Get}. When $L$ is minimal, $\gamma_\bullet(L)$ is also minimal \cite[Prop.\ 5.7]{Get}, and, using the theory of $k$-invariants developed below, one can show the two models are (non-canonically) isomorphic.
\end{example}

\mbox{}

\paragraph{{\bf Main Results Revisited: Lie's 3rd Theorem for Lie $n$-algebras}}
Given a finite type Lie $n$-algebra $L$, i.e.\  a Lie $\infty$-algebra concentrated in the first $n$ non-positive degrees, one might hope to improve our main theorem and an obtain a finite dimensional integrating Lie $n$-group. In some cases, this is indeed achievable, for example when $L$ is the string Lie 2-algebra $\str_k(\g)$. As discussed by Henriques \cite{Hen}
there are, in general, obstructions preventing a truncated infinite-dimensional integration \cite[Thm.\ 7.5, Ex.\ 7.10]{Hen}. More precisely, Henriques constructed a natural map
\[
    \partial_n\maps \pi_{n+1} G\to H_{n-1}L
\]
where $G$ is the simply connected Lie group integrating $H_0(L)$, and he 
proved that $\tau_{\le n}\int L$ is representable if and only if $\partial_n$ has discrete image. This same issue arises in our context for constructing a finite-dimensional Lie $n$-group.

In Sec.\ \ref{s:lie3n}, we introduce a notion of Lie $n^\ast$-groupoid, extending to the $n > 1$ case the notion of effective Weinstein groupoid introduced by Tseng-Zhu \cite{TZ} in their study of integration of Lie 1-algebroids. 
See Section \ref{s:lie3n} for the precise definition. In particular,  a Lie $n^\ast$-group is, roughly speaking, a Lie $(n+1)$-group which is an $n$-group up to some discrete ambiguity in the $(n+1)$-cells. As established in \cite[Thm.\ 1.2]{TZ}, every Lie algebroid integrates to an effective Weinstein groupoid as defined in \cite[Sec.\ 4.2]{TZ}, and, thus, by Zhu's nerve construction \cite[Thm.\ 1.4]{Z} to a Lie $1^\ast$-groupoid in our sense. We establish an analogue of this for finite type Lie $n$-algebras.

\begin{corollary}[Lie's 3rd Theorem for Lie $n$-algebras\footnote{In the early stages of this work, the authors incorrectly claimed a stronger result.  Fixing this led to the present definition of Lie $n^\ast$-group.}]\label{c:lien*}
Let $n>0$ and let $L$ be a finite type Lie $n$-algebra.  Let $\cG^{\ft}_\bl(L) \xto{\sim} \sint L$ be a locally minimal model as in Thm.\ \ref{t:main4}. 
\begin{enumerate}
\item The Moore truncation $\tau_{<n+1}\cG^{\ft}(L)_\bl$ is a finite dimensional Lie $n^*$-group, and a locally minimal model for $\tau_{<n+1}\sint L$.  

\item If $\partial_n(\pi_{n+1}G)\subset H_{n-1}L$ is a discrete subgroup, then $\tau_{\le n}\cG^{\ft}(L)_\bl$ is a finite dimensional Lie $n$-group, and a locally minimal model for $\tau_{\le n}\sint L$.
\end{enumerate}
\end{corollary}
\mbox{}

\paragraph{{\bf Towards Lie's 2nd Theorem for Lie $n$-algebras}}
With our main theorem in hand, it is natural to ask for the full Lie theory package, i.e. a differentiation functor $G\mapsto \g$, along with a ``Lie's Second Theorem'' establishing that it is fully faithful on a subcategory of suitably connected higher Lie groups. Recent work by the first author \cite{Diff}  provides such a differentiation functor with good homotopical properties, along with explicit canonical identifications to the classical differentiation of simplicial Lie groups. Prior to this, a tractable differentiation functor for Lie $\infty$-groups had been a persistent stumbling block in the literature.

In the sequel to the present work, we show that the differentiation functor constructed in \cite{Diff} is left inverse to the integration functor considered here, and using this, we establish the expected analogue of Lie's Second Theorem.\\

\paragraph*{\bf Outline of Paper}
In this final subsection, we guide the reader through the main technical developments needed to prove Theorems \ref{t:main1}, \ref{t:main3}, \ref{t:main4}, and thereby our main results Theorem \ref{t:lie3} and Corollary \ref{cor:FinString}.

Henriques' results stem from a careful analysis of the Postnikov tower of a finite type Lie $\infty$-algebra, and we build on this approach here.  Recall that classically, Postnikov theory analyzes a connected space $X$ in terms of a tower of fibrations
\[
    X\to \cdots\to X^{(n)}\to X^{(n-1)} \to \cdots\to X^{(1)}
\]
where $X^{(1)}$ is (the classifying space of) the fundamental groupoid of $X$ and the homotopy fiber of $X^{(n)}\to X^{(n-1)}$ is a $K(\pi_n X,n)$. Moreover, each stage of the tower can be classified by a ``$k$-invariant'' sitting in a homotopy pullback
\[
    \xymatrix{
     X^{(n)} \ar[rr] \ar[d] && EK(\pi_nX,n)//\pi_1 X\ar[d] \\
     X^{(n-1)} \ar[rr]^-{k_n} && K(\pi_nX,n+1)//\pi_1X
    }.
\]
Henriques' proved \cite[Thm.\ 6.4]{Hen} that the homotopy sheaves of his integration $\sint L$ are finite dimensional diffeological groups, i.e. quotients of finite dimensional Lie groups by finitely generated subgroups. This suggests that, given an analogous theory of $k$-invariants for Lie $\infty$-groups, one could construct a finite dimensional model by carefully inducting up the Postnikov tower. We carry out this approach here, constructing $k$-invariants for Lie $\infty$-groups and $L_\infty$-algebras.  

Whereas $k$-invariants for spaces can be constructed using the axiom of choice (cf. \cite{BGM}), this breaks down in geometric settings.  In Section~\ref{s:post}, we recall the Moore-Duskin Postnikov tower and establish the relevant properties for smooth $\infty$-groups. We recall the theory of minimal fibrations in this setting.  We then build on recent work of the second author \cite{W2} to establish the existence of {\em functorial} local $k$-invariants of a smooth $\infty$-group (see Theorem~\ref{t:kinvar}); for readability, we defer the full proof, along with the explicit local formulas to Section~\ref{app:kinvar}. We use this to analyze the geometry of the stages in the (Moore-Duskin) Postnikov tower of a smooth $\infty$-group. Note that this relies only minimally on the geometric context at hand, and we develop this theory in the generality of a ``category with covers'' as in \cite{W,RZ,W2}. 

In Section~\ref{s:lien}, we carry out an analogous development of $k$-invariants for $L_\infty$-algebras (see Section~\ref{sec:Linf-post}).  We build on work of the first author \cite{R} and of the first author and Zhu \cite{RZ} to analyze how these $k$-invariants behave under integration, and we explicitly relate them to the global $k$-invariants developed in Section~\ref{s:post} (see Example~\ref{ex:LEA-int}). Again, for readability we defer full proofs to Section~\ref{app:Linf}, where we establish the $k$-invariants as a consequence of Quillen's classification theorem for principal dg coalgebra bundles \cite[\Sec B.5]{Quillen}, and a coalgebraic analog of the associated bundle construction developed by Prigge \cite[\Sec 2]{Prigge}.

With good theories of $k$-invariants in hand for the infinitesimal and global settings, we turn our attention to differential topology proper and the setting of Lie $\infty$-groups. In Section~\ref{s:liegroup}, we assemble the results we will need on Lie $\infty$-groups. We recall the formalism of descent categories introduced by \cite{BG}, and show that the category $\Sh(\sm)$ of sheaves on smooth manifolds with surjective local diffeomorphisms of sheaves gives an example (Lemma~\ref{l:desccat}). As we observe in Remark~\ref{r:desc<=>catcov}, a careful reading of \cite{W} shows that all of the results we established in Section~\ref{s:prelim} and in Sections~\ref{s:post1} and~\ref{s:min} carry over unchanged to the setting of descent categories (as one would expect from \cite{BG,W2}). Further, with the addition of a single assumption, which holds in our intended application (see Remark~\ref{r:discgroupkinvar}), the same applies to the results of Section~\ref{sec:smooth-kinv}.  With this in hand, we are now in a position to develop the framework of locally minimal models in analogy with the treatment of minimal models above. Then, using Quinn's theory of manifold stratified spaces \cite{Q} at a key step, we develop a (partial) analog of Artin and Mazur's \cite{AM} theory of split hypercovers (in Proposition~\ref{p:hypexist}).  This suffices to give the first statement in Theorem~\ref{t:main1}.

In Section~\ref{s:base}, we use the results of Sections~\ref{s:liegroup} and~\ref{s:post} to establish statement \eqref{t:main1.2} of Theorem \ref{t:main1} 
and hence Theorem \ref{t:main4} for the special case of an ordinary Lie algebra $\g$.  This provides the base case for our inductive proof of Theorem \ref{t:main3} and Theorem \ref{t:main4}, which we carry out in Section \ref{sec:main3-pf}, and 
Section \ref{sec:main4-pf}, respectively. Together, this completes the proof of our main theorem. 

In Section~\ref{s:lie3n}, we introduce the notion of a Lie $n^\ast$-group/oid, we prove Lie's Third Theorem for Lie $n$-algebras (Corollary~\ref{c:lien*}).  We close by explicitly working out two examples. The first, Example~\ref{ex:s2g}, is the content of Corollary \ref{cor:FinString}. We give a finite dimensional locally minimal model of the string 2-group of a simple Lie algebra of compact type. The second, Example~\ref{ex:h2g}, describes the finite dimensional Lie $2^\ast$-group integrating Henriques' \cite[Example 7.10]{Hen} of a Lie 2-algebra with a nontrivial obstruction to integrating it to a 2-group.

Necessary background is recalled in Section~\ref{s:prelim}, detailed proofs of the construction of $k$-invariants for smooth $\infty$-groups are deferred to Appendix~\ref{app:kinvar}, while detailed proofs of the same for $L_\infty$-algebras are deferred to Appendix~\ref{app:Linf}.  

Our exposition contains many detailed examples.  In an effort to aid readability, we use a $\blacklozenge$ symbol to denote the conclusion of an example, analogous to how the symbol $\qed$ denotes the end of a proof.

\paragraph*{\bf Acknowledgments}
CLR thanks the organizers and participants of the 2022 BIRS workshop on Poisson Geometry, Lie Groupoids and Differentiable Stacks for their feedback on the early stages of this work. JW thanks Andr\'e Henriques for helpful correspondence, and Ezra Getzler for introducing him to this question and for many helpful conversations.

\section{Preliminaries}\label{s:prelim}

\subsection{Categories with covers}\label{s:catcovs}
As in \cite{W,W2}, we work in a (sub-canonical) {\em category with covers}, i.e. we fix a (locally small) category $\C$ with a subcategory of ``covers'' satisfying the following axioms
\begin{enumerate}
    \item\label{a:term} $\C$ has a terminal object $\ast$ and for all $X\in\C$, the map $X\to\ast$ is a cover,
    \item\label{a:basechange} pullbacks of covers along arbitrary maps in $\C$ exist and are covers,
    \item\label{a:rightcancellation} if $f$ and $gf$ are covers, then so is $g$,
    \item\label{a:subcan} covers are effective epimorphisms, i.e. for a cover $f\colon X\to Y$, the diagram
    \[
        X\times_Y X\rightrightarrows X\to^f Y
    \]
    is a coequalizer.
\end{enumerate}

We record some elementary consequences of the axioms. 
\begin{lemma}\label{l:prods}\mbox{}
   $\C$ has finite products and projections along factors are covers.
\end{lemma}
\begin{proof}
    By Axiom~\ref{a:term} and Axiom~\ref{a:basechange}, the pullback square
    \[
        \xymatrix{
            X\times Y \ar[r] \ar[d] & Y \ar[d]\\
            X \ar[r] & \ast
        }
    \]
    exists in $\C$ and all maps in it are covers. The general result now follows by induction.
\end{proof}

\begin{lemma}\label{l:covepi}\mbox{}
    \begin{enumerate}
        \item Let $f\colon X\to Y$ and $g\colon Y\to X$ be covers such that $gf=1_X$.  Then $fg=1_Y$, and $f$ and $g$ are inverse isomorphisms.  
        \item A map $f\colon X\to Y$ is a cover if and only if its pullback along some cover $g\colon Z\to Y$ is a cover.
        \item A map $f\colon X\to Y$ is an isomorphism if and only if its pullback along some cover $g\colon Z\to Y$ is an isomorphism.  
        \item Isomorphisms are covers.
        \item Products of covers are covers.
    \end{enumerate}
\end{lemma}
\begin{proof}
    For the first statement, by assumption, $1_Y\circ f=f=f\circ 1_X=(fg)\circ f$. By Axiom~\ref{a:subcan}, covers are epimorphisms. The statement follows.

    For the second statement, the ``only if'' is just Axiom~\ref{a:basechange}, while the ``if'' follows from Axiom~\ref{a:basechange} and our assumption that covers form a subcategory.

    For the third statement, note that the pullback exists by Axiom~\ref{a:basechange}. The ``only if'' is trivial.  For the ``if'', suppose that $g^\ast f\colon Z\times_Y X\to Z$ is an isomorphism. We have an commuting diagram
    \[
        \xymatrix{
        (Z\times_Y X)\times_X (Z\times_Y X) \ar[d]^\cong_{g^\ast f\times g^\ast f} \ar@<-.5ex>[r] \ar@<.5ex>[r] & Z\times_Y X\ar[d]^\cong_{g^\ast f}\ar[r]^-{f^\ast g} & X\ar[d]^f\\
        Z\times_Y Z  \ar@<-.5ex>[r] \ar@<.5ex>[r] & Z \ar[r]^g \ar[r] & Y
        }
    \]
    By Axiom~\ref{a:basechange}, the map $f^\ast g$ is a cover. By Axiom~\ref{a:subcan}, both commutative forks are coequalizers. Our assumption that $g^\ast f$ is an isomorphism combined with the universal property of coequalizers implies that $f$ is an isomorphism as claimed.

    For the fourth statement, let $f\colon X\to^\cong Y$ be an isomorphism. The square
    \[
        \xymatrix{
            X \ar[r]^f \ar[d] & Y \ar[d]\\
            \ast \ar[r] & \ast
        }
    \]
    is a pullback, and the map $\ast$ to $\ast$ is a cover by Axiom~\ref{a:term}. We conclude $f$ is a cover by Axiom~\ref{a:basechange}.

    For the fifth statement, let $f\colon X\to Y$ and $g\colon W\to Z$ be covers. The map
    \[
        f\times g\colon X\times W\to Y\times Z
    \]
    factors as
    \[
        X\times W\to^{f\times 1} Y\times W\to^{1\times g}Y\times Z.
    \]
    Because covers compose, it suffices to show each of these maps is a cover. The first map sits in a pullback square
    \[
        \xymatrix{
            X\times W \ar[d]^{f\times 1} \ar[r] & X \ar[d]_f\\
            Y\times W \ar[r] & Y
        }
    \]
    By Axiom~\ref{a:basechange} and the previous statement, we conclude $f\times 1$ is a cover. The analogous argument shows the same for $1\times g$.
\end{proof}

One perspective on categories with covers is that they provide a setting for working with ``smooth'' objects in geometric contexts, with the covers being ``smooth'' morphisms.  Relevant examples include:
\begin{enumerate}
    \item the category $\Set$, with covers being surjections,
    \item the category $\sm$ of smooth Banach manifolds with surjective submersions as covers,
    \item the category of complex manifolds with surjective (holomorphic) submersions as covers.
\end{enumerate}

Let $\PSh(\C)$ denote the category of presheaves on $\C$. As in \cite{W,W2}, the axioms ensure that the covers define a (singleton) Grothendieck topology on $\C$. Let 
\[
    \Sh(\C)\subset\PSh(\C)
\]
denote the category of sheaves with respect to this topology.  As noted in {\em loc. cit.}, our assumption that covers are effective epimorphisms implies that the Yoneda embedding takes values in sheaves, i.e. we have canonical fully faithful embeddings
\[
    \C\subset\Sh(\C)\subset\PSh(\C)
\]
We identify an object in $\C$ with its functor of points (the sheaf it represents) and use the language of sheaves and generalized elements in our study of objects in $\C$.

\begin{example}
    Let $\C$ be a category with covers.  Recall (e.g. \cite[p. 15]{DHI}) that a map $f\colon F\to G$ of presheaves on $\C$ is a {\em generalized cover} or {\em local epimorphism} if for all sections $U\to G$, there exists a cover $V\to U$ and a lift of $V\to G$ through $f$
        \[
            \xymatrix{
            && F \ar[d]_f\\
            V \ar[r] \ar@{..>}[urr]^\exists & U\ar[r] & G
            }
        \]
    The category $\Sh(\C)$ with local epimorphisms as covers satisfies the axioms of a category with covers.
\end{example}

\begin{remark}
    If $\Sh(\C)$ has enough points, then a map of presheaves is a local epimorphism if and only if it is surjective on all stalks. In general, being a local epimorphism of representables is a strictly weaker condition than being a cover.  For example, every retract is a local epimorphism, but retracts are not generally covers.  Indeed, let $\C=\sm$ with surjective submersions as covers, let $X$ be a positive dimensional connected manifold, let $x\in X$ be a point. Then the map $X\sqcup\{x\}\to X$ is a retract, and thus a local epimorphism of sheaves, but it is not a surjective submersion.
\end{remark}

\subsection{Simplicial objects and smooth $\infty$-groupoids}\label{s:sobj}
Let $\Delta$ denote the ordinal category, i.e. the category of finite, non-empty, linearly ordered sets and nondecreasing maps between them. Let $\sC$ denote the category of simplicial objects in $\C$, i.e. 
\[
    \sC:=\Fun(\Delta^{\op},\C).
\]
The inclusion of constant diagrams gives a fully faithful embedding 
\[
    \C\subset\sC.
\]
Similarly, observe that the assignment $\ast_{\Set}\mapsto \ast_{\C}$ extends to a fully faithful, cocontinuous embedding
\[
    \Set\subset\Sh(\C). 
\]
This extends to a fully faithful embedding
\[
    \sSet\subset\sSh(\C).
\]
We will treat simplicial sets as simplicial sheaves on $\C$ without further comment. 

For $n\ge 0$, let $\Delta_{\le n}\subset\Delta$ be the full subcategory on ordinals of cardinality at most $n$. Let $\sC_{\le n}$ denote the category of $n$-truncated simplicial objects in $\C$, i.e.
\[
    \sC_{\le n}:=\Fun(\Delta^{\op}_{\le n},\C).
\]
Restriction along the inclusion $\Delta_{\le n}\subset\Delta$ gives a functor
\begin{align*}
    \sC&\to\sC_{\le n}\\
    X&\mapsto X_{\le n}.
\end{align*}
For $n=0$, this makes $\C$ into a retract of $\sC$.

We write $\Delta^n\in\sSet$ for the standard $n$-simplex, $\partial\Delta^n\subset\Delta^n$ for its boundary, and $\Lambda^n_i\subset\partial\Delta^n$ for its $i^{th}$ horn.  We write
\begin{align*}
    s^i&\colon \Delta^{n+1}\to \Delta^n\\
    d^i&\colon \Delta^{n-1}\to \Delta^n
\end{align*}
for the codegeneracy and coface maps. For a simplicial object $X\in\sC$, we write
\begin{align*}
    s_i&\colon X_n\to X_{n+1}\\
    d_i&\colon X_n\to X_{n-1}
\end{align*}
for the degeneracy and face maps.

Let $X\in\sC$ and $K$ a finite simplicial set.  We write $\hom(K,X)$ to denote the internal hom sheaf, defined in the usual manner.\footnote{If $\C$ has finite coproducts and coproducts of covers are covers, as in the examples listed above, then this is just the sheaf which assigns to $U\in\C$, the set of simplicial maps $\hom(U\times K,X)$.  }  Note that the construction $(K,X)\mapsto \hom(K,X)$ is covariant in $X$ and contravariant in $K$. This is a frequent and familiar construction.
\begin{example} 
    Let $K=\Delta^n$ and let $X\in\sC$.  Then
    \[
        \hom(\Delta^n,X)=X_n 
    \]
    and $s_i=(s^i)^\ast\colon X_n\to X_{n+1}$ and similarly $d_i=(d^i)^\ast$.
\end{example}

Recall that the {\em spine} of an $n$-simplex $\Delta^n$ is the union of its edges $\Delta^n|_{\{i<i+1\}}$ for $i=0,\ldots,n-1$.  Given a simplicial object $X$, denote by 
\begin{align*}
        \spine\colon X_n&\to X_1^{\times n}\intertext{the map which sends a $n$-simplex to its spine (for any $n$)}
        x&\mapsto (x|_{\{0<1\}},\ldots,x|_{\{n-1<n\}})
\end{align*}

We write
\begin{align*}
    \Lambda^n_iX&:=\hom(\Lambda^n_i,X),\\
    M_nX&:=\hom(\partial\Delta^n,X),\\
    \Csk_{k,n}X&:=\hom(\sk_k\Delta^n,X).
\end{align*}
Note that $\sk_{n-1}\Delta^n=\partial\Delta^n$, and the inclusions $\sk_k\Delta^n\subset\Lambda^n_i\subset\partial\Delta^n\subset\Delta^n$ (for $k<n-1$) induce canonical maps
\[
    \xymatrix{
        X_n\ar[d]_{\mu_n} \ar[dr]_{\lambda^n_i} \ar[drr]^{\sigma_k^n} & & \\
        M_nX \ar[r] & \Lambda^n_i X \ar[r] & \Csk_{k,n} X
    }
\]
For $0<k$, the inclusion $\spine \Delta^n\subset \sk_k\Delta^n$ induces a canonical map
\[
    \xymatrix{
        X_n\ar[d]_{\sigma_k^n} \ar[dr]^{\spine} \\
        \Csk_{k,n} X \ar[r] & X_1^{\times n}
    }
\]
We will also need relative versions of several of the above.  Let $f\colon X\to Y$ be a map of simplicial objects in $\C$.  Following \cite{W}, we write
\begin{align*}
    \Lambda^n_i(f)&:=\Lambda^n_i X\times_{\Lambda^n_i Y} Y_n,\\
    M_n(f)&:=M_nX\times_{M_nY}Y_n,\intertext{and also}
    \Csk_{k,n}(f)&:=\Csk_{k,n} X\times_{\Csk_{k,n} Y}Y_n.
\end{align*}
We denote the natural maps induced from $\sk_k\Delta^n\subset \Lambda^n_i\subset \partial\Delta^n\subset \Delta^n$ (for $k<n-1$) as follows
\[
    \xymatrix{
        X_n \ar[d]_{\mu_n(f)} \ar[dr]_{\lambda^n_i(f)} \ar[drr]^{\sigma_k^n(f)} & &\\
         M_n(f) \ar[r] & \Lambda^n_i(f) \ar[r] & \Csk_{k,n}(f)
    }
\]
Note that for a map $f\colon X\to Y$ in $\sC$ and a fixed $k$, as $n$ varies, the objects $\Csk_{k,n}(f)$ naturally define a simplicial sheaf $\Csk_k(f)$. We refer to this as the {\em relative $n$-coskeleton} of $f$.  Similarly, for an object $X\in\sC$, we obtain the {\em $n$-coskeleton} $\Csk_k X$ with $n$-simplices given by $\Csk_{k,n} X$. From the definitions above, we have
\begin{align*}
    \Csk_n(f)&\cong \Csk_n X\times_{\Csk_n Y}Y.
\end{align*}
We say that a map $f\colon X\to Y$ is {\em $n$-coskeletal} if the natural map $X\to \Csk_n(f)$ is an isomorphism. Note that if $f$ is $n$-coskeletal, then for any $g\colon Z\to Y$, we have a natural isomorphism
\[
    \hom_{\sC/Y}(Z,X)\cong \hom_{\sC_{\le n}/Y_{\le n}}(Z_{\le n},X_{\le n}).
\]

\begin{definition} \label{def:basics}
    Let $f\colon X\to Y$ be a map of simplicial objects in $\C$. We say that $f$ is:
    \begin{enumerate}
        \item a {\em hypercover} if $\mu_n(f)$ is a cover for all $n\ge 0$;
        \item an {\em $n$-hypercover} if $f$ is a hypercover and an $\mu_k(f)$ is an isomorphism for all $k\ge n$; 
        \item a {\em Kan fibration} if $\lambda^n_i(f)$ is a cover for all $n,i$;
        \item a {\em covering Kan fibration} if $f$ is a Kan fibration and $f_0\colon X_0\to Y_0$ is a cover.
        \item an {\em $n$-stack} if $f$ is a Kan fibration and $\lambda^k_i(f)$ is an isomorphism, for all $k>n$.
        \item a {\em covering $n$-stack} if $f$ is an $n$-stack and $f_0$ is a cover.
    \end{enumerate} 
    We say a simplicial object $X\in\sC$ is {\em Kan} if $X\to \Delta^0$ is a Kan fibration.
\end{definition}

\begin{remark}
    We will frequently want to distinguish between a map $f\colon X\to Y$ in $\sC$ being e.g. a hypercover in $\sC$, and it satisfying the corresponding property in the category $\sSh(\C)$ with covers the local epimorphisms.  As a global convention, we will refer to the former as a hypercover and the latter as a ``hypercover of (representable) simplicial sheaves''.  Similarly for Kan fibrations, etc.  We will also refer to hypercovers of simplicial sheaves as ``local trivial fibrations'', Kan fibrations of simplicial sheaves as ``local Kan fibrations'', $n$-stacks of simplicial sheaves as ``local $n$-fibrations'', etc. 
\end{remark} 

\begin{definition}\label{d:inftygpd}
    A {\em smooth $\infty$-groupoid in $\C$} is a Kan simplicial object, i.e. an $X\in \sC$ for which the maps 
    \[
        \lambda^n_i\colon X\to \Lambda^n_i X
    \]
    are covers for all $n,i$.  A {\em smooth $n$-groupoid in $\C$} is a Kan simplicial object for which the maps $\lambda^k_i$ are isomorphisms for $k>n$. A {\em smooth $\infty$-group} is a reduced Kan simplicial object, i.e. a smooth $\infty$-groupoid $X$ with $X_0=\ast$.  A {\em smooth $n$-group} is a smooth $n$-groupoid $X$ with $X_0=\ast$.
\end{definition}

The following lemmas are standard and will be used repeatedly below. 
\begin{lemma}\label{l:stackcosk}
    Let $f\colon X\to Y$ be an $n$-stack. Then $f$ is $(n+1)$-coskeletal.
\end{lemma}
\begin{proof}
    This follows immediately from \cite[Definition 4.4, Proposition 4.5]{W}.
\end{proof}

\begin{lemma}\label{l:covkan}
	Let $f\colon X\to Y$ be a covering Kan fibration. Then for all $n>0$, the map 
	\[
		f_n\colon X_n\to Y_n
	\]
    is a cover in $\C$.
\end{lemma}
\begin{proof}
    This is a special case of \cite[Lemma 2.8]{W2}.
\end{proof}

We will need a mild strengthening of a standard result (see \cite[Lemma 2.10]{Hen}, \cite[Lemmas 2.15, 2.16]{W}, \cite[Lemma 3.9]{BG}).  Recall that an inclusion of simplicial sets $S\into T$ is an {\em expansion} if there exists a filtration
\begin{align*}
    S=S_0&\into\cdots\into S_N=T\intertext{where for all $0<\ell\le N$}
    S_\ell&\cong S_{\ell-1}\cup_{\Lambda^{n_\ell}}\Delta^{n_\ell}.
\end{align*}
A finite simplicial set is {\em collapsible} if the inclusion of some vertex is an expansion. Given an inclusion of simplicial set $S\into T$ and a map $f\colon X\to Y$ in $\sC$, we write
\[
    (S\into T)(f):=\hom(S,X)\times_{\hom(S,Y)}\hom(T,Y).
\]
Note that this mildly differs from our notation for $\Lambda^k_i(f)$ and $M_n(f)$.

\begin{lemma}\label{l:expand}
    Let $S\into T$ be an inclusion of finite simplicial sets. Let $S$ be collapsible, with filtration
    \[
        \Delta^0=S_0\into\cdots\into S_N=S
    \]
    with
    \[
        S_{\ell+1}\cong S_\ell\cup_{\Lambda^{n_\ell}_{i_\ell}}\Delta^{n_\ell}
    \]
    for all $0\le \ell<N$. Let
    \[
        \xymatrix{
        X\ar[rr]^f \ar[dr]_g && Y \ar[dl]^h\\
        & Z
        }
    \]
    be a diagram in $\sC$.  Suppose that $f_0\colon X_0\to Y_0$ is a cover, that $\lambda^k_i(f)$ is a cover for all $k\le\dim S$. Suppose that $(S_\ell\into T)(h)$ exists in $\C$ for all $0\le \ell\le N$. 
    
    Then for all $\ell$, $(S_\ell\into T)(g)$ exists in $\C$ and the map
    \[
        f_\ast\colon (S_\ell\into T)(g)\to (S_\ell\into T)(h)
    \]
    is a cover.
\end{lemma}
\begin{proof}
    We prove the lemma by inducting along the length of the filtration. For the base case, $S=\Delta^0$. We have a pullback square
    \[
        \xymatrix{
            X_0\times_{Z_0}\hom(T,Z) \ar[r]^{f_\ast} \ar[d] & Y_0\times_{Z_0}\hom(T,Z) \ar[d] \\
            X_0 \ar[r]^{f_0} & Y_0
        }
    \]
    Because covers pull back (Axiom~\ref{a:basechange}), our assumption on $f_0$ implies that the pullback exists in $\C$ and the top horizontal map is a cover.

    Now assume we have shown that the map
    \[
        f_\ast\colon (S_\ell\into T)(g)\to (S_\ell\into T)(h)
    \]
    is a cover in $\C$. We have a pullback square
    \[
        \xymatrix{
            (S_\ell\into T)(g)\times_{(S_\ell\into T)(h)}(S_{\ell+1}\into T)(h) \ar[r] \ar[d] & (S_{\ell+1}\into T)(h) \ar[d] \\
            (S_\ell\into T)(g) \ar[r]^{f_\ast} & (S_\ell\into T)(h)
        }
    \]
    By our inductive hypothesis and Axiom~\ref{a:basechange}, the top horizontal map is a cover in $\C$. Next, consider the pullback square
    \[
        \xymatrix{
           (S_{\ell+1}\into T)(g) \ar[r] \ar[d] & (S_\ell\into T)(g)\times_{(S_\ell\into T)(h)}(S_{\ell+1}\into T)(h)  \ar[d] \\
            X_{n_\ell} \ar[r]^{\lambda^{n_\ell}_{i_\ell}(f)} & \Lambda^{n_\ell}_{i_\ell}(f)
        }
    \]
    The bottom map is a cover by assumption. This implies that the top horizontal map is a cover. Because covers compose, we see that the map
    \[
        (S_{\ell+1}\into T)(g)\to (S_{\ell+1}\into T)(h)
    \]
    is a cover as claimed.  This completes the inductive step and thus the proof.
\end{proof}

\begin{lemma}\label{l:kanngpd}
    Suppose we have a diagram in $\sC$
    \[
        \xymatrix{
        X\ar[rr]^f \ar[dr]_g && Y \ar[dl]^h \\
        & Z
        }
    \]
    where $g$ is an $n$-stack.
    \begin{enumerate}
        \item If $f$ is a Kan fibration, then $f$ is an $n$-stack.
        \item If $f$ is a covering Kan fibration and $d_i\colon Z_k\to Z_0$ is a cover for all $k$ and $i$, then $h$ is also an $n$-stack.
        \item If $g$ and $h$ are $n$-stacks, the map $f$ is an isomorphism if and only if it induces an isomorphism on $n$-skeleta.
    \end{enumerate}
\end{lemma}
\begin{proof}
    Observe that the for all $k$ and $i$, we have the following commuting diagram
    \begin{equation}\label{e:filler}
        \xymatrix{
            X_k \ar[r]^{\lambda^k_i(f)} \ar[dr]_{\lambda^k_i(g)} & \Lambda^k_i(f) \ar[d] \ar[r] & Y_k \ar[d]^{\lambda^k_i(h)} \\
            & \Lambda^k_i(g) \ar[r]^{f_\ast} & \Lambda^k_i(h)
        }
    \end{equation}
    where the square is a pullback.  

    Suppose $f$ is a Kan fibration. By Axiom~\ref{a:rightcancellation}, the left vertical map in the diagram above is a cover. For $k>n$, the map $\lambda^k_i(g)$ is an isomorphism. The composition 
    \[
        \Lambda^k_i(f)\to\Lambda^k_i(g)\to^{\lambda^k_i(g)^{-1}} X_k
    \]
    is thus a cover and a left inverse to $\lambda^k_i(f)$.  We conclude by Lemma~\ref{l:covepi} that $f$ is an $n$-stack. 

    Now assume that $f$ is a covering Kan fibration, and that $d_0,d_1\colon Z_1\to Z_0$ are covers. We show by a simple induction that $\lambda^k_i(h)$ and $f_\ast$ are covers for all $k$, $i$. Note that for all $k,i$, the simplicial set $\Lambda^k_i$ is collapsible. Observe also that the map $f_k\colon X_k\to Y_k$ is a cover for all $k$ by Lemma~\ref{l:covkan}, and $\lambda^k_i(g)$ is a cover for all $k$ and $i$ by assumption.
    
    For the base case, $k=1$, we have a pullback square
    \[
        \xymatrix{
        \Lambda^1_i(h) \ar[r] \ar[d] & Z_1 \ar[d]^{d_{1-i}}\\
        Y_0 \ar[r] & Z_0
        }
    \]
    By our assumption on $Z$ and the fact that covers pull back (Axiom~\ref{a:basechange}), the square exists in $\C$. Lemma~\ref{l:expand} implies that $f_\ast$ is a cover. By Axiom~\ref{a:rightcancellation} (and the fact that covers form a subcategory), we see that $\lambda^1_i(h)$ is a cover for all $i$. 
    
    Now suppose we have shown that $f_\ast$ and $\lambda^m_i(h)$ are covers for all $m<k$. By \cite[Lemma 2.15]{W2}, our assumption on $d_i\colon Z_k\to Z_0$ combines with our inductive hypothesis to imply that $\Lambda^k_i(h)$ exists in $\C$. By Lemma~\ref{l:expand}, the map 
    \[
        f_\ast\colon \Lambda^k_i(g)\to \Lambda^k_i(h)
    \]
    is a cover. We conclude by the argument above that $\lambda^k_i(h)$ is a cover as well, thus completing the inductive step and showing that $h$ is a Kan fibration. To see that $h$ is in fact an $n$-stack, observe that $g$ and $f$ are $n$-stacks. By the 2-of-3 property of isomorphisms, for $k>n$, the pullback of $\lambda^k_i(h)$ along the cover $f_\ast\colon \Lambda^k_i(g)\to\Lambda^k_i(h)$ as in the diagram~\eqref{e:filler} is an isomorphism. By Lemma~\ref{l:covepi}, we see that $\lambda^k_i(h)$ is an isomorphism for $k>n$ as claimed.  
    
    For the last statement, the ``if'' is trivial.  For the ``only if'', assume that $f$ induces an isomorphism on $n$-skeleta. By Lemma~\ref{l:stackcosk}, it suffices to show $f$ induces an isomorphism on $(n+1)$-skeleta. Taking $k=n+1$ in the square~\eqref{e:filler} above, we see that $f_{n+1}$ is the composite of the top horizontal maps. The bottom horizontal map is an isomorphism by our assumption on $f$. The last statement now follows.
\end{proof}

\begin{remark}
    Note that $n$-stacks are preserved under composition \cite[Theorem 2.17]{W}. Together with this, we can view the above lemma as establishing a ``2-of-3'' property for $n$-stacks.
\end{remark}

Following Henriques \cite[Section 3]{Hen} (which itself follows van Osdol \cite{Van}), for $X$ a smooth $\infty$-groupoid in $\C$ and $x\in X_0$, we have group sheaves $\pi_i(X,x)\in\Sh(\C)$; when $X$ is a smooth $\infty$-group, we will just write $\pi_i X$. These have the expected properties, including functoriality, invariance under hypercovers \cite[Theorem 7.1]{RZ} and the usual long exact sequence associated to a fibration \cite[Propositon 3.4]{Hen}. 

\begin{definition}
    A map $f\colon X\to Y$ of smooth $\infty$-groupoids in $\C$ is a {\em local weak equivalence} if it induces an isomorphism of sheaves $f_\ast\colon \pi_0X\to^\cong \pi_0 Y$ and for all $x\in X$ and $i\ge 1$, it induces isomorphisms $f_\ast\colon\pi_i(X,x)\to^\cong \pi_i(Y,f(x))$.
\end{definition}

\begin{remark}
As shown by the first author's joint work with Zhu \cite{RZ}, if the topos $\Sh(\C)$ satisfies an additional, yet geometrically natural, hypothesis \cite[Def.\ 6.5]{RZ}, then the category of smooth $\infty$-groupoids admits an internal homotopy theory as an ``incomplete category of fibrant objects'' or ``iCFO'', in which the Kan fibrations serve as fibrations and the hypercovers as trivial fibrations. This framework makes explicit the homotopy theory underlying the arguments of Henriques \cite{Hen} and the second author \cite{W}, analogous to how Quillen's model categories made explicit the structure underlying classical homological algebra or the homotopy theory of simplicial sets. In what follows, we tacitly assume that the category of $\infty$-groupoids in $\C$ is equipped with this iCFO structure.
\end{remark}

\begin{lemma}\label{l:rightproper}
    The iCFO of smooth $\infty$-groupoids in $\C$ is right proper: given a Kan fibration $f\colon X\to Y$ of smooth $\infty$-groupoids in $\C$, and a local weak equivalence $g\colon Z\to Y$ such that the pullback $X\times_Y Z$ exists in $\sC$, then $p\colon X\times_Y Z\to X$ is a local weak equivalence.
\end{lemma}
\begin{proof}
    By \cite[Theorem 2.17(4)]{W}, the map $X\times_Y Z\to Z$ is a Kan fibration. The lemma now follows mutatis mutandis from the standard argument in simplicial sets.  To wit, by the functoriality of the long exact sequence in homotopy groups \cite[Proposition 3.4]{Hen} (which follows \cite[Section II]{Van}) and the Five Lemma, we see that for all basepoints $w\in X\times_YZ$, the map induces isomomorphisms $\pi_i(X\times_YZ,w)\cong \pi_i(X,p(w))$ for $i>0$. The classical argument via path-lifting shows that $p$ induces an isomorphism on $\pi_0$ as well.
\end{proof}

\subsection{Nerves, Eilenberg-Mac Lane spaces, and homotopy quotients} \label{sec:nerves}
We close this section by reviewing some classical constructions on groups and simplicial groups. We first recall two isomorphic presentations of the nerve $NG$ of a group object $G$ in $\C$. The most common presentation involves ``homogeneous'' coordinates and is given by
\begin{align*}
    NG_n&:=G^n\\
    d_0(g_1,\ldots,g_n)&:=(g_2,\ldots,g_n)\\
    d_i(g_1,\ldots,g_n)&:=(g_1,\ldots,g_ig_{i+1},\ldots,g_n)\tag{for $0<i<n$}\\
    d_n(g_1,\ldots,g_n)&:=(g_1,\ldots,g_{n-1})\\
    s_i(g_1,\ldots,g_n)&:=(g_1,\ldots,g_i,e,g_{i+1},\ldots,g_n)
\end{align*}
Less frequently in the literature, but appearing naturally in the context of integration, is the following.
\begin{definition}[Inhomogenous Coordinates on $NG$]
    Let $G$ be a group object in $\C$.  Define $NG^\iota$ to be the simplicial object with
    \begin{align*}
        NG^\iota_n&:=G^n\\
        d_0(g_1,\ldots,g_n)&:=(g_1^{-1}g_2,\ldots,g_1^{-1}g_n)\\
        d_i(g_1,\ldots,g_n)&:=(g_1,\ldots,g_{i-1},g_{i+1},\ldots,g_n)\tag{for $0<i$}\\
        s_0(g_1,\ldots,g_n)&:=(e,g_1,\ldots,g_n)\\
        s_i(g_1,\ldots,g_n)&:=(g_1,\ldots,g_i,g_i,g_{i+1},\ldots,g_n)
    \end{align*}
    We refer to $NG^\iota$ as the nerve of $G$ with {\em inhomogeneous coordinates}.
\end{definition}

An exercise in the formulas shows the following.
\begin{lemma}\label{l:nervecoord}
    Let $G$ be a group object in $\C$.  For $n\ge 0$, the assignments
    \begin{align*}
        (g_1,\ldots,g_n)&\mapsto (g_1,g_1g_2,\ldots,g_1\cdots g_n)\intertext{define a simplicial isomorphism}
        NG&\to^\cong NG^\iota.
    \end{align*}
\end{lemma}

\begin{remark}
    In light of the above lemma, we can omit the superscript $\iota$ and instead just specify which coordinates we are working with.
\end{remark}

If $A$ is an abelian group object in $\C$, we also have the Eilenberg-Mac Lane objects $K(A,n)\in\sC$ for $n\ge 0$ \cite{EM} (for a more recent reference, see \cite[Ch. III.2]{GJ}). Just as for nerves of groups, we have homogeneous and inhomogeneous coordinates on them, which we now review.

The most common presentation of $K(A,n)$ arises via the Dold--Kan correspondence (cf. \cite[III.2]{GJ}), and gives a smooth $n$-group with
\begin{align*}
    K(A,n)_{<n}&=\Delta^0_{<n}\\
    K(A,n)_n&\cong A\\
    K(A,n)_{n+1}&\cong A^{n+1}\intertext{and nontrivial face and degeneracy maps given by}
    s_ia&=(0,\ldots,\overbrace{a}^i,\ldots,0)\\
    d_0(a_0,\ldots,a_n)&=a_0\\
    d_{n+1}(a_0,\ldots,a_n)&=a_n\intertext{and for $i<n+1$}
    d_i(a_0,\ldots,a_n)&=a_{i-1}+a_i.
\end{align*}
Note that, as a smooth $n$-group, $K(A,n)$ is $(n+1)$-coskeletal (Lemma~\ref{l:stackcosk}), so the formulas above determine $K(A,n)\in\sC$ up to isomorphism. We also have ``inhomogeneous'' coordinates, which often appear in geometric contexts.
\begin{definition}\label{d:Ainhomog}[Inhomogenous Coordinates on $K(A,n)$]
    Let $A$ be an abelian group object in $\C$. Fix $n\ge 0$.  Define $K(A,n)^\iota$ to be the $(n+1)$-coskeletal simplicial abelian group object with
    \begin{align*}
        K(A,n)^\iota_{<n}&=\Delta^0_{<n}\\
        K(A,n)^\iota_n&=A\\
        K(A,n)^\iota_{n+1}&=A^{n+1}\intertext{and nontrivial face and degeneracy maps given by}
        s_ia&=(0,\ldots,0,\overbrace{a,a}^{i,i+1},0,\ldots,0)\\
        d_{n+1}(a_0,\ldots,a_n)&=(-1)^n\sum_{i=0}^n (-1)^ia_i\intertext{and for $i<n+1$}
        d_i(a_0,\ldots,a_n)&=a_i.
    \end{align*}
\end{definition}
Using the $(n+1)$-coskeletality, an exercise in the formulas shows the following.
\begin{lemma}
    Let $A$ be an abelian group object in $\C$.  For $n\ge 0$, the assignments
    \begin{align*}
        \varphi_n(a)&:=a\\
        \varphi_{n+1}(a_0,\ldots,a_n)&:=(a_0,a_0+a_1,a_1+a_2,\ldots,a_{n-1}+a_n)\intertext{define a simplicial isomorphism}
        \varphi\colon K(A,n)&\to^\cong K(A,n)^\iota.
    \end{align*}
\end{lemma}

\begin{remark}
    Just as for nerves of groups the above lemma justifies us in dropping the superscript $\iota$ from $K(A,n)$ and instead just specify which coordinates we are working with.
\end{remark}

We will need one more perspective on $K(A,n)$, this time viewing it as a ``nerve'' of the simplicial group object $K(A,n-1)$.\footnote{Because $A$ is an {\em abelian} group object, $K(A,n-1)$ is a simplicial group object with the group operations given coordinate-wise.}  Recall that for a simplicial group $G$ in $\sC$,  its {\em nerve} $\LW G\in \sC$ (first introduced by \cite{EM}) can be defined as follows:\footnote{Our conventions here differ from those in \cite[Section 6]{W} by an orientation reversal: the ordinal category $\Delta$ has the ``order-reversing'' involution which sends $\{0<\ldots<n\}$ to $\{n<\ldots<0\}$. Pulling back along this involution converts the present convention into that of \cite{W} and vice versa.}
\begin{align*}
    \LW G_n&=G_{n-1}\times\cdots\times G_0\intertext{with face and degeneracy maps given by}
    s_i(g_{n-1},\ldots,g_0)&=(s_{i-1}g_{n-1},\ldots,s_0g_{n-i},e,g_{n-i-1},\ldots,g_0)\\
    d_0(g_{n-1},\ldots,g_0)&=(g_{n-2},\ldots,g_0)\\
    d_n(g_{n-1},\ldots,g_0)&=(d_{n-1}g_{n-1},\ldots,d_1g_1)\intertext{and, for $0<i<n$}
    d_i(g_{n-1},\ldots,g_0)&=(d_{i-1}g_{n-1},\ldots,d_1g_{n-i+1},d_0g_{n-i}\cdot g_{n-i-1},g_{n-i-2},\ldots,g_0).
\end{align*}
For $G$ a constant simplicial group (i.e. $G\in \C\subset\sC$), this recovers the nerve, i.e.
\[
    \LW G = NG.
\]

\begin{lemma}\label{l:W=K[1]}
    Let $A$ be an abelian group object in $\C$. For $n\ge 0$, let $K(A,n)$ denote the Eilenberg-Mac Lane object with inhomogeneous coordinates. The assignments
    \begin{align*}
        \varphi_{n+1}(a)&:=a\\
        \varphi_{n+2}((a_0,\ldots,a_n),b)&:=(b,a_0+b,a_1,\ldots,a_n)\intertext{define a simplicial isomorphism}
        \varphi\colon \LW K(A,n)&\to^\cong K(A,n+1).
    \end{align*}
\end{lemma}
\begin{proof}
    A straightforward check verifies that the formulas define a simplicial map. To complete the proof, we must verify that $\LW K(A,n)$ is $(n+1)$-coskeletal.  This follows from \cite[Theorem 6.7]{W} by induction on $n$, noting for the base case that a ``strict Lie $0$-group'' is precisely a group object $G\in C\subset \sC$. 
\end{proof}

\begin{remark}
    Note that the induction on $n$ in the proof shows that $K(A,n)$ is a smooth $n$-group for all $n\ge 1$ (and a smooth $0$-groupoid for $n=0$). 
\end{remark}

We conclude this subsection by recalling homotopy quotients.\footnote{Our treatment continues to follow \cite[Section 6]{W}, but with the orientation reversal on simplices mentioned above. Note also that {\em loc. cit.} writes ``$WG\times_{G}X$'' for our ``$X//G$''.}
\begin{definition}\label{d:hoquot}
    Let $X\in \sC$ and let $G\in sC$ be a group object acting on $X$. Define the {\em homotopy quotient} 
    \[
        X//G\to \LW G
    \]
    as follows:
    \begin{align*}
        (X//G)_n&:=X_n\times \LW G_n,\intertext{and writing $(x,\underline{g})$ for an $n$-simplex of $X//G$, we define degeneracy and face maps as follows:}
        s_i(x,\underline{g})&:=(s_ix,s_i\underline{g})\\
        d_0(x,g_{n-1},\ldots,g_0)&=(g_{n-1}^{-1}\cdot d_0x,d_0\underline{g})\intertext{and for $i>0$}
        d_i(x,\underline{g})&:=(d_ix,d_i\underline{g}).
    \end{align*}
\end{definition}

\begin{remark}\label{r:hoquot}
    By \cite[Theorem 6.7(2)]{W}, if $X$ is a smooth $n$-groupoid, then $X//G\to \LW G$ is a smooth $n$-stack, and thus, by Lemma~\ref{l:stackcosk}, the map is $(n+1)$-coskeletal. In particular, for $A$ an abelian group object in $\C$ and $G\in\C$ acting on $A$ by automorphisms, we see that $K(A,n)//G\to NG$ is a smooth $n$-stack.
\end{remark}

We close with a final example from \cite{W2} which will provide the starting point for our treatment of smooth $k$-invariants in Section~\ref{sec:smooth-kinv}. Recall that, for a cover $X\to Y$ in $\C$, its {\em nerve} $N(X/Y)\to Y$ is the simplicial object over $Y$ with 
\begin{align}\label{e:nerve}
    N(X/Y)_n=\overbrace{X\times_Y\cdots\times_Y X}^{n+1}
\end{align}
and $i^{th}$ face (degeneracy) maps given by omitting (repeating) the $i^{th}$ factor. Now let $f\colon X\to Y$ in $\sC$ be a {\em covering Kan fibration}, i.e. a Kan fibration for which $f_0\colon X_0\to Y_0$ is a cover. In \cite{W2}, the second author introduced a functorial nerve construction
\[
    \LW(f)\colon \LW(X/Y)\to Y.
\]
Concretely,
\[
    \LW(X/Y):=\Tot\circ N(X/Y)
\]
where $\Tot$ denotes the Artin-Mazur totalization \cite{AMTot} and $N$ denotes the level-wise nerve of $f\colon X\to Y$ (using that $f_n$ is a cover for all $n$, by Lemma~\ref{l:covkan}).

In analogy with the classical nerve of a cover, we have the following.
\begin{proposition}\cite[Proposition 1.3]{W2}\label{p:nerve}
    Let $f\colon X\to Y$ be a covering Kan fibration. Then 
    \[
        \LW(f)\colon \LW(X/Y)\to Y
    \]
    is a hypercover. Further, if $f$ is an $n$-stack, then $\LW(f)$ is an $(n+1)$-hypercover.
\end{proposition}
\begin{proof}
    The first statement is \cite[Proposition 1.3]{W}.  For the second, observe that the proof of \cite[Proposition 1.3]{W} proceeds by constructing, for all $k\ge 0$, an isomorphism of maps
    \[
		\xymatrix{
			\bar{W}(X/Y)_k \ar[rr]^{\mu_k(\LW(f))} \ar[d]_\cong && M_k\bar{W}(f) \ar[d]^\cong \\
			X_k\times^{fd_0,f}_{Y_{k-1}} \bar{W}(X/Y)_{k-1} \ar[rr]^{\lambda^k_0(f)\times 1} && \Lambda^k_0 (f)\times_{M_{k-1} Y}^{f\partial d_0,f\partial} \bar{W}(X/Y)_{k-1}.
		}
	\]
    The bottom arrow sits in a pullback diagram
    \[
		\xymatrix{
			X_k\times^{fd_0,f}_{Y_{k-1}} \bar{W}(X/Y)_{k-1} \ar[rr]^{\lambda^k_0(f)\times 1} \ar[d] && \Lambda^k_0 (f)\times_{M_{k-1} Y}^{f\partial d_0,f\partial} \bar{W}(X/Y)_{k-1} \ar[d]\\
            X_k \ar[rr]^{\lambda^k_0(f)} && \Lambda^k_0(f)
		}
	\]
    If $f$ is an $n$-stack, then the bottom arrow is an isomorphism for $k\ge n+1$.  We conclude that $\mu_k(\LW(f))$ is an isomorphism for $k\ge n+1$ as claimed.
\end{proof}

\noindent
Now fix some $B\in\sC$, e.g. $B=NG$ as below.  We work in the category $\sC_{/B}$ and refer to group objects in $\sC_{/B}$ as {\em $B$-groups}. Given a $B$-group $G\to B$, the classical $W$ and $\LW$ constructions can be extended as in \cite{W2} to give a functorial ``universal bundle''
\[
    W(G/B)\to \LW(G/B)
\]
where the functors $W$ and $\LW$ are given by
\begin{align*}
    \LW&:=\Tot\circ N\\
    W&:=\Tot\circ\Dec_1 \circ N. 
\end{align*}
Here $\Dec_1$  is Illusie's ``level-wise'' initial d\'ecalage \cite{Illusie} corresponding to the endofunctor on $\Delta$ given by
\[
    [n]\mapsto [0]\sqcup[n].
\]
The natural transformation $W\to \LW$ is obtained by applying $\Tot$ component-wise to the natural transformation $\Dec_1\to\Id$.  We record a lemma for later use.

\begin{lemma}\label{l:relw=k[1]}
    Let $G$ be a group object in $\C$, $A$ an abelian group object in $\C$ and let $G$ act on $A$ by automorphisms. Then there exists commuting square over $NG$, a natural in $G$ and $A$,
    \[
        \xymatrix{
        W((K(A,n)//G)/NG) \ar[d]  \ar[r]^-\cong & WK(A,n)//G\ar[d] \ar[r]^-\cong & WK(A,n)//G \ar[d]   \\
        \LW((K(A,n)//G)/NG)\ar[r]^-\cong  &  \LW K(A,n)//G \ar[r]^-\cong & K(A,n+1)//G
        },
    \]
    in which the lower right horizontal map is the homotopy quotient of the map in Lemma~\ref{l:W=K[1]}.
\end{lemma}
\begin{proof}
    The existence of the right square follows by inspection from the formula for Lemma~\ref{l:W=K[1]}. The existence bottom left horizontal isomorphism follows immediately from unpacking the definition, as in \cite[Example 3.6]{W2}. The existence of the upper left horizontal isomorphism, fitting into the left square, follows similarly.
\end{proof}

\section{Postnikov theory for smooth $\infty$-groups}\label{s:post}
In this section, we develop a theory of Postnikov sections and $k$-invariants for smooth $\infty$-groups in a category with covers.  Classically, Postnikov theory classifies a reduced Kan complex $X$ in terms of sequences of homotopy groups and $k$-invariants, i.e. maps from a given stage of the Postnikov tower to a twisted Eilenberg-Mac Lane space as in the triangle
\begin{equation*}
    \xymatrix{
    X^{(n-1)} \ar[dr] \ar[rr]^-{k_n} &&  K(\pi_nX,n+1)//\pi_1X \ar[dl]\\
    & N\pi_1 X
    }
\end{equation*}
We aim to develop an analogous treatment for smooth $\infty$-groups in a category with covers $\C$. 

\subsection{Postnikov sections}\label{s:post1}
We begin by recalling two models $\tau_{<n}$ and $\tau_{\le n}$ of relative Postnikov sections, due to Moore \cite{Moore} and Duskin \cite{Duskin}. While each gives a model of a Postnikov tower, a key insight is that, used in tandem, they give a functorial tower
\[
    X\to\cdots\to\tau_{\le n}X\to\tau_{<n}X\to \tau_{\le n-1}X\to\cdots\to\cdots\tau_{\le 0}X=\pi_0X
\]
where the maps $\tau_{\le n}X\to \tau_{<n}X$ are minimal Kan fibrations and the maps $\tau_{<n}X\to\tau_{\le n-1}X$ are trivial fibrations (see \cite[Section 5]{W2} for a recent discussion for simplicial sets).  Our goal here is to show that the core features carry over into $\sC$ for any category with covers.

\begin{definition}\label{d:moore}
    Let $f\colon X\to Y$ be a Kan fibration in $\sC$. Let $n>0$.  Define the {\em relative Moore $n$-truncation} $\tau_{<n}(X,f)$ to be the simplicial sheaf where
    \begin{enumerate}
        \item $\tau_{<n}(X,f)_{\le n-1}=X_{\le n-1}$, and for $k\ge n$
        \item $\tau_{<n}(X,f)_k=\Image(\sigma_{n-1}^k(f)\colon X_k\to \Csk_{n-1,k}(f))$,
    \end{enumerate}
    Note that the image is taken sheaf-theoretically.
\end{definition}
This construction is functorial and fits into a natural commuting diagram
\[
    \xymatrix{
        X \ar[r]^-{p} \ar[dr]_f & \tau_{<n}(X,f) \ar[d]^{\tau_{<n}f} \\
        & Y
        }
\]
When $Y=\ast$, we write $\tau_{<n}X$ and refer to this as the {\em Moore $n$-truncation}.

The relative Moore $n$-truncation gives a model of the $(n-1)^{st}$ relative Postnikov section of the Kan fibration $f$; see \cite{Moore,May,GJ} for classical references, and Henriques \cite{Hen} for a treatment in the present context. We now recall a second model, due to Duskin \cite{Duskin}, and studied in the present context in \cite{Hen,W}. 

In \cite[Definition 3.5]{W}, the second author introduced a {\em relative higher morphism space} $P^{\ge k}(f)$ associated to a Kan fibration $f\colon X\to Y$.  The object $P^{\ge k}(f)$ is defined using the simplicial join operation, and efficiently models the space of $k$-dimensional paths in $X$ which are ``vertical'' with respect to $f$.  In particular, \cite[Theorem 3.6]{W} established that $P^{\ge k}(f)$ is a smooth $\infty$-groupoid when $f$ is Kan, and if $f$ is an $n$-stack, then $P^{\ge k}(f)$ is a smooth $(n-k)$-groupoid (as expected).  Similarly, there exists an augmentation
\[
    P^{\ge k}(f)\to M_k(f),
\]
and if $f$ is a hypercover, this augmentation is a also hypercover \cite[Theorem 3.7]{W} (i.e. the space of such morphisms is contractible rel boundary).  Note that for $Y=\ast$ and $k=1$, this construction was introduced earlier by Duskin \cite{Duskin2} and used to effect by Lurie \cite{Lurie}.\footnote{For a simplicial set $X$ and $x,y\in X_0$, Lurie uses the notation $\hom^L_X(x,y)$ for the fiber at $(x,y)$ of the augmentation $P^{\ge 1}X\to M_1X$.} 

In \cite[Section 4]{W}, the second author defined the relative Duskin truncation via this relative higher morphism space. Note that the presentation in \cite[Section 4]{W} extends the treatment in Glenn \cite[Example 5]{Glenn} (which treats the case $\C=\Set$ and $Y=\ast$). An exercise in the definitions confirms that it agrees with other treatments in the literature such as \cite[3.5 Scholium]{Beke}, \cite[p. 279]{Get} and \cite[Definition 3.5]{Hen}.\footnote{At first glance, Henriques' definition appears different from that of the other sources. To see that the two are canonically isomorphic, note that $n$-groupoids are $(n+1)$-coskeletal, that the two constructions are identical on $n$-skeleta, and that this extends to an isomorphism of $(n+1)$-skeleta using the uniqueness of filling $(n+1)$-horns.}

\begin{definition}\label{d:duskin} 
    Let $f\colon X\to Y$ be a Kan fibration in $\sC$. Let $n\ge 0$.  Define the {\em relative Duskin $n$-truncation} $\tau_{\le n}(X,f)$ to be the simplicial sheaf where
    \begin{enumerate}
        \item $\tau_{\le n}(X,f)_{\le n-1}=X_{\le n-1}$, 
        \item $\tau_{\le n}(X,f)_n=\pi_0(P^{\ge n}(f))$,
        \item $\tau_{\le n}(X,f)_{n+1}=\Lambda^{n+1}_1\tau_{\le n}(X,f)_n \times_{\Lambda^{n+1}_1 Y}Y_{n+1}$, and for $k>n+1$
        \item $\tau_{\le n}(X,f)_k=\hom(\sk_{n+1}\Delta^k,\tau_{\le n}(X,f)_{\le n+1})\times_{\Csk_{n+1,k}Y}Y_k$.
    \end{enumerate}
\end{definition}
 As shown in \cite[Section 4]{W}, this defines a simplicial sheaf; in particular, by \cite[Axiom 4]{W}, \cite[Lemma 4.2]{W} determines the ``missing face map'' 
 \[
    d_1\colon \tau_{\le n}(X,f)_{n+1}=\Lambda^{n+1}_1\tau_{\le n}(X,f)_n\to \pi_0(P^{\ge n}(f))=\tau_{\le n}(X,f)_n.
 \]
 By inspection, the assignment $f\mapsto \tau_{\le n}(X,f)$ is functorial and fits into a natural commuting diagram
\begin{equation} \label{eq:dusk-diag}
  \xymatrix{
    X \ar[r]^-{q} \ar[dr]_f & \tau_{\le n}(X,f) \ar[d]^{\tau_{\le n}f} \ar[r] & \tau_{<n}(X,f) \ar[dl]^{\tau_{<n}f}  \\
    & Y
  }
\end{equation}
When $Y=\ast$, we write $\tau_{\le n}X$ and refer to this as the {\em Duskin $n$-truncation}. 

\begin{remark}
    Let $f\colon X\to Y$ be a Kan fibration as above.  For $y\in Y_0$, denote the fiber of $f$ by $F_y:=\{y\}\times_Y X$. For $Y$ a smooth $\infty$-groupoid, we show below in Lemma~\ref{l:trhomgrps} that the homotopy groups of $\tau_{<n+1}(X,f)$ and $\tau_{\le n}(X,f)$ behave exactly as one would expect from the classical theory of Postnikov sections.
\end{remark}

\begin{remark}\label{r:infpost}
    Note that, taking $n=\infty$ in Definitions~\ref{d:moore} and \ref{d:duskin}, we get $X=\tau_{\le \infty}X=\tau_{< \infty}X$, as expected.
\end{remark}

\begin{example}
    Let $X$ be a smooth $\infty$-group in $\sC$. Consider the map 
    \[
        X\to \tau_{\le 1}X.
    \]
    On vertices, this is the identity $\ast\to\ast$.  On 1-simplices, $\pi_0(P^{\ge 1}X)$ is just the sheaf of based homotopy classes of 1-simplices in $X$, which, because $X$ is reduced, is just the sheaf $\pi_1 X$. The map
    \begin{align*}
        X_1&\to (\tau_{\le 1}X)_1=\pi_1 X\intertext{sends a 1-simplex to its based homotopy class}
        x&\mapsto [x].
    \end{align*} 
    Continuing to unpack the definition, on 2-simplices, the map
    \begin{align*}
        X_2&\to (\tau_{\le 1}X)_2=\Lambda^2_1 (\tau_{\le 1}X)\cong \pi_1X^{\times 2}\intertext{sends a 2-simplex to the based homotopy class of its spine}
        x&\mapsto [\spine(x)]=([x|_{\{0<1\}}],[x|_{\{1<2\}}]).
    \end{align*}
    As $\tau_{\le 1}X$ is a smooth 1-group in $\sSh(\C)$ and thus 2-coskeletal (Lemma~\ref{l:stackcosk}), the above fully determines the map of simplicial objects.  Unpacking this, we see that, on $n$-simplices, the map $X\to\tau_{\le 1}X$ is given by
    \begin{align*}
        X_n &\to \pi_1X^{\times n}\\
        x&\mapsto [\spine(x)]
    \end{align*}
    and that this gives an isomorphism $\tau_{\le 1}X\cong N\pi_1 X$ with homogeneous coordinates on $N\pi_1X$.
\end{example}

We will need the following lemmas.
\begin{lemma}\label{l:postlim}
    Let $f\colon X\to Y$ be a Kan fibration. The natural maps over $Y$ give isomorphisms
    \begin{align*}
        X\cong \varprojlim_n \tau_{\le n}(X,f)\cong \varprojlim_n \tau_{<n}(X,f).
    \end{align*}
\end{lemma}
\begin{proof}
    This follows immediately from the observation that the maps
    \[
        X\to \tau_{\le n}(X,f)\to\tau_{<n}(X,f)
    \]
    are the identity on $(n-1)$-skeleta.
\end{proof}

Recall that by \cite[Proposition 4.5]{W}, if $\pi_0(P^{\ge n}(f))$ is representable, then $\tau_{\le n}f\colon \tau_{\le n}(X,f)\to Y$ is an $n$-stack. We can say more.

\begin{lemma}\label{l:kan}
    Let $f\colon X\to Y$ be a Kan fibration in $\sC$. 
	\begin{enumerate}
		\item Suppose that the orbit space of $P^{\ge n}(f)$ exists in the sense of \cite[Definition 2.8]{W}, i.e. that the map $X_n\to \pi_0 (P^{\ge n}(f))$ exists in $\C$ and is a cover. Then $\tau_{\le n}f\colon \tau_{\le n}(X,f)\to Y$ is an $n$-stack and the map $q\colon X\to \tau_{\le n}(X,f)$ is a Kan fibration.
		\item\label{l:kanmoore} Suppose that the map $X_n\to \Image(\mu_n(f))$ is a cover in $\C$. Then $\tau_{<n}(X,f)\to Y$ is an $n$-stack, and the map $p\colon X\to \tau_{<n}(X,f)$ is a Kan fibration.
	\end{enumerate}
\end{lemma}
\begin{proof}
	For the first statement, note that by \cite[Proposition 4.5]{W}, that $X\to \pi_0(P^{\ge n}(f))$ is a cover implies that $\tau_{\le n}f\colon \tau_{\le n}(X,f)\to Y$ is an $n$-stack. To see that $q$ is a Kan fibration, observe that for $m<n$, the map 
	\[
	   \lambda^m_i(q)\colon X_m\to \Lambda^m_iX\times_{\Lambda^m_i \tau_{\le n}(X,f)} \tau_{\le n}(X,f)_m=:\Lambda^m_i(q)
	\]
	is the identity.  This implies that for $m=n$, the map $\lambda^n_i(q)$ is isomorphic to
    \[
        X_n\to \tau_{\le n}(X,f)_n=\pi_0(P^{\ge n}(f))
    \]
    for all $i$, which is a cover, by assumption. Next, observe that for all $m$ and $i$, we have a commuting diagram
    \[
        \xymatrix{
            X_m \ar[r] \ar[d] & \tau_{\le n}(X,f)_m \ar[r] \ar[d] & Y_m \ar[d] \\
            \Lambda^m_i X \ar[r] & \Lambda^m_i\tau_{\le n}(X,f) \ar[r] & \Lambda^m_i Y
        }
    \]
    By \cite[Proposition 4.5]{W}, for $m>n$, the right square in this diagram is a pullback. By the universal property of pullbacks, the above diagram determines an isomorphism
    \[
        \Lambda^m_i(q):=\Lambda^m_iX\times_{\Lambda^m_i\tau_{\le n}(X,f)}\tau_{\le n}(X,f)_m\to^\cong \Lambda^m_i X\times_{\Lambda^m_i Y}Y_m=\Lambda^m_i(f)
    \]
    which sits in a commuting triangle
    \[
        \xymatrix{
            X_m \ar[rr]^-{\lambda^m_i(q)} \ar[drr]_-{\lambda^m_i(f)} && \Lambda^m_i(q) \ar[d]^\cong \\
            && \Lambda^m_i(f)
        }
    \]
    Because $f$ is Kan, all these maps are covers, i.e. $q$ is Kan as well.

	For the second statement, first note that the proof of \cite[Proposition 7.9]{Hen} shows, {\em mutatis mutandis}, that $\tau_{< n}(X,f)$ is Kan if a) $f$ is a Kan fibration and b) $X\to \tau_{<n}(X,f)$ is a map of objects in $\C$.  Indeed, the only changes we need to make are a) in lieu of Henriques' $B_n$, we need to consider the presheaf
    \[
        B_n(f):=\Image(\mu_n(f))
    \]
    which specializes to $B_n$ in the case when $Y=\ast$,\footnote{Note that there is a typo in \cite[Section 7]{Hen} so that what should be $B_n$, e.g. as in the proof of \cite[Proposition 7.9]{Hen}, is instead written as $B_{n-1}$.} and b) our assumption that $X_n\to \Image(\mu_n(f)):=B_n(f)$ is a cover allows us to use Axiom~\ref{a:rightcancellation} for categories with covers in lieu of \cite[Lemma 4.8]{Hen}. 
 
    That $p$ is a Kan fibration now follows by a similar argument as the above.  We begin by observing that $\tau_{<n}(X,f)$ is relatively $n$-coskeletal, i.e.
    \[
        \tau_{<n}(X,f)\cong \Csk_n\tau_{<n}(X,f)\times_{\Csk_nY}Y.
    \]
    As a result, for $m>n$, just as above, we have a canonical isomorphism
    \[
        \Lambda^m_i(p):=\Lambda^m_iX\times_{\Lambda^m_i\tau_{<n}(X,f)}\tau_{<n}(X,f)_m\to^\cong \Lambda^m_iX\times_{\Lambda^m_iY}Y_m=:\Lambda^m_i(f).
    \]
    We conclude that the morphism $\lambda^m_i(p)$ is isomorphic to $\lambda^m_i(f)$, and thus a cover because $\lambda^m_i(f)$ is. Similarly, for $m<n$, the map $\lambda^m_i(p)$ is the identity.  

    For $m=n$, unraveling the definition, we see that the map $\lambda^n_i(p)$ is isomorphic to the map 
    \[
        X_n\to B_n(f).
    \]
    By assumption this is a cover in $\C$.  
\end{proof}

We will want the following standard lemmas.
\begin{lemma}\label{l:trhomgrps}
    Let $Y$ be a smooth $\infty$-groupoid and $f\colon X\to Y$ a Kan fibration.  For $y\in Y_0$, let $F_y$ denote the fiber $F_y:=\{y\}\times_Y X$.  For any $x\in F_{y,0}$, The map 
    \[
        X\to\tau_{\le n}(X,f)
    \]
    induces an isomorphism
    \begin{align*}
        \pi_{n+1}(X,x)/\pi_{n+1}(F,x)&\cong \pi_{n+1}(\tau_{\le n}(X,f),x)\intertext{and, for $i\le n$, induces isomorphisms}
        \pi_i(X,x)&\cong \pi_i(\tau_{\le n}(X,f),x).\intertext{For $i>n+1$, the map $\tau_{\le n}(X,f)\to Y$ induces isomorphisms}
        \pi_i(\tau_{\le n}(X,f),x)&\cong \pi_i(Y,y).
    \end{align*}
    The same is true with $\tau_{<n+1}$ in lieu of $\tau_{\le n}$ above.
\end{lemma}
\begin{proof}
    For $\C=\Set$ this is classical.  Indeed, the map $X\to\tau_{\le n}(X,f)$ gives a map of fibrations over $Y$. By the five-lemma applied to the associated map of long exact sequences, we reduce to showing that $\pi_i(F_y,x)=0$ for $i>n$, and this follows immediately from the definitions.  A similar argument shows this for $\tau_{<n+1}$.
    
    For general $\C$, by a standard argument \cite[Section `For Logical Reasons']{Beke}, we reduce to the case where $\Sh(\C)$ has enough points. Because $\tau_{\le n}$ and $\tau_{<n+1}$ commute with passing to stalks, we see from the corresponding result for simplicial sets, that all of the maps in question become isomorphisms on all stalks.  They are thus isomorphisms of sheaves as claimed.
\end{proof}

\begin{lemma}\label{l:trtr}
    Let $f\colon X\to Y$ be a Kan fibration. 
    \begin{enumerate}
        \item For $k>0$, the natural map $\tau_{\le n}(X,f)\to^\cong\tau_{<n+k}(\tau_{\le n}(X,f),\tau_{\le n}f)$ is an isomorphism over $Y$.
        \item For $k\ge 0$, the natural map $\tau_{<n}(X,f)\to^\cong \tau_{\le n+k}(\tau_{<n}(X,f),\tau_{<n}f)$ is an isomorphism over $Y$.
        \item Let $n<m$.  Then the natural maps
            \begin{align*}
                \tau_{\le n}(X,f)&\to^\cong \tau_{\le m}(\tau_{\le n}(X,f),\tau_{\le n} f)\\
                \tau_{\le n}(\tau_{\le m}(X,f),\tau_{\le m}f)&\to^\cong \tau_{\le n}(X,f)
            \end{align*}
            are isomorphisms over $Y$. Moreover, the same statement holds with $\tau_{<n}$ and $\tau_{<m}$ in lieu of $\tau_{\le n}$ and $\tau_{\le m}$.
    \end{enumerate}
\end{lemma}
\begin{proof}
    For $\C=\Set$ this is a standard exercise in the definitions. For the general case, we reason just as in the proof of Lemma~\ref{l:trhomgrps}.
\end{proof}

\begin{lemma}\label{l:duskmoorepull}
    Let $f\colon X\to Y$ be a Kan fibration, and let $Z\to Y$ be a map in $\sC$ such that the pullback $Z_0\times_{Y_0} X_0\to Z_0$ exists in $\sC$. 
    \begin{enumerate}
        \item Suppose $\tau_{\le n}(X,f)\to Y$ exists in $\sC$. Then the canonical map over $Z$ gives an isomorphism
        \[
            \tau_{\le n}(Z\times_Y X,\pi_Z)\to^\cong Z\times_Y\tau_{\le n}(X,f).
         \]
        \item Suppose $\tau_{<n}(X,f)\to Y$ exists in $\sC$. Then the canonical map over $Z$ gives an isomorphism
        \[
            \tau_{<n}(Z\times_Y X,\pi_Z)\to^\cong Z\times_Y\tau_{<n}(X,f).
         \]
    \end{enumerate}
\end{lemma}
\begin{proof}
    By \cite[Theorem 2.17(4)]{W}, the assumptions guarantee that the pullback 
    \[
        Z\times_Y X\to Z
    \]
    exists in $\sC$ and is a Kan fibration. 
    
    For the first statement, the functoriality of $\tau_{\le n}$ and the universal property of pullbacks give the canonical map
    \[
        \tau_{\le n}(Z\times_Y X,\pi_Z)\to Z\times_Y\tau_{\le n}(X,f).
    \]
    By inspection, this map is an isomorphism on $(n-1)$-skeleta. By Lemma~\ref{l:kanngpd}, the result follows once we show it's an isomorphism on $n$-simplices, i.e. once we show that the canonical map induces an isomorphism
    \[
        \pi_0 P^{\ge n}(\pi_Z)\to^\cong Z_n\times_{Y_n}\pi_0 P^{\ge n}(f).
    \]
    But this follows directly from the definition of the relative path-space, i.e. from unpacking the formulas on \cite[p. 545]{W}.

    For the second statement, the analogous reasoning gives a canonical map
    \[
        \tau_{<n}(Z\times_Y X,\pi_Z)\to Z\times_Y \tau_{< n}(X,f).
    \]
    By inspection, this map is an isomorphism on $(n-1)$-skeleta. By Lemmas~\ref{l:trtr}, both maps are $n$-stacks over $Z$, thus by Lemma~\ref{l:kanngpd}, it suffices to show the canonical map is an isomorphism on $n$-simplices, i.e. that it induces an isomorphism
    \[
        \Image(\mu_n(\pi_Z))\to^\cong Z_n\times_{Y_n}\Image(\mu_n(f)),
    \]
    and this follows directly from inspection.
\end{proof}

\begin{lemma}\label{l:mooredusk}
    Let $n>0$ and let $f\colon X\to Y$ be a Kan fibration such that $\tau_{<n+1}(X,f)$ and $\tau_{\le n}(X,f)$ exist in the sense of Lemma~\ref{l:kan}. Then the map $p\colon \tau_{<n+1}(X,f)\to \tau_{\le n}(X,f)$ is an $(n+1)$-hypercover.
\end{lemma}
\begin{proof}
    The map $p$ is a local weak equivalence, by inspection (see e.g. the discussion on \cite[p. 1026]{Hen}). By Lemmas~\ref{l:kan} and~\ref{l:trtr}, our assumptions on $f$ imply that the map $\tau_{<n+1}(X,f)\to \tau_{\le n}(X,f)$ is a Kan fibration. By \cite[Theorem 6.7]{RZ}, it is therefore a hypercover. By Lemma~\ref{l:kanngpd} and Lemma~\ref{l:kan}, it is also an $(n+1)$-stack. Therefore, by \cite[Theorem 2.17(2)]{W}, it is an $(n+1)$-hypercover.
\end{proof}

The following provides an analogue for $\tau_{<n}$ of \cite[Proposition 4.6]{W}.
\begin{lemma}\label{l:moorehyp}
    Let $f\colon X\to Y$ be a hypercover. 
    \begin{enumerate}
        \item The natural map $\tau_{\le n}(X,f)\to\tau_{<n}(X,f)$ is an isomorphism over $Y$. In particular, $\tau_{<n}f\colon \tau_{<n}(X,f)\to Y$ is an $n$-hypercover.
        \item The map $p\colon X\to \tau_{<n}(X,f)$ is a hypercover.
    \end{enumerate}
\end{lemma}
\begin{proof}
    Observe that by \cite[Proposition 4.6(1)]{W}, the map $\tau_{\le n}(X,f)\to Y$ is an $n$-hypercover, i.e. for all $k\ge n$, the maps
    \[
        \mu_k(\tau_{\le n}f)\colon \tau_{\le n}(X,f)_k\to M_k(\tau_{\le n}f)
    \]
    are isomorphisms. From the definition of $\tau_{<n}$, this implies that
    \[
        \tau_{\le n}(X,f)\cong \tau_{<n}(\tau_{\le n}(X,f),\tau_{\le n}f).
    \]
    The first statement now follows by Lemma~\ref{l:trtr}, and the second follows by \cite[Proposition 4.6(2)]{W}.
\end{proof}

\begin{lemma}\label{l:tauexact}
    Let $f\colon X\to Y$ be map of smooth $\infty$-groupoids. 
    \begin{enumerate}
        \item  Suppose the conditions of Lemma~\ref{l:kan} hold for $X$ and $Y$ so that $\tau_{\le n} f\colon \tau_{\le n} X\to \tau_{\le n} Y$ exists in $\sC$. 
            \begin{enumerate}
                \item If $f$ is a Kan fibration, then $\tau_{\le n} f\colon \tau_{\le n} X\to \tau_{\le n} Y$ is a Kan fibration.
                \item If $f$ is a weak equivalence, then so is $\tau_{\le n} f\colon \tau_{\le n}X\to \tau_{\le n}Y$.
            \end{enumerate}
        \item  Similarly, suppose the conditions of Lemma~\ref{l:kan} hold for $X$ and $Y$ so that $\tau_{<n} f\colon \tau_{<n} X\to \tau_{<n} Y$ exists in $\sC$. 
            \begin{enumerate}
                \item If $f$ is a Kan fibration, then $\tau_{<n} f\colon \tau_{<n} X\to \tau_{<n} Y$ is a Kan fibration.
                \item If $f$ is a weak equivalence, then so is $\tau_{<n} f\colon \tau_{<n}X\to \tau_{<n}Y$.
            \end{enumerate}
    \end{enumerate}           
\end{lemma}
\begin{proof}
    Assume $\tau_{\le n}f$ exists.  Consider the commuting square 
    \[
        \xymatrix{
        X \ar[r]^f \ar[d] & Y\ar[d]^{q} \\
        \tau_{\le n}X \ar[r]^-{\tau_{\le n}f} & \tau_{\le n}Y
        }
    \]
    First, suppose $f$ is a Kan fibration.  The vertical maps are covering Kan fibrations by Lemma~\ref{l:kan}. Because Kan fibrations compose \cite[Theorem 2.17]{W}, $qf$ is a Kan fibration. By Lemma~\ref{l:kanngpd}, $\tau_{\le n}f$ is therefore a Kan fibration.

    Next suppose $f$ is a weak equivalence. We reduce, as in the proof of Lemma~\ref{l:trtr} to the case where $\Sh(\C)$ has enough points. Because $\tau_{\le n}$ commutes with passing to stalks, the result now follows from the classical result for simplicial sets (see\cite[Proposition 3.5]{Beke} and \cite{Glenn}).
    
    The statements for $\tau_{<n}$ follow {\em mutatis mutandis} from the same argument as above.
\end{proof}

\begin{example}
    Note that neither $\tau_{<n}$ nor $\tau_{\le n}$ preserve arbitrary pullbacks. To see this for simplicial sets, let $G$ be a group and $X$ a free $G$-set. Consider the pullback square of Kan fibrations
    \[
        \xymatrix{
            X   \ar[r] \ar[d] & X//G \ar[d]  \\
            \Delta^0 \ar[r] & NG
            }.
    \]
   Applying $\tau_{<1}$ to the above square, we get the square
    \[
        \xymatrix{
            X \ar[r] \ar[d] & X//G \ar[d]  \\
            \Delta^0 \ar[r] & \Delta^0
            }.
    \]
    where the identification $\tau_{<1}X//G=X//G$ follows from our assumption that the action is free, and $\tau_{<1}NG=\Delta^0$ and $\tau_{<1}X=X$ by inspection.  In particular, this square is not a pullback.

    Similarly, applying $\tau_{\le 0}$ to the above square, we get
    \[
        \xymatrix{
            X \ar[r] \ar[d] & \Delta^0\ar[d]  \\
            \Delta^0 \ar[r] & \Delta^0
            }.
    \]
    which is also not a pullback.
\end{example}

\subsection{Minimal Kan fibrations}\label{s:min}
\begin{definition}\label{d:min}
    Let $f\colon X\to Y$ be a Kan fibration.  We say $f$ is {\em minimal} if $X_0\cong \pi_0(Y_0\times_Y X)$, and for all $n>0$, the natural map gives an isomorphism
    \[
        X_n\to^\cong \pi_0(P^{\ge n}(f)).
    \]
\end{definition}

Just as for simplicial sets, we have the following.
\begin{lemma}\label{l:duskmoore}
    Let $n>0$ and let $f\colon X\to Y$ be a Kan fibration such that $\tau_{\le n}(X,f)$ and $\tau_{<n}(X,f)$ exist in the sense of Lemma~\ref{l:kan}. Then the map 
    \[
        q\colon \tau_{\le n}(X,f)\to \tau_{<n}(X,f)
    \]
    is a minimal Kan fibration.
\end{lemma}
\begin{proof}
    The map $q$ is an isomorphism on $(n-1)$-skeleta, so it satisfies the minimality below dimension $n-1$. Moreover, unpacking the definition of $P^{\ge n-1}(q)$, we see that $\pi_0P^{\ge n-1}(q)=P^{\ge n-1}(q)_0=X_{n-1}$. Last, by Lemmas~\ref{l:kanngpd},~\ref{l:kan}, and~\ref{l:trtr}, $q$ is an $n$-stack. Therefore, by \cite[Theorem 3.6]{W} $P^{\ge n}(q)$ is a smooth $0$-groupoid, i.e. it is the object $\tau_{\le n}(X,f)_n$ views as a constant simplicial object.  We conclude that $q$ is a minimal Kan fibration as claimed.
\end{proof}

We will need the following lemmas.
\begin{lemma}\label{l:minpost}
    A Kan fibration $f\colon X\to Y$ is minimal if and only if the maps $\tau_{<n+1}(X,f)\to \tau_{\le n}(X,f)$ are isomorphisms for all $n$. 
\end{lemma}
\begin{proof}
    By Lemma~\ref{l:mooredusk}, the map $\tau_{<n+1}(X,f)\to \tau_{\le n}(X,f)$ is an $(n+1)$-hypercover. It is the identity on $(n-1)$-skeleta, so we see that it is an isomorphism if and only if it is an isomorphism on $n$-simplices. But on $n$-simplices, this map is precisely the map
    \[
       \tau_{<n+1}(X,f)=X_n\to \pi_0P^{\ge n}(f)=\tau_{\le n}(X,f).
    \]
    This map is an isomorphism for all $n$ if and only if $f$ is minimal.
\end{proof}

\begin{lemma}\label{l:duskmooremin}
    Let $f\colon X\to Y$ be a minimal Kan fibration. For all $n\ge 0$, $\tau_{<n+1}(X,f)\cong\tau_{\le n}(X,f)$ and the map $\tau_{<n+1}f\cong \tau_{\le n}f\colon \tau_{\le n}(X,f)\to Y$ is a minimal Kan fibration in $\sC$. 
\end{lemma}
\begin{proof}
    The isomorphism $\tau_{<n+1}f\cong\tau_{\le n}f$ is Lemma~\ref{l:minpost}. The minimality of $f$ plus Lemma~\ref{l:kan} implies that $\tau_{\le n}f\colon \tau_{\le n}(X,f)\to Y$ is a Kan fibration in $\sC$. That it is minimal, follows from the minimality of $f$ and \cite[Theorem 3.6]{W}.
\end{proof}

\begin{lemma}\label{l:duskmin}
    Let $f\colon X\to Y$ be a minimal Kan fibration of smooth $\infty$-groupoids and suppose that $\tau_{\le n}f\colon \tau_{\le n}X\to \tau_{\le n}Y$ exists in $\sC$. Then $\tau_{\le n}f$ is a minimal Kan fibration. 
\end{lemma}
\begin{proof}
    If $n=0$, then there is nothing to check.  Suppose $n>0$. By definition, the map $\tau_{\le n}f$ agrees with $f$ on $(n-1)$-skeleta.  Moreover, by inspection, $\pi_0P^{\ge n-1}(f)\cong \pi_0P^{\ge n-1}(\tau_{\le n}f)$. We conclude that $\tau_{\le n}f$ satisfies the minimality condition on $(n-1)$-skeleta.  Because $\tau_{\le n}X$ and $\tau_{\le n}Y$ are smooth $n$-groupoids, the map $\tau_{\le n}f$ is an $n$-stack. By \cite[Theorem 3.6]{W}, $P^{\ge n}(\tau_{\le n}f))$ is a smooth $0$-groupoid, i.e. it is the object $(\tau_{\le n}X)_n$ viewed as a constant simplicial object.  Therefore $\tau_{\le n}f$ is an $n$-stack which satisfies the minimality condition on $n$-skeleta.  We conclude that it is minimal.
\end{proof}

\begin{lemma}\label{l:minkan}\mbox{}
    \begin{enumerate}
        \item Minimal Kan fibrations pull back: i.e. given a minimal Kan fibration $f\colon X\to Y$ and a map $g\colon Z\to Y$ in $\sC$ such that the pullback $Z\times_Y X$ exists in $\sC$, then $Z\times_Y X\to Z$ is a minimal Kan fibration.
        \item  Suppose we have a commuting triangle in $\sC$
              \[
                \xymatrix{
                X \ar[rr]_\sim^\varphi \ar[dr] && Y \ar[dl] \\
                & Z
                }
            \]
            in which both maps to $Z$ are minimal Kan fibrations and $\varphi$ is a local weak equivalence. Then $\varphi$ is an isomorphism.
        \item Any square in $\sC$
        \[
            \xymatrix{
                X \ar[r]^\sim \ar[d] & W \ar[d] \\
                Y \ar[r]^\sim & Z
            }
        \]
        in which the vertical arrows are minimal Kan fibrations and the horizontal arrows are local weak equivalences is a pullback square.
    \end{enumerate}
\end{lemma}
\begin{proof}
    For the first two statements, by a standard argument \cite[Section `For Logical Reasons']{Beke}, we reduce to the case where $\Sh(\C)$ has enough points. Our axioms for covers guarantee that representables are sheaves; in particular, to check that a map in $\sC$ is an isomorphism, it suffices to check it for all stalks. Observe that a map is a minimal Kan fibration if and only if it is a minimal Kan fibration on all stalks (because for all $n$ and $f$, the map of sheaves $X_n\to \pi_0P^{\ge n}(f)$ is an isomorphism if and only if it is one on all stalks, and similarly for $X_0\to \pi_0(Y_0\times_Y X_0)$). The first statement now follows from the classical fact that minimal Kan fibrations of simplicial sets pull back \cite[Lemma III.5.4]{BGM}, while the second follows from the classical fact that weak equivalences between minimal Kan fibrations of simplicial sets are isomorphisms \cite[Lemma I.10.4]{GJ}. 

    The third statement follows from the first two.  Indeed, consider the map  (in $\sSh(\C)$ if need be)
    \[
        X\to Y\times_Z W
    \]
   By the first statement, this is a map between minimal fibrations over $Y$. It is a local weak equivalence over $Y$ by Lemma~\ref{l:rightproper} and the 2-of-3 property for weak equivalences. By the second statement, it is an isomorphism.
\end{proof}

\begin{example}\mbox{}
    Note that the (absolute) Moore truncation does not preserve minimal Kan fibrations, i.e. the first statement in the lemma is not true if we take $\tau_{<n}$ in lieu of $\tau_{\le n}$. To see this for simplicial sets, let $G$ be a group and consider the minimal Kan fibration $G//G\to NG$.  Applying $\tau_{<1}$ we obtain the Kan fibration $G//G\to \Delta^0$ which is not minimal!

    Note also that minimal fibrations do not compose (contrary to some assertions in the literature, e.g. \cite[Proposition 4.5.20(i)]{FP}). In the example above, $G//G\to NG$ and $NG\to \Delta^0$ are both minimal fibrations, but the composite $G//G\to \Delta^0$ is not!
\end{example}

\subsection{Smooth $k$-invariants} \label{sec:smooth-kinv}
We are now in a position to construct $k$-invariants in $\sC$. Classically, one begins with the observation above that for a Kan fibration $f\colon X\to Y$, the map
\[
    \tau_{\le n}(X,f)\to\tau_{<n}(X,f)
\]
is a minimal Kan fibration (Lemma~\ref{l:duskmoore}).  Using the axiom of choice, Barratt, Gugenheim and Moore showed that every minimal Kan fibration of simplicial sets is a fiber bundle \cite[Proposition IV.2.2]{BGM} and they constructed a trivializing cocycle for this fiber bundle \cite[Proposition IV.3.2, pp. 655-656]{BGM}, which in the above context is just the $k$-invariant.  The arguments of \cite{BGM} do not generalize because we cannot appeal to arguments based on choice in arbitrary categories with covers $\C$. 

In \cite[Theorem 1.4]{W2}, the second author showed that for any category with covers $\C$, given a principal bundle in $\sC$, we can functorially extract a trivializing cocycle. Here, we take the formulas from \cite{W2} as an ansatz to define the $k$-invariants directly.  We use this to prove that the above story for Kan complexes applies to any Kan fibration in $\sC$.

We can now state the main theorem of this section. We defer the complete proof until Section~\ref{app:kinvar}. 
\begin{theorem}\label{t:kinvar} 
    Let $Y$ be a smooth $\infty$-group in $\sC$, and assume that $\pi_1Y$ is a group in $\C$. Let $f\colon X\to Y$ be a Kan fibration such that $\tau_{<n}f$ is an isomorphism for some $n>1$, and let $F:=X\times_Y \Delta^0$. Suppose that $\pi_nF$ is a group object in $\C$.  Then there exists a commuting diagram of smooth $\infty$-groups, natural in $f$, 
    \begin{equation}\label{e:kinvar1}
        \xymatrix{
         X\times_Y \LW(X/Y) \ar[r]^-\psi \ar[d]_{f|_{\LW(X/Y)}} & WK(\pi_nF,n)//\pi_1 Y \ar[d] \\
         \LW(X/Y) \ar[r]^-{\varphi} & K(\pi_nF,n+1)//\pi_1 Y
        }
    \end{equation}
    such that the restriction of $\psi$ to the fibers
    \[
        F\cong X\times_Y\LW(X/Y)\times_{\LW(X/Y)}\Delta^0\to^\psi WK(\pi_nF,n)//\pi_1Y\times_{K(\pi_nF,n+1)//\pi_1Y}\Delta^0\cong K(\pi_n F,n)
    \]
    induces the canonical isomorphism 
    \[
        \psi_\ast\colon \pi_nF\to^\cong \pi_n K(\pi_nF,n).
    \]
    In particular, the square above factors through a pullback square in $\sC$
    \begin{equation}\label{e:kinvar2}
        \xymatrix{
         \tau_{\le n}(X,f)\times_Y \LW(X/Y) \ar[r]^-{\bar{\psi}} \ar[d] & WK(\pi_nF,n)//\pi_1 Y \ar[d] \\
         \LW(X/Y) \ar[r]^-{\varphi} & K(\pi_nF,n+1)//\pi_1 Y
        }.
    \end{equation}
\end{theorem}

\begin{remark}\mbox{}
    \begin{enumerate}
        \item The classical construction shows that $\pi_1Y$ acts on $\pi_nF$ by automorphisms. Our assumption that $n>1$ further guarantees that $K(\pi_nF,n)$ is a simplicial group, and thus that right vertical map in the diagram above is indeed well-defined.
        \item If $\pi_1Y$ is not a group in $\C$, we can still interpret the theorem in $\Sh(\C)$. In this case, by the Yoneda lemma, our assumption on $\pi_nF$ will still guarantee that the triangle in the diagram above exists in $\C$.
        \item The key aspects of the proof are explicit formulas defining the local $k$-invariant $\varphi$, and the map $\psi$ fitting into the square~\eqref{e:kinvar1} and restricting to the fiber as described. Given these, the pullback square~\eqref{e:kinvar2} follows from a straightforward argument using the properties of minimal fibrations we established above.  As remarked above we defer the complete proof to Section~\ref{app:kinvar}.
    \end{enumerate}
\end{remark}

\subsubsection{Postnikov Sections are Principal Bundles}
To conclude our treatment of Postnikov theory in $\sC$, we now use results of \cite{W} to show that for smooth $\infty$-groups $X$ in many geometric contexts, the Postnikov section
\[
    \tau_{\le n}X\to \tau_{<n}X
\]
is a principal bundle for the $N\pi_1X$-group $K(\pi_n X,n)//\pi_1X\to N\pi_1X$. 

To set this up, recall (cf. \cite[Definition 1.2]{W2}) that, for $B\in\sC$ a (Kan) {\em principal bundle} for a $B$-group $G\to B$ consists of a covering Kan fibration $f\colon P\to Y$ in $\sC_{/B}$ equipped with an action $\rho\colon G\times_B P\to P$ over $Y$ fitting into the following commuting diagram
\begin{equation*}
    \xymatrix{
        G\times_B P \ar@<-.5ex>[r]_{\pi_P} \ar@<.5ex>[r]^\rho \ar[d]^\cong_{\rho\times \pi_P} & P \ar@{=}[d] \ar[r] & Y \ar@{=}[d]\\
        P\times_Y P \ar@<-.5ex>[r]_{\pi_2} \ar@<.5ex>[r]^{\pi_1} & P \ar[r] & Y
    }.
\end{equation*}
Note that our assumption that covers are effective epimorphisms implies that both forks in the above diagram are coequalizers.  Recall that a {\em twisted Cartesian product} (cf. \cite{BGM,W}) is a map $f\colon X\to Y$ in $\sC$ for which there exists some $Z\in \sC$ and isomorphisms for all $k\ge 0$
\[
    \xymatrix{
        X_k \ar[rr]^{\varphi_k}_\cong \ar[dr]_f && Y_k\times Z_k \ar[dl]^{\pi_1} \\
        & Y_k
    }
\]
such that, for all $i$,  
\begin{align*}
    \varphi_k s_i&=(s_i\times s_i)\varphi_{k-1}\intertext{and, for $i>0$,}
    \varphi_{k-1}d_i&=(d_i\times d_i)\varphi_k.
\end{align*}
A {\em local $n$-bundle} \cite[Definition 5.2]{W} is a twisted Cartesian product which is also an $n$-stack..\footnote{Our treatment continues to follow \cite{W}, but with the orientation reversal on simplices mentioned above.}

Recall from \cite[Section 5]{W} that the subcategory of {\em regular embeddings of $\C$} is the smallest subcategory which contains all graphs of morphisms $g\colon W\to Z$
\[
    \Gamma_g\colon W\times_Z Z\into W\times Z
\]
and is closed under pullback along covers. A {\em regular smooth $n$-groupoid in $\C$} is a smooth $n$-groupoid $Z$ such that the map
\[
    \mu_1\colon Z_1\to M_1 Z\cong Z_0\times Z_0
\]
is a regular embedding. We say that {\em Godemont's Theorem holds in $\C$} (cf. \cite[Axiom 5]{W}) if for every regular smooth $n$-groupoid $X\in sC$, the orbit space $\pi_0 X$ exists in $\C$.  

We will need the following, which follows {\em mutatis mutandis} from the proof of \cite[Proposition 4.7]{W}.
\begin{lemma}\label{l:hyptrun}
    Let Godemont's Theorem hold in $\C$.  Let $B\in\sC$ and let $Y\to B$ be an $n$-stack.  Let $X\to B$ and let $f\colon W\to X$ be a hypercover.  Then any map $W\to Y$ in $\sC_{/B}$ factors uniquely through $\tau_{\le n}(W,f)$ as in 
    \[
        \xymatrix{
            W \ar[dr]_f \ar[r] & \tau_{\le n}(W,f) \ar[d]^{\tau_{\le n}(f)} \ar[r] & Y \ar[d]\\
            & X \ar[r] & B
        }.
    \]
\end{lemma}

\begin{proposition}\label{p:postbund}
   Let Godemont's Theorem hold in $\C$. Let $Y\in\sC$ be a smooth $\infty$-group. Let $f\colon X\to Y$ be a Kan fibration such that $\tau_{<n}f$ is an isomorphism for some $n>1$ and let $F:=X\times_Y Y_0$. Suppose that $\pi_nF$ is a group object in $\C$.  Then $\tau_{\le n}(X,f)\to Y$ is a principal bundle in $\sC$ for the $N\pi_1Y$-group $K(\pi_nF,n)//\pi_1 Y\to N\pi_1Y$.  In particular, for each $k$, $\tau_{\le n}(X,f)_k\to Y_k$ is a locally trivial $\pi_1 Y$-twisted principal $\pi_n F^{\binom{k}{n}}$-bundle in $\C$.
\end{proposition}
\begin{proof}
    Note that the statement about the maps on $k$-simplices follows immediately from the description of $K(\pi_n F,n)//\pi_1 Y$ above and the claim that $\tau_{\le n}(X,f)\to Y$ is a principal $K(\pi_n F,n)//\pi_1Y\to N\pi_1Y$ bundle.  We establish this now.
    
    Observe that the map 
    \[
        WK(\pi_nF,n)//\pi_1 Y\to K(\pi_nF,n+1)//\pi_1Y
    \]
    is a principal $K(\pi_nF,n)//\pi_1Y\to N\pi_1Y$ bundle, and a minimal $n$-stack.

    We now use this to show that $\tau_{\le n}(X,f)\to Y$ exists in $\sC$ (as opposed to just $\sSh(\C)$). By Theorem~\ref{t:kinvar}, we have a commuting diagram 
    \[
        \xymatrix{
            \tau_{\le n}(X,f)\times_Y\LW(X/Y) \ar[d] \ar[r]^{\bar{\psi}} & W K(\pi_nF,n)//\pi_1 Y\ar[d] \\
            \LW(X/Y) \ar[r]^-{\varphi} \ar[d] & K(\pi_nF,n+1)//\pi_1 Y \ar[d]\\
            Y \ar[r] & N\pi_1 Y
     }
    \]
    in which the top square is a pullback diagram. By \cite[Theorem 6.7]{W}, the upper right vertical map is an $n$-stack, while the lower right vertical map is an $(n+1)$-stack.  Therefore, by Lemma~\ref{l:hyptrun}, the map $\varphi$ factors through $\tau_{\le n+1}(\LW(X/Y),\LW(f))$ as the composition of the horizontal maps in the diagram
    \[
        \xymatrix{
            \LW(X/Y) \ar[r]^-\sim \ar[dr] & \tau_{\le n+1}(\LW(X/Y),\LW(f)) \ar[r] \ar[d] & K(\pi_nF,n+1)//\pi_1Y \ar[d]\\
            & Y \ar[r] & N\pi_1Y
        }.
    \]
    By \cite[Proposition 4.6]{W}, the middle vertical map is a $(n+1)$-hypercover.
    
    Now, let $p\colon \tau_{\le n}(X,f)\times_Y \LW(X/Y)\to\tau_{\le n+1}(\LW(X/Y),\LW(f))$ denote the composite
    \[
        p\colon \tau_{\le n}(X,f)\times_Y\LW(X/Y)\to \LW(X/Y)\to^\sim\tau_{\le n+1}(\LW(X/Y),\LW(f)).
    \]
   By Lemma~\ref{l:trtr}, we see that the map $\bar{\psi}$ factors through $\tau_{\le n}(\tau_{\le n}(X,f)\times_Y\LW(X/Y),p)$ as in the diagram
      \[
        \xymatrix{
            \tau_{\le n}(X,f)\times_Y\LW(X/Y) \ar[d] \ar[r]^-\sim & \tau_{\le n}(\tau_{\le n}(X,f)\times_Y\LW(X/Y),p) \ar[d]_{\tau_{\le n}p} \ar[r] & WK(\pi_nF,n)//\pi_1Y \ar[d]\\
            \LW(X/Y) \ar[r]^-\sim \ar[dr] & \tau_{\le n+1}(\LW(X/Y),\LW(f)) \ar[r] \ar[d] & K(\pi_nF,n+1)//\pi_1Y \ar[d]\\
            & Y \ar[r] & N\pi_1Y
        }.
    \]
    By Theorem~\ref{t:kinvar}, we conclude that the top right horizontal map induces an isomorphism on $\pi_n$ of the fibers of the upper middle and upper right vertical maps, and thus, by Lemma~\ref{l:trhomgrps} and our assumption on $f$, a weak equivalence on the fibers. By inspection, the upper right vertical map is a minimal Kan fibration. Our assumptions on $f$ imply the same for $\tau_{\le n}p$. By Lemma~\ref{l:minkan}, we obtain a weak equivalence, and thus an isomorphism, of minimal Kan fibrations over $\tau_{\le n+1}(\LW(X/Y),\LW(f))$
     \begin{align*}
        \tau_{\le n}(\tau_{\le n}(X,f)\times_Y\LW(X/Y),p)&\to^\cong\\
        \tau_{\le n+1}(\LW(X/Y),\LW(f))&\times_{K(\pi_nF,n+1)//\pi_1Y}WK(\pi_nF,n)//\pi_1Y.\nonumber
    \end{align*}
    By the same reasoning, the canonical map gives an isomorphism over $\tau_{\le n+1}(\LW(X/Y),\LW(f))$
    \begin{align*}
        \tau_{\le n}(\tau_{\le n}(X,f)\times_Y\LW(X/Y),p)&\to^\cong \tau_{\le n}(X,f)\times_Y\tau_{\le n+1}(\LW(X/Y),\LW(f)).
    \end{align*}
    Together, these give an isomorphism over $\tau_{\le n+1}(\LW(X/Y),\LW(f))$
    \begin{align}
        \tau_{\le n}(X,f)&\times_Y\tau_{\le n+1}(\LW(X/Y),\LW(f))\to^\cong\label{e:prinbuniso}\\
        &\tau_{\le n+1}(\LW(X/Y),\LW(f))\times_{K(\pi_nF,n+1)//\pi_1Y}WK(\pi_nF,n)//\pi_1Y\nonumber.
    \end{align}
    We conclude that 
    \[
        \tau_{\le n}(X,f)\times_Y\tau_{\le n+1}(\LW(X/Y),\LW(f))\to\tau_{\le n+1}(\LW(X/Y),\LW(f))
    \]
    is a local $n$-bundle in the sense above. Let $q$ denote the map
    \[
        \tau_{\le n}(X,f)\times_Y\tau_{\le n+1}(\LW(X/Y),\LW(f))\to\tau_{\le n+1}(\LW(X/Y),\LW(f))\to Y.
    \]
    We can now invoke the main theorem \cite[Theorem 5.7]{W} to conclude that 
    \[
        \tau_{\le n}\left(\tau_{\le n}(X,f)\times_Y\tau_{\le n+1}(\LW(X/Y),\LW(f)),q\right)\to Y
    \]
    exists in $\C$. But, considering the diagram (in $\sSh(\C)$) 
    \[
        \xymatrix{
            \tau_{\le n}(X,f)\times_Y\tau_{\le n+1}(\LW(X/Y),\LW(f)) \ar[rr] \ar[dr] && \tau_{\le n}(X,f)\ar[dl]\\
            & Y
        }
    \]
    we see, by \cite[Theorem 4.5]{W} as above, that 
    \[
        \tau_{\le n}\left(\tau_{\le n}(X,f)\times_Y\tau_{\le n+1}(\LW(X/Y),\LW(f)),q\right)\cong \tau_{\le n}(X,f)
    \]
    over $Y$. In particular, $\tau_{\le n}(X,f)\to Y$ is an $n$-stack in $\sC$. 
    
    To see that it is a principal $K(\pi_n F,n)//\pi_1Y\to N\pi_1Y$-bundle, consider the action on $\tau_{\le n}(X,f)\times_Y\tau_{\le n+1}(\LW(X/Y),\LW(f))$ 
    \[
        \xymatrix{
            K(\pi_nF,n)//\pi_1 Y\times_{N\pi_1 Y}\tau_{\le n}(X,f)\times_Y\tau_{\le n+1}(\LW(X/Y),\LW(f)) \ar[dr] \ar[dd]_\rho \\
            & \LW(X/Y)\\
            \tau_{\le n}(X,f)\times_Y\tau_{\le n+1}(\LW(X/Y),\LW(f)) \ar[ur]
        }
    \]
    determined by the isomorphism~\ref{e:prinbuniso} above and the action of $K(\pi_n F,n)//\pi_1Y\to N\pi_1Y$ on $WK(\pi_nF,n)//\pi_1Y\to N\pi_1Y$.  Observe that both source and target are local $n$-bundles over $Y$.  Indeed, we have already shown this for the target.  Using that 
    \[
        K(\pi_nF,n)//\pi_1\times_{N\pi_1 Y} WK(\pi_nF,n)//\pi_1 Y\to K(\pi_n+1,F)//\pi_1Y
    \]
    is a local $n$-bundle and that local $n$-bundles are preserved under base change, {\em mutatis mutandis}, the above argument shows this for the source too.  Postcomposing with $\LW(f)$, we obtain a diagram
    \[
        \xymatrix{
            K(\pi_nF,n)//\pi_1 Y\times_{N\pi_1 Y}\tau_{\le n}(X,f)\times_Y\tau_{\le n+1}(\LW(X/Y),\LW(f)) \ar[dr] \ar[dd]_\rho \\
            & Y\\
            \tau_{\le n}(X,f)\times_Y\tau_{\le n+1}(\LW(X/Y),\LW(f)) \ar[ur]
        }
    \]
    and by \cite[Theorem 4.5]{W}, the relative $n$-truncation of both diagonal maps exists in $\sC$.  The same argument as above now shows that the relative $n$-truncation of the top diagonal map is isomorphic to $K(\pi_n F,n)//\pi_1 Y\times_{N\pi_1 Y}\tau_{\le n}(X,f)$. We thus obtain an action
    \[
        \bar{\rho}\colon K(\pi_n F,n)//\pi_1 Y\times_{N\pi_1 Y}\tau_{\le n}(X,f)\to \tau_{\le n}(X,f).
    \]
    By the functoriality of the constructions above, we see that it makes $\tau_{\le n}(X,f)\to Y$ into a principal bundle as claimed.  
\end{proof}

We can also prove the following without any assumption of Godemont's theorem.
\begin{corollary}\label{c:minbun}
    Let $f\colon X\to Y$ be a minimal Kan fibration of smooth $\infty$-groups. Let $F=X\times_Y Y_0$ have $\pi_1F=0$.  Then $\pi_nF$ is a group object in $\C$ for all $n$, and $f\colon X\to Y$ is an inverse limit of principal bundles for the $N\pi_1Y$-groups $K(\pi_nF,n)//\pi_1Y\to N\pi_1Y$.
\end{corollary}
\begin{remark}
   This is a partial analogue of \cite[Proposition III.2.2]{BGM} which states that every minimal Kan fibration of simplicial sets is a fiber bundle.
\end{remark}
\begin{proof}[Proof of Corollary~\ref{c:minbun}]
    By Lemma~\ref{l:postlim}, $X=\varprojlim_n \tau_{\le n}(X,f)$. By Lemmas~\ref{l:duskmoore} and~\ref{l:minpost}, the map
    \[
        \tau_{\le n}(X,f)\to\tau_{\le n-1}(X,f)
    \]
    is a minimal Kan fibration in $\sC$ for all $n$. In particular, the pullback 
    \[
        \tau_{\le n}(X,f)\times_{\tau_{\le n-1}(X,f)}X_0
    \]
    is a smooth $\infty$-group. Further by Proposition~\ref{p:postbund}, for $n>1$, the map
    \[
        \tau_{\le n}(X,f)_n\to\tau_{\le n-1}(X,f)_n
    \]
    is a principal $\pi_1X$-twisted $\pi_nF$ bundle. In particular, we conclude that the natural map
    \[
        \hom(\Delta^n/\partial\Delta^n,F)\to^\cong\pi_nF\to\tau_{\le n}(X,f)\times_{\tau_{\le n-1}(X,f)}X_0
    \]
    is an isomorphism (where the first map is an isomorphism by our assumption of minimality). 
    
    The result now follows from Proposition~\ref{p:postbund}, noting that the sole role of the assumption of Godemont's theorem was to guarantee the existence of $\tau_{\le n}(X,f)$ in $\sC$, which as noted above, exists here as a consequence of minimality.
\end{proof}

\begin{example}\label{ex:nomin}
    Let $\g$ be a simple Lie algebra of compact type, with integrating simple compact Lie group $G$. There exist a wide variety of (equivalent) models of the String 2-group $\String(G)$, including as a Lie 2-group \cite{Hen,W}, as a ``2-group in stacks'' \cite{SP}, or as a simplicial Banach Lie group \cite{Hen,BSCS} to name a few. Many of these examples produce infinite dimensional objects \cite{Hen,BSCS}, and none produces a minimal model. As we now explain, this is not a defect in the methods: no minimal model of $\String(G)$ exists as a Lie 2-group. 

    Indeed, suppose that $X$ were a minimal Lie 2-group (in $\ssm$) equivalent to $\String(G)$. Then $\tau_{\le 1}X\cong NG$, and by Corollary~\ref{c:minbun}, $X\to \tau_{\le 1}X$ would be a principal $K(S^1,2)$-bundle. Further, because $H^1(G^{\times k};S^1)\cong H^2(G^{\times k};\Z)=0$, the corresponding $(S^1)^{\times \binom{k}{2}}$-bundles on the spaces of $k$-simplices $G^k$ would all be topologically trivial. Thus, the extension of 2-groups corresponding to $X\to NG$ would be classified by an element in the ``naive'' group cohomology of $G$, $H^3_{naive}(G;S^1)$. But, classical results of Hu \cite{Hu1,Hu2}, van Est \cite{vE1,vE2} and Hochschild-Mostow \cite{HM} imply that $H^3_{naive}(G;S^1)=0$ (for all simply connected compact simple Lie groups $G$).  We conclude that $X\cong K(S^1,2)\times NG$, as simplicial manifolds. But, by \cite[Theorem 8.4]{Hen} and the surrounding discussion, this is impossible!  Therefore, no minimal model for $\String(G)$ exists!
\end{example}

\begin{remark}\label{r:nomin} 
There is a second obstruction for the existence of minimal models which is closely related to Henriques' obstruction to the integrability of a Lie $n$-algebra $L$ to a Lie $n$-group \cite[Thm.\ 7.5]{Hen}. Recall, as was discussed in the introduction, that for any finite type Lie $\infty$-algebra $L$, the homotopy sheaves $\pi_i\sint L$ are finite dimensional ``diffeological groups'' \cite[Thm.\ 6.4]{Hen}, i.e. quotients of finite dimensional Lie groups by finitely generated subgroups. Henriques showed \cite[Ex.\ 7.10]{Hen} that there exists a Lie $n$-algebra $L$ for which $\pi_{n+1}\sint L$ is non-Hausdorff and hence not a Lie group. From this, he deduces that $\tau_{\leq n}\sint L$ cannot be a simplicial manifold, and hence the ``$n$-integration'' of $L$ is obstructed.
Now let $\cG$ be a minimal Lie $\infty$-group. Then, we have an embedding of spaces
    \[
        \pi_n\cG\cong \hom(\Delta^n/\partial\Delta^n,\cG)\into \cG_n.
    \]
In particular, the homotopy groups of $\cG$ must all be Hausdorff topological groups, as a Hausdorff space cannot contain a non-Hausdorff subspace. 
In particular, if $L$ is a Lie $n$-algebra with an obstructed $n$-integration,  then $\sint L$ has no minimal model. 
\end{remark}

%

\section{$L_\infty$-algebras}\label{s:lien}
We begin this section by summarizing the background material needed on $L_\infty$-algebras, Lie $n$-algebras, and their respective homotopy theories. The main reference is \cite{R}. 
Then, in \Sec \ref{sec:uni-fib} and \Sec \ref{sec:tower} we introduce universal fibrations and the theory of relative $k$-invariants for the Postnikov sections of a minimal fibration between Lie $n$-algebras. This is in direct analogy with the theory of smooth $k$-invariants presented earlier in Sec.\ \ref{sec:smooth-kinv}. Finally, we recall the homotopical properties of Henriques' integration functor \cite{Hen} in \Sec \ref{sec:MC}, as developed in the relevant sections of \cite{RZ}. Along the way, we analyze some key examples of integration which we will need later in the proof of our main theorem.

\subsection{Notation and conventions} \label{sec:notate}
We work over a field of characteristic zero, specializing to $\R$ or $\mathbb{C}$ in Sec.\ \ref{sec:MC}. By a graded vector space, we mean a $\Z$-graded vector space, not necessarily bounded in either direction. We use \underline{homological} conventions for all differential graded (dg) structures. In particular, objects traditionally thought of as non-negatively graded cochain complexes, such as the de Rham complex, will be considered as cochain complexes concentrated in non-positive degrees. As in \cite[Sec.\ 2.1]{R}, for a graded vector space $V$, we denote by $\bs V$ (resp. by $\bs^{-1} V$) the suspension (resp. the desuspension) of $V$\,. In other words,
\[
(\bs V)_{i}:=V_{i-1} \qquad (\bs^{-1} V_{i}):=V_{i+1}.
\]

We adopt the remaining notation and conventions of \cite[Sec.\ 2.1]{R} for graded linear algebra with the following \und{exception}: For $n \in \N$ we denote by $V[n]$ and $V[-n]$ the graded vector spaces
\[
V[n]:=\bs^nV, \qquad V[-n]:=\bs^{-n}V.
\] 
The reason for the small discrepancy is to insure compatibility with the usual notation for the Eilenberg-Mac Lane objects previously appearing in  Sec.\ \ref{sec:nerves}.

\subsection{$L_\infty$-algebras and their morphisms} \label{sec:Linf-morph}
Following \cite[Sec.\ 3]{R}, recall that an \df{$L_\infty$-algebra} $(L, \el)$
is a graded vector space $L$ equipped with 
a collection $\el = \{\el_1, \el_2, \el_3, \ldots \}$
of graded skew-symmetric linear maps (or brackets)
\[
\el_k \maps \Alt^k L \to L, \quad 1 \leq k < \infty 
\]
with  $\deg{\el_k} = k-2$, satisfying an infinite sequence of Jacobi-like identities of the form:
\begin{align} \label{eq:Jacobi}
   \sum_{\substack{i+j = m+1, \\ \sigma \in \Sh(i,m-i)}}
  (-1)^{\sigma}\epsilon(\sigma)(-1)^{i(j-1)} l_{j}
   (l_{i}(x_{\sigma(1)}, \dots, x_{\sigma(i)}), x_{\sigma(i+1)},
   \ldots, x_{\sigma(m)})=0
\end{align}
for all $m \geq 1$. Above, the permutation $\sigma$ ranges over all $(i,m-i)$ shuffles, 
and $\epsilon(\sigma)$ denotes the Koszul sign. In particular, Eq.\ \ref{eq:Jacobi} implies that $(L, \el_1)$ is a chain complex. 

Equivalently, an $L_\infty$-structure on a graded vector space $L$ is a degree $-1$ codifferential $\delta$ on the reduced cocommutative coalgebra $\bar{S}(\bs L) = \bigoplus_{i \geq 1} S^{i}(\bs L)$. (This is the ``Chevalley-Eilenberg coalgebra'' of $L$).  Later on, we will adopt this perspective in order to prove some of the more computationally intensive statements in this section.

A (weak) \df{$L_\infty$-morphism} $f \maps (L,\el) \to (L',\el')$
is a collection $f=\{f_1,f_2,\ldots\}$
of graded skew-symmetric linear maps
$$f_k \maps \Alt^k L \to L' \quad 1 \leq k < \infty$$
with $\deg{f_k} = k-1$, satisfying an infinite sequence of equations: For all $m \geq 1$
\begin{equation} \label{eq:Linf-morph}
\begin{split}
& \sum_{\substack{j+k = m+1, \\ \sigma \in \Sh(k,m-k)}}(-1)^{\sigma}\epsilon(\sigma)(-1)^{k(j-1)+1} f_j \bigl( l_k(x_\sigma(1), \ldots, 
x_{\sigma(k)}), x_{\sigma(k+1)},\ldots x_{\sigma(m)} \bigr)  
  \\
 & + \sum_{\substack{ 1 \leq t \leq m \\ i_{1} + \cdots i_{t} = m}}
\sum_{\tau} (-1)^{\ta}\epsilon(\ta) \chi(\ta,f)\, l'_{t} \bigl( f_{i_1}(x_{\tau(1)},\ldots,x_{\tau(i_1)}),
f_{i_2}(x_{\tau(i_1 + 1 )},\ldots,x_{\tau(i_1 + i_2)}), \\
& \quad \dots, f_{i_t}(x_{\tau(i_1 + \cdots + i_{t-1} +1 )},\ldots,x_{\tau(m)})\bigr) =0.
\end{split}
\end{equation}
Above  $\tau$ ranges through all $(i_1,\ldots,i_t)$ shuffles satisfying  
$\ta(1) < \ta(i_1 +1) < \ta (i_1 + i_2 +1) < \cdots < \ta(m-i_t +1)$, and $\chi(\ta,f) := \ti{\eps}(\ta,f)(-1)^{\frac{t(t-1)}{2} + \sum_{s=1}^{t-1}i_s(t-s)} $, where $\ti{\eps}(\ta,f)$ is the Koszul sign associated to the evaluation of $f_{i_1} \tensor f_{i_2} \tensor \cdots f_{i_t}$ with
$x_{\ta(1)} \tensor x_{\ta(2)} \tensor \cdots \tensor x_{\ta(n)}$. 

The above equations for $m=1$ imply that
$f_1 \maps (L,\el_1) \to (L',\el'_1)$ is a morphism of chain complexes, and this defines the so-called \df{tangent functor} 
\cite[Sec.\ 3.1]{R}
\begin{align} \label{eq:tan}
\tan \maps \Linf &\to \Chain\\ 
(L,\el) &\mapsto (L,\el_1).\nonumber
\end{align}
Thankfully, we can largely ignore Eq.\ \ref{eq:Linf-morph} by using the convenient fact that a morphism between $L_\infty$-algebras $L$ and $L'$ is equivalently a degree 0 morphism of dg coalgebras 
\[
F \maps \bigl(\S(\bs L), \delta \bigr) \to \bigl(\S(\bs L'), \delta' \bigr). 
\]
We will adopt this approach in some of the proofs in this section. See Appendix \ref{app:Linf} for further details. In particular, treating $L_\infty$-morphisms as dg coalgebra morphisms gives us a clear way to compose them
\cite[Eq.\ 2.6]{R}. It is typical to consider the category $\Linf$ of $L_\infty$-algebras and $L_\infty$-morphisms as a full subcategory of the category of dg cocommutative coalgebras.

As in \cite{R}, we will write morphisms in $\Linf$ using a single lower-case letter, e.g. 
\[
f \maps (L,\el) \to (L',\el'),
\]
and the $k$-ary map in the collection $f$ will always be denoted by $f_k$. It is a basic fact that $f$ is an isomorphism in $\Linf$, i.e. an \df{$L_\infty$-isomorphism}, if and only if $\tan(f)$ is an isomorphism of complexes. More generally, we say $f$ is a \df{$L_\infty$-quasi-isomorphism} if 
$\tan(f)$ is a quasi-isomorphism. This is a property on morphisms which is strictly stronger than inducing a quasi-isomorphism on the Chevalley-
Eilenberg dg coalgebras.

We say $f$ is a \df{strict $L_\infty$-morphism} if $f_k =0$ for all $k \geq 2$. In this case we write 
\[
f=f_1 \maps (L,\el) \to (L',\el')
\]
and it follows from Eq.\ \ref{eq:Linf-morph} that all $k$-ary brackets are preserved by the chain map $f_1$:
\[
\el'_k \circ f_1^{\tensor k} = f_1 \circ \el_k \quad \text{for all $k \geq 1$}.
\] 
In analogy with morphisms between Lie algebras, note that the kernel of any strict morphism $f=f_1 \maps (L,\el) \to (L',\el')$
between $L_\infty$-algebras is necessarily an \df{$L_\infty$-ideal}:
\begin{equation} \label{eq:ideal}
(\ker f_1, \el \vert_{\ker f_1}) \ideal  L.
\end{equation}
Finally, recall that for any $L_\infty$-morphism $f \maps (L,\el) \to (L',\el')$, the induced map on homology 
\[
H(f_{1}) \maps \bigl( H(L),[\cdot,\cdot]) \to (H(L'),[\cdot,\cdot]' \bigr)
\]
is a morphism of graded Lie algebras, where the above Lie brackets are
induced by the bilinear brackets $\el_2$ and $\el'_2$, respectively. In particular $H(f_1) \maps H_0(L) \to H_0(L')$ is a morphism of Lie algebras.

\subsection{Lie $n$-algebras} \label{sec:LnA}
Let $n \in \N \cup \{\infty\}$. A \df{Lie $n$-algebra} is an $L_\infty$-algebra $(L,\el)$
whose underlying graded vector space $L$ is concentrated in degrees $[0,n-1]$, i.e.\
$L= \bigoplus_{i \geq 0}^{n-1} L_i$. Hence, Lie $1$-algebras are precisely Lie algebras. A Lie $n$-algebra is \df{finite-type} if each $L_i$ is a finite-dimensional vector space.  A Lie $n$-algebra is \df{strict} if $\el_{\geq 3}=0$. Hence, strict Lie $n$-algebras are nothing but (homological) dg Lie algebras concentrated in non-negative degrees. 

\begin{example}\label{ex:semi-direct}
The following  construction provides an important supply of strict Lie $n$-algebras. Given a Lie algebra $(\g,[\cdot,\cdot]_\g)$, thought of as Lie 1-algebra concentrated in degree 0, and a dg $\g$-module $(C,d)$, we denote by $C//\g$ the \df{semi-direct product} dg Lie algebra whose underlying chain complex is $(\g \dsum C, 0 \dsum d)$ and whose $\el_2$ bracket is $[\cdot,\cdot]_\g$ in degree zero, and
\begin{equation} \label{eq:semi}
\el_2(x,c):= x \cdot c = - \el_2(c,x) \qquad \forall x \in \g \quad  \forall c \in C.
\end{equation}
Note that if $C$ is concentrated in degrees $0,\ldots,n$, then $C//\g$ is a strict Lie $(n+1)$-algebra.
\end{example}

For degree reasons, the structure maps of an $L_\infty$-morphism $f \maps (L,\el) \to (L',\el')$ between Lie $n$-algebra will satisfy the equalities $f_{k} =0$ for all $k \geq n+1$. We denote by $\LnA{n}$ and $\lnaft$ the full subcategories of $\Linf$ whose objects are Lie $n$-algebras, and finite type Lie $n$-algebras, respectively. As shown by the first author in \cite[Thm.\ 5.2]{R}, for a fixed $n$, both $\LnA{n}$ and $\lnaft$ form categories of fibrant objects (CFO), in the sense of K.\ Brown \cite{Brown}. The \df{weak equivalences} in these CFOs are the $L_\infty$-quasi-isomorphisms, and 
the \df{fibrations} are those $L_\infty$-morphisms $f \maps (L,\el) \fib (L',\el')$ such that $\tan(f)$ is a surjection in all positive degrees.

Recall that an $L_\infty$-algebra $(L,\el)$ is \df{minimal} if $\el_1=0$. The following relative version of minimal models will play a crucial role in the next subsection.
\begin{definition} \label{def:minfib}
A strict fibration $f = f_1 \maps (L,\el) \fib (L',\el') \in \LnA{n}$ is {\em minimal} if its kernel \eqref{eq:ideal} is a minimal Lie $n$-algebra.
\end{definition}

\subsection{Universal fibrations, Postnikov sections, and relative $k$-invariants} \label{sec:Linf-post}
As observed in \cite{Hen} and further developed in \cite[Sec.\ 7]{R}, every Lie $n$-algebra $(L, \el)$ admits a functorial Postnikov tower
\begin{equation} \label{eq:Ptower}
\cdots \tlt{m+1}L \to \tleq{m}L \to \tlt{m}L \to   \cdots \to  \tlt{1}L \to  \tleq{0}L = H_0(L) 
\end{equation}
whose zeroth section, or stage, is the the Lie 1-algebra $H_0(L)$, and with subsequent sections 
consisting of fibrations between suitably truncated Lie $n$-algebras. We recall the details of this construction below in Example \ref{ex:Ptower}. In particular, the morphisms
$\tlt{m+1}L \afib \tleq{m}L$ are acyclic fibrations, while the $\tleq{m}L \fib \tlt{m}L$ are minimal fibrations, in the sense of Def.\ \ref{def:minfib}. For this reason, we generalize and complete this theory by introducing the $L_\infty$-analog of relative $k$-invariants for an arbitrary minimal fibration $f = f_1 \maps (L,\el) \fib (L',\el') \in \LnA{n}$. 

These invariants are realized as explicit $L_\infty$-morphisms over $L'$ whose targets are certain universal fibrations. Each morphism classifies a corresponding section of a Postnikov tower of minimal fibrations over $L'$:
\begin{equation} \label{eq:rPtower}
\cdots \to \tleq{m}(L,f) \fib \tleq{m-1}(L,f) \fib  \cdots \fib \tleq{1}(L,f) \fib \tleq{0}(L,f).
\end{equation}
As we demonstrate below, the limit of the above tower in the category over $L'$ recovers $f$. Furthermore, we will see in Prop.\ \ref{prop:classify} that the $i^{th}$ section of the Postnikov tower for $f$ can be presented as the pullback over $L'$ of a universal fibration along the $(i-1)^{th}$ $k$-invariant. As a consequence, for any Lie $n$-algebra $L$,
the minimal fibrations $\tleq{m}L \fib \tlt{m}L$ appearing in its Postnikov tower \eqref{eq:Ptower} can be exhibited as a sequence of pullbacks of universal fibrations over $\tlt{m}L$. This result provides a key ingredient in the proof of our main theorem 
given in Sec.\ \ref{s:ind}.   
  
\newcommand{\relLnA}[1]{\LnA{\infty}_{/^{\mathrm{fib}}#1}}
In the following subsections, we fix a minimal fibration 
\[
f = f_1 \maps (L,\el) \fib (L',\el')
\] 
in $\LnA{\infty}$, and we work in the category $\relLnA{L'}$ whose objects are arbitrary fibrations over $L'$ in $\LnA{\infty}$, and whose morphisms are the usual commuting triangles. It follows from \cite[\Sec 4]{Brown} that the CFO structure on $\LnA{\infty}$ pulls back along the forgetful functor 
$\relLnA{L'} \to \LnA{\infty}$ to a CFO structure on $\relLnA{L'}$.

\subsubsection{Universal fibrations} \label{sec:uni-fib}
For $m \geq 1$ and a finite-dimensional non-graded vector space $A$, let $EA[m]$ denote the acyclic complex $A[m+1] \xto{\id_A} A[m]$. Using the notation in Example \ref{ex:semi-direct}, we define Lie $\infty$-algebras 
\[
\LEO{m}{L'}:= \bigl(EA[m]//\End(A) \bigr) \times L', \quad \LBO{m}{L'}:=\bigl(A[m+1]//\End(A) \bigr) \times L'. 
\]
The projections
\[
\End(A) \dsum A[m] \dsum A[m+1] \to \End(A), \qquad \End(A) \dsum A[m+1] \to \End(A), 
\]    
induce strict fibrations $\pi_E \maps \LEO{m}{L'} \to \End(A) \times L'$, and $\pi_B \maps \LBO{m}{L'} \to \End(A) \times L'$ in $\relLnA{L'}$. 
Furthermore, the projection
\begin{equation} \label{eq:uni-fib}
\End(A) \dsum A[m] \dsum A[m+1] \to \End(A) \dsum A[m],
\end{equation}
induces a minimal fibration $p_{A(m)/L'} \maps \LEO{m}{L'} \to \LBO{m}{L'}$. We call $p_{A(m)/L'}$ the \df{universal fibration} with fiber $A[m]$.

\subsubsection{Construction of Postnikov sections and $k$-invariants}\label{sec:tower}
We construct the Lie $\infty$-algebras $\bigl(\tleq{m}(L,f), \tleq{m}\el \bigr)$ which form the sections of Postnikov tower \eqref{eq:rPtower} for our fixed minimal fibration $f=f_1$. Let $m \geq 0$. As a chain complex, $\tleq{m}(L,f)$ is the graded vector space:
\[
\tleq{m}(L,f)_i = 
\begin{cases}
L_i & \text{if $i \leq m$,}\\
L'_i & \text{if $i >  m$.}\\
\end{cases}
\]
with the differential
\[
\tleq{m}^f\el_1(y) =
\begin{cases}
\el_1(y)& \text{if $\deg{y} \leq m$,} \\
\el_1(x) & \text{if $\deg{y} = m+1$, and $y=f_1(x)$,}\\ 
\el'_1(y) & \text{if $\deg{y} > m+1$.}
\end{cases}
\]
Note that $\tleq{m}\el_1$ is well defined in degree $m+1$ since $\ker f$ is a minimal Lie $\infty$-algebra, and $f$ is surjective in all positive degrees. Let 
$r_{\leq m} \maps (L,\el_1) \to (\tleq{m}(L,f), \tleq{m}\el_1)$
denote the surjective chain map 
\begin{equation}
r_{\leq m}(x) = 
\begin{cases}
x & \text{if $\deg{x} \leq m$}\\
f(x) & \text{if $\deg{x} > m$}.
\end{cases}
\end{equation}
The higher arity brackets on $\tleq{m}(L,f)$ are defined as:
\begin{equation} \label{eq:tleqm-brack}
\tleq{m}^f\el_k \bigl(r_{\leq m}(x_1), \ldots, r_{\leq m}(x_k):= r_{\leq m} \cc \el_k(x_1,\ldots,x_k),
\end{equation}  
which promote $r_{\leq m}$ to a strict fibration.  

By construction, the Lie $\infty$-algebras $\tleq{m}(L,f)$ assemble into a tower of fibrations under $L$ in $\relLnA{L'}$ 
\begin{equation} \label{diag:rPtowermap}
\begin{tikzdiag}{3}{3}
{
       \&          L      \&    \&     \\
\cdots \&  \tleq{m}(L,f) \& \tleq{m-1}(L,f) \& \cdots \\
       \&          L'      \&    \&     \\
};
\path[->>,font=\scriptsize]
(m-1-2) edge node[auto,swap] {$r_{\leq m}$} (m-2-2)
(m-1-2) edge node[auto] {$r_{\leq m-1}$} (m-2-3)
(m-2-1) edge node[auto] {$q^f_{\leq m +1}$} (m-2-2)
(m-2-2) edge node[auto] {$q^f_{\leq m}$} (m-2-3)
(m-2-3) edge node[auto] {$q^f_{\leq m-1}$} (m-2-4)
(m-2-2) edge node[auto,swap] {$\tleq{m}(f)$} (m-3-2)
(m-2-3) edge node[auto] {$\tleq{m-1}(f)$} (m-3-2)
;
\end{tikzdiag}
\end{equation}
Above $q^{f}_{\leq m} \maps \tleq{m}(L,f) \to \tleq{m-1}(L,f)$ and  $ \tleq{m}(f) \maps \tleq{m}(L,f) \to L'$ are the strict fibrations
\begin{equation} \label{eq:rPtowermap}
q^{f}_{\leq m} (x) = 
\begin{cases}
f(x) & \text{if $\deg{x} = m$} \\
x & \text{if $\deg{x} \neq m$} \\
\end{cases}
\qquad
\tleq{m}(f)(x) = 
\begin{cases}
f(x) & \text{if $\deg{x} \leq  m$} \\
x & \text{if $\deg{x} > m$} \\
\end{cases}
\end{equation}
Note that we have $\tleq{m}(f) \cc r_{\leq m} = f$. Furthermore, $q^f_{\leq m}$ is, in fact, a minimal fibration, surjective in all degrees, whose fiber is concentrated in degree $m$:
\begin{equation} \label{eq:kerq}
\ker q^f_{\leq m} = A[m], \qquad A:= (\ker f)_{m}.
\end{equation} 
We emphasize that $A$ is considered as a non-graded vector space, in agreement with the notation introduced in Sec.\ \ref{sec:uni-fib} above. It is clear that we recover the fibration $f$ from the tower as $m \xto{} \infty$:
\[
(L,\el) \xto{f} (L',\el')  \cong \plim_m  \tleq{m}(L,f) \quad \text{in $\relLnA{L'}$}
\]

Now we fix $m \geq 1$, and construct the $k$-invariant for the $m$th stage of \eqref{diag:rPtowermap}.  Let $\si_m \maps L'_m \to L_m$ denote a linear section of $f \maps L \to L'$ in degree $m$. This induces a section $\eta$ of $q^f_{\leq m} \maps \tleq{m}(L,f) \to \tleq{m-1}(L,f)$ at the level of graded vector spaces:
\begin{equation} \label{eq:eta}
\eta(x)=
\begin{cases}
\si_m(x) & \text{if $\deg{x} =m$},\\ 
x & \text{if $\deg{x} \neq m$}.
\end{cases}
\end{equation} 
and, hence, an isomorphism of graded vector spaces
\begin{equation} \label{eq:prA}
\phi \maps \tleq{m-1}(L,f) \dsum A[m] \xto{\cong} \tleq{m}(L,f).  
\end{equation}

Define for each $k \geq 1$, linear maps $\psi_k \maps \Alt^k \tleq{m-1}(L,f) \to \End(A) \dsum A[m+1]$
\begin{equation} \label{eq:class1}
\begin{split}
\psi_1(x)&:= 
\begin{cases}
\bs^{-m} \cc \tleq{m}^f\el_2 \bigl (x, \bs^{m}(-) \bigr) \in  \End(A) &
\text{if $\deg{x}=0$},\\
\bs \cc \pr_{A[m]} \cc  \phi^{-1}  \cc \tleq{m}^f\el_1(x) &
\text{if $\deg{x}=m+1$}\\
0 &\text{otherwise}
\end{cases}\\
\psi_k(x_1,\ldots, x_k)& := 
\begin{cases}
\bs \cc \pr_{A[m]} \cc  \phi^{-1}  \cc  \tleq{m}^f\el_k \bigl(\eta(x_1),\ldots, \eta(x_k) \bigr) & \text{if $\sum_{i} \deg{x_i} = m-k+2$},\\
0 &\text{otherwise}.
\end{cases}
\end{split}
\end{equation}

\begin{proposition} \label{prop:classify}
\mbox{}
\begin{enumerate}
\item The linear maps \eqref{eq:class1} define an $L_\infty$-morphism 
\[
\psi \maps \tleq{m-1}(L,f) \to A[m+1]//\End(A).
\]
\item There exists an $L_\infty$-morphism 
\[
\ti{\psi} \maps \tleq{m}(L,f) \to EA[m]//\End(A)
\]
inducing a pullback square 
\begin{equation} \label{diag:classify}
\begin{tikzdiag}{3}{5}
{
\tleq{m}(L,f)\&   \LEO{m}{L'} \\
\tleq{m-1}(L,f)\&   \LBO{m}{L'} \\
};
\path[->,font=\scriptsize]
(m-1-1) edge node[auto] {$\bigl(\ti{\psi},\tleq{m}(f)\bigr)$} (m-1-2)
(m-2-1) edge node[auto] {$\bigl(\psi,\tleq{m-1}(f)\bigr)$} (m-2-2)
;
\path[->>,font=\scriptsize]
(m-1-1) edge node[auto,swap] {$q^{f}_{\leq m}$} (m-2-1)
(m-1-2) edge node[auto] {$p_{A(m)/L'}$} (m-2-2)
;
\pbdiag[.3]
\end{tikzdiag}
\end{equation}
in $\relLnA{L'}$. 
\end{enumerate}
\end{proposition}

\begin{proof} 
Since $m \geq 1$, the minimal fibration $q^{f}_{\leq m}$ is surjective in all degrees, even if $f$ is not surjective in degree $0$. Hence, via the Chevalley-Eilenberg coalgebra functor, $q^{f}_{\leq m}$ corresponds to a coalgebra bundle, in the sense of Def.\ \ref{def:coalg-bun}. Furthermore, the Chevalley-Eilenberg coalgebra functor creates pullbacks by \cite[Prop.\ 4.1]{R}.
Therefore, the proof then follows by Quillen's classification theorem for principal dg coalgebra bundles \cite[\Sec B.5]{Quillen}, and a coalgebraic analog of the associated bundle construction developed by Prigge in \cite[\Sec 2]{Prigge}.
We recall the necessary details in Appendix \ref{app:Linf}, and defer the complete proof until Section 
\ref{sec:class-proof}.
\end{proof}

\begin{example}[Postnikov sections of Lie $\infty$-algebras]\label{ex:Ptower}
We discussed in the beginning of this section Henriques' construction of the Postnikov tower \eqref{eq:Ptower} for a Lie $\infty$-algebra $(L,\el)$. Here, following \cite[Sec.\ 7]{R}, we recall the details, which provide an important example of the classification of minimal fibrations given in Prop.\ \ref{prop:classify}. 
The sections of the tower for $(L,\el)$  
\begin{equation} \label{diag:Ptowermap}
\begin{tikzdiag}{3}{3}
{
       \&                \&  L   \& \&    \\
\cdots \&  \tlt{m+1}L \& \tleq{m}L \& \tlt{m}L \& \cdots \\
};
\path[->>,font=\scriptsize]
(m-1-3) edge node[auto,swap] {$p_{< m +1}$} (m-2-2)
(m-1-3) edge node[auto] {$p_{\leq m}$} (m-2-3)
(m-2-1) edge node[auto] {$q_{\leq m +1}$} (m-2-2)
(m-2-2) edge node[auto] {$q_{< m+1}$} node[below]{$\sim$} (m-2-3)
(m-2-3) edge node[auto] {$q_{\leq m}$} (m-2-4)
(m-2-4) edge node[auto] {$q_{< m}$} node[below]{$\sim$} (m-2-5)
;
\end{tikzdiag}
\end{equation}
lie under $L$, and consist of two flavors of truncations, in analogy with the Duskin and Moore truncations for smooth $\infty$-groups discussed earlier in Sec.\ \ref{s:post}. By construction, the tower recovers $(L,\el)$ as $m \to \infty$: 
\[
(L,\el) = \plim ( \cdots \to \tleq{1}L \to  \tlt{1}L \to  \tleq{0}L \bigr).
\]
For each $m \geq 0$, both $\tau_{\leq m}L$ and $\tau_{< m}L$ are Lie $(m+1)$-algebras with underlying complexes:
\begin{equation*} 
\begin{split}
(\tau_{\leq m} L)_i=
\begin{cases}
L_i & \text{if $i < m$,}\\
\coker(d_{m+1}) & \text{if $i=m$,}\\
0 & \text{if $i >m$,}
\end{cases}
\qquad 
(\tau_{< m} L)_i=
\begin{cases}
L_i & \text{if $i < m$,}\\
\im (d_{m}) & \text{if $i=m$,}\\
0 & \text{if $i >m$.}
\end{cases}
\end{split}
\end{equation*}
In degree $m$, the differentials for $\tau_{\leq m}L$ and $\tau_{<  m}L$ are $d_{m} \maps L_m/\im(d_{m+1}) \to L_{m-1}$, and the inclusion $ \im(d_m) \emb L_{m-1}$, respectively. Hence, the homology complexes of
$\tau_{\leq m}L$ and $\tau_{<  m}L$ are
\begin{equation} \label{eq:L-homology}
H_{i}(\tau_{\leq m}L) = 
\begin{cases}
H_i(L) & \text{if $i \leq m$,}\\
0 & \text{if $i>m$,}
\end{cases}
\qquad
H_{i}(\tau_{< m}L) = 
\begin{cases}
H_i(L) & \text{if $i <  m$,}\\
0 & \text{if $i \geq m$.}
\end{cases}
\end{equation}

The projections $p_{\leq m} \maps L \to \tau_{\leq m}L$ and $p_{< m} \maps L \to \tau_{< m}L$ are strict fibrations of Lie $\infty$-algebras. In particular, in degree $m$, $p_{\leq m}$ is the surjection $L_m \to \coker(d_{m+1})$, and $p_{< m}$ is the differential $d_{m} \maps L_m \to \im(d_{m})$. The higher arity brackets $\tleq{m}\el_k$ and $\tlt{m}\el_k$ are defined using the brackets on $L$, along with the surjections $p_{\leq m}$ and $p_{< m}$, respectively, as was done in the relative case \eqref{eq:tleqm-brack} above.    

The morphism $q_{< m+1} \maps \tau_{< m +1}L \to \tau_{\leq m}L$ is the strict acyclic fibration defined to be the projection $L_m \to \coker d_{m+1}$ in degree $m$, the identity in all degrees $<m$, and the zero map in degree $m+1$. By \cite[Lem.\ 7.4]{R}, acyclic fibrations in $\LnA{\infty}$ admit right inverses. Hence, for all $m\geq 0$, there exists $L_\infty$-morphisms $\si_{< m+1}$ splitting the $q_{< m+1}$:
\begin{equation} \label{eq:qlt-split}
\begin{tikzdiag}{2}{3}
{
\tau_{< m +1}L \& \tau_{\leq m}L\\
};
\path[->>,font=\scriptsize]
(m-1-1) edge node[auto] {$q_{<m+1}$} node[below]{$\sim$}  (m-1-2)
;
\path[->,dashed,font=\scriptsize]
(m-1-2) edge [bend right=45] node[auto,swap] {$\si_{< m+1}$} (m-1-1)
;
\end{tikzdiag}
\end{equation}
Next, the morphism 
\begin{equation} \label{eq:q-fib}
q_{\leq m} \maps \tleq{m}L \to \tlt{m}L
\end{equation}
appearing in tower \eqref{diag:Ptowermap} is the strict fibration given in degree $m$ by the differential $d_m \maps \coker{d_{m+1}} \to \im d_m$, and equal to the identity in all other degrees. We observe that Eq.\ \ref{eq:L-homology} implies that, for each $m \geq 0$, the fibration $q_{\leq m} \maps \tleq{m}L \to \tlt{m}L$ is minimal with fiber concentrated in degree $m$:
\begin{equation} \label{eq:Hm-fiber}
\ker q_{\leq m} = A[m], \quad A:= H_{m}(L)\
\end{equation}
It then follows from Prop.\ \ref{prop:classify} that $q_{\leq m}$ is a pullback in the category $\relLnA{\tlt{m}L}$.  

\end{example}

\subsection{Maurer-Cartan theory and integration} \label{sec:MC}
Let us first recall, following \cite[Sec.\ 8.6]{RZ}, the definition of Henriques' integration functor \cite{Hen} for finite-type Lie $\infty$-algebras using the language of Maurer-Cartan theory. Throughout, $S$ denotes a submanifold (possibly with corners) of $\R^N$. Fix an integer $r \geq 1$. As in \cite[Sec.\ 5.1]{Hen}, we denote by 
\[
\bigl(\Omega(S), d_{\dR} \bigr)
\] 
the differential graded Banach algebra
of \df{$C^r$-differentiable forms}. By definition, a $k$-form $\alpha$ on $S$ is an element of $\Omega(S)$ if and only if both $\alpha$ and the $(k+1)$-form $d_{\dR}\alpha$ are $r$-times continuously differentiable. 

Let $(L,\el) \in \lnaft$ be a finite type Lie $n$-algebra, and denote by 
\[
(L \tensor \OS, \el^{\Om}) 
\]
the $\Z$--graded $L_\infty$-algebra whose underlying \underline{chain complex} is $(L \tensor \OS , \el^{\Om}_1)$ where 
\[
\begin{split}
(L \tensor \OS)_m &:= \bigoplus_{i-j=m} L_{i} \tensor \OS^{j}\\
\el^{\Om}_1 &:= \el_1 \tensor \id_{\OS} + \id_{L} \tensor d_{\dR}
\end{split}
\]
and whose higher brackets are defined as:
\begin{equation} \label{eq:Ome-bracket}
\el^{\Om}_k \bigl(x_1 \tensor \omega_1, \ldots, x_k \tensor \omega_k \bigr):= (-1)^{\varepsilon} \el_k(x_1,\ldots,x_k) \tensor \omega_1\omega_2 \cdots \omega_k,
\end{equation}
with 
\[
\varepsilon :=  \sum_{1 \leq i < j \leq k} \deg{\omega_i}\deg{x_j}.
\]
Note that the \df{curvature}\footnote{The signs $\sgn{k}$ appear in the definition of curvature in order to stay consistent with the conventions used in \cite{RZ}. They simplify certain computations in the Chevalley-Eilenberg coalgebra. This forces an aesthetically unpleasant modification of the classical Maurer-Cartan equation, namely $dx - \frac{1}{2}[x,x] =0$, which agrees with the Maurer-Cartan equation \cite[p.\ 1032]{Hen} used by Henriques.}  
of an element $a \in (L \tensor \OS)_{-1}$:
\begin{equation} \label{eq:curv}
\begin{split}
\curv^{\Om}(a)  = \ell^\Om_1(a) + \sum_{k \geq 2} \sgn{k} \frac{1}{k!} \ell^\Om_k(a,a,\ldots,a) \in (L \tensor \OS)_{-2}
\end{split}
\end{equation}
is well-defined, i.e., the above summation is finite, since $\OS$ is bounded. Therefore, we may consider the set of 
\df{Maurer--Cartan elements} of the $L_\infty$-algebra $(L \tensor \OS, \el^{\Om})$:
\[
\MC \bigl(L \tensor \OS \bigr):= \bigl \{ a \in (L \tensor \OS)_{-1} ~\vert ~ \curv^\Om(a)=0 \bigr \}.
\]

Next, we recall the functorial properties of the Maurer-Cartan set. If $f \maps (L, \el) \to (L',\el')$ is a morphism in $\lnaft$, then the maps $f^{\Om}_k\maps \Lambda^{k} \bigl(L \tensor \OS \bigr) \to L' \tensor \OS$: 
\[
f^{\Om}_k\bigl(x_1 \tensor \omega_1, \ldots, x_k \tensor \omega_k \bigr):= (-1)^{\varepsilon} f_k(x_1,\ldots,x_k) \tensor \omega_1\omega_2 \cdots \omega_k,
\]
define a $L_\infty$-morphism $f^\Om \maps (L\tensor \OS, \el^\Om) \to (L' \tensor \OS, \el^{ \prime \Om})$.
Furthermore, $f^\Om$ induces a well defined function 
\[
f^\Om_\ast \maps (L \tensor \OS)_{-1} \to (L'\tensor \OS)_{-1} 
\]
where
\begin{equation} \label{eq:Fstar}
\begin{split}
f^\Om_\ast(a)  = f^\Om_1(a) + \sum_{k \geq 2} \sgn{k} \frac{1}{k!} f^\Om_k(a,a,\ldots,a).
\end{split}
\end{equation}
for all $a \in \bigl(L \tensor \OS \bigr)_{-1}$. Again, as above, the summation appearing in the 
formula for $f^\Om_\ast$ is finite. For reference later, we note that if $f=f_1 \maps (L, \el) \to (L',\el')$ is a strict morphism, then Eq.\ \ref{eq:Fstar} simplifies to
\begin{equation} \label{eq:strict_Fstar}
f^\Om_\ast = f \tensor \id_{\OS}. 
\end{equation} 
As demonstrated, for example, in \cite[Prop.\ 8.16]{RZ}, the assignments $(L\tensor \OS, \el^\Om) \mapsto \MC \bigl( (L\tensor \OS \bigr)$, and 
$f^\Om \mapsto f^{\Om}_\ast$ define a functor
\[
\MC(-\tensor \OS) \maps \lnaft \to \Set,
\]
natural in $S \subseteq \R^{N}$, which leads us to the definition of Henriques' integration functor:
\begin{pdef}[Def.\ 5.2, Thm.\ 5.10 \cite{Hen}] \label{pdef:int}
Let $L \in \LnA{n}^{\ft}$ be a finite type Lie $n$-algebra. 
The assignment
\[
L \mapsto \left(\sint L \right)_m:= \MC\bigl ( L \tensor \Omega(\Delta^m) \bigr)
\]
induces a functor
\begin{equation} \label{eq:int}
\sint \maps \LnA{n}^{\ft} \to \LnG{\infty} 
\end{equation}
from the category of finite type Lie $n$-algebras to the category of
Lie $\infty$-groups.
\end{pdef}

\begin{example}[Integrating Lie algebras] \label{ex:intnerve} 
  
  Let $\g$ be the Lie algebra of a simply connected Lie group $G$.  It follows from Prop-Def.\ \ref{pdef:int} that 
  \[
    (\sint \g)_n=\MC(\g\otimes\Omega^1(\Delta^n)),
  \]  
  i.e.\ the $n$-simplices of $\sint\g$ are by definition the space of flat connections on the trivial $G$-bundle on the $n$-simplex $\Delta^n$. Note that, by the convention established in Eq.\ \ref{eq:curv}, a flat connection $\theta$ is one which satisfies the Maurer-Cartan equation $d_{\dR}\theta - \frac{1}{2}[\theta,\theta]=0$.  As observed by Henriques \cite[Example 5.5]{Hen} (building on Sullivan \cite[Theorem 8.1]{Sull}), the map which sends a flat connection to its space of flat sections modulo translation gives a natural isomorphism 
  \begin{align}
    \MC(\g\otimes\Omega^1(\Delta^n))&\to^\cong \Map(\Delta^n,G)/G \label{e:flatsec},
  \end{align}
where $\Map(\Delta^m,G)$ denotes $G$-valued $C^{r+1}$ maps. Post-composing with the isomorphism
  \begin{align*}
    \Map(\Delta^n,G)/G&\to^\cong \Map_\ast(\Delta^n,G)\\
    \left[g\right]&\mapsto g(1)^{-1}\cdot g
  \end{align*}
  we obtain an isomorphism 
  \begin{align}
    (\sint\g)_n&\to^\cong \Map_\ast(\Delta^n,G) \label{e:sullbased} \\
    \theta&\mapsto g_\theta\nonumber    
  \end{align}
  where $g_\theta$ denotes the unique section satisfying 
\begin{equation} \label{eq:flat}
\tha = -g_\tha^{-1}d_{\dR}g_\tha,
\end{equation}
and sending $0\in\Delta^n$ to $e\in G$. 
  
  Under these isomorphisms, and suppressing the subscript $\theta$ for ease of reading, the simplicial structure on $\sint\g$ is given by
  \begin{align*}
    d_0g&:=g(1)^{-1}\cdot g\circ d^0\\
    d_ig&:=g\circ d^i\tag{for $i>0$}\\
    s_ig&:=g\circ s^i.
  \end{align*}
  Now consider the 1-truncation 
  \[
    \sint \g\to \tau_{\le 1}\sint\g.
  \]
  As observed in \cite[Example 7.2]{Hen}, the isomorphism~\ref{e:flatsec} induces an isomorphism
  \[
    (\tau_{\le 1}\sint\g)_n\cong \Map(\sk_0\Delta^n,G)/G.
  \]
  Similarly, the isomorphism~\ref{e:sullbased} induces an isomorphism
  \begin{align*}
    (\tau_{\le 1}\sint\g)_n&\to^\cong \Map_\ast(\sk_0\Delta^n,G)=G^n\\
    [g]&\mapsto (g(1),\ldots,g(n))
  \end{align*}
  which identifies $\tau_{\le 1}\sint\g$ with $NG$ with inhomogeneous coordinates. Tracing through the above, the 1-truncation of $\sint\g$ is isomorphic to the map
  \begin{align*}
    \sint\g & \to NG\intertext{given on $n$-simplices by}
    g&\mapsto (g(1),\ldots,g(n))
  \end{align*}
  with inhomogeneous coordinates on $NG$.
\end{example}

\begin{example}[Integrating abelian Lie {$n$}-algebras] \label{ex:ab-int}
Recall that an abelian Lie $n$-algebra in $\lnaft$ is finite-type chain complex $V \in \Chain$ concentrated in non-negative degrees. In this case, it follows directly from the definition that we have an equality of simplicial Banach spaces $\sint V = Z_{-1}(V \tensor \Om(\Del^\bl) )$, where the latter denotes the level-wise degree $-1$ cycles in the simplicial chain complex $V \tensor \Om(\Del^\bl)$. Via this identification, 
one can exhibit an injective morphism of simplicial Banach spaces 
\begin{equation} \label{eq:DK-inc}
\iota_V \maps K_\bl(V[1]) \hookrightarrow \sint V,  
\end{equation}
natural in $V$, where $K_\bl(-)$ is the Dold-Kan inverse to the normalized chains functor $N_\ast \maps s\mathsf{Vect} \to \Chain$.
Let us recall the construction. For each $m \geq 0$, let $C_{W}(\Del^m)$ denote the cochain complex of ``elementary Whitney forms'' on the $m$-simplex, in the sense of \cite[Sec.\ 3]{Get}. Every differential form in $C_{W}(\Del^m)$ has polynomial coefficients, hence we have an
inclusion of simplicial cochain complexes $C_W(\Del^\bl) \sse \Omega(\Del^\bl)$. As shown by Whitney \cite{Whit}, there is also an isomorphism of $C_{W}(\Del^\bl) \cong N^{\ast}(\Del^\bl)$, where $N^{\ast}(\Del^m) = \Hom(N_\ast(\Del^m),\R)$ denotes the normalized $\R$-valued cochains on the $m$-simplex. As a result, we obtain \eqref{eq:DK-inc} via the composition
\[
K_\bl(V[1]) = \hom_{\Chain}(N_\ast(\Del^\bl),V[1]) \cong Z_{-1}(V \tensor N^\ast(\Del^\bl) ) \hookrightarrow  Z_{-1}(V \tensor \Om^\ast(\Del^\bl) ) \cong
\sint V.  
\]     

\end{example}

\begin{example}[Integrating {$EA[n]$} and {$A[n+1]$}] \label{ex:A-int}
We focus here on two special cases of abelian $L_\infty$-algebras. Let $n \geq 1$ and $A$ a finite-dimensional vector space. Given a submanifold (with corners) $S \sse \R^{m+1}$,  denote by $(\Om^\ast(S ; A), d_{\dR})$ the de Rham complex\footnote{Here, by abuse of notation, $d_{\dR}$ is the trivial extension of the de Rham differential $d_{\dR}$ to $\Om^\ast(\Delta^m ; A)=A \tensor \Omega^\ast$ with no new Koszul signs.} of $C^r$-differential forms with values in the \und{non-graded} $\End(A)$-module $A$. As in Sec.\ \ref{sec:uni-fib}, let $EA[n]$ denote the chain complex $A[n+1] \xto {\id} A[n]$.
Then there are identifications \cite[Examples 5.3, 5.4]{Hen} for the integrations of $EA[n]$ and $A[n+1]$ in terms of the complex $\Omega(\Delta^m;A)$:  
\begin{align*}
  (\sint EA[n])_m&=\Omega^{n+1}(\Delta^m;A)\\
  (\sint A[n+1])_m&=\Omega^{n+2}_{\text{closed}}(\Delta^m;A).
\end{align*}  
Given $\mu \in \Omega^{n+1}(\Delta^m;A)$, the signed de Rham differential $(-1)^{n+1} d_{\dR} \mu \in \Omega^{n+2}(\Delta^m;A)$ defines a morphism of abelian Lie $(n+2)$-algebras, and hence a morphism of Lie $\infty$-groups 
\begin{equation} \label{eq:sign-dR}
\sint EA[n] \xto{(-1)^{n+1} d_{\dR}}  \sint A[n+1].    
\end{equation}
Moreover, by the functoriality of integration, this morphism is $\wti{\GL(A)}$-equivariant, with respect to the actions on source and target induced from the action on $A$, where $\wti{\GL(A)}$ denotes the 1-connected Lie group integrating $\End(A)$. As a result, \eqref{eq:sign-dR} induces a morphism between the homotopy quotients, in the sense of Def.\ \ref{d:hoquot}.  
\begin{equation} \label{eq:diag-hoq}
   \xymatrix{
     (\sint EA[n])//\widetilde{\GL(A}) \ar[rr] \ar[dr] && (\sint A[n+1])//\widetilde{\GL(A)} \ar[dl]\\
     & N\widetilde{\GL(A)}
   }
\end{equation}
\end{example}

\begin{example}[Integrating $\LE{n}$ and $\LB{n}$]\label{ex:LEA-int}
Keeping the notation and conventions of the previous example, let $\LE{n}$ and $\LB{n}$ denote the Lie $(n+2)$-algebras whose underlying chain complexes are $\LE{n} = \End(A) \dsum EA[n]$, and $\LB{n} = \End(A) \dsum A[n+1]$. Hence, in the notation of Sec.\ \ref{sec:uni-fib}, we are considering the domain and codomain of the universal fibration $p_{A(n)/L'} \maps \LEO{n}{L'} \to  \LBO{n}{L'}$ with trivial base $L'=0$. As we will now show, the $m$-simplices of $\sint L_{EA}$ and $\sint L_{BA}$ also have a convenient description in terms of the de Rham complex $(\Om^\ast(\Delta^m ; A), d_{\dR})$.

Let $(S,v_0)$ denote a pointed submanifold (with corners) of $\R^{m+1}$, and $\tha \in \MC(\End(A) \tensor \Om^1(S))$ a Maurer-Cartan element with $g_\tha \maps S \to \wti{\GL(A)}$ the unique $C^{r+1}$-function satisfying Eq.\ \ref{eq:flat}. Let $\omega \in \Om^{k}(S;A)$, be an $A$-valued $k$-form, and denote by $[\tha \wedge \omega] \in \Om^{k+1}(S;A)$ the action of $\tha(-) \wedge -$ on $\omega$.
A direct calculation \cite[Eqs.\ 30--32]{Hen} using Eq.\ \ref{eq:flat} shows that $\kappa \in A \tensor \Om^{k+1}(S;A)$ satisfies the equation         
$d_\dR \omega - [\tha \wedge \omega] = \ka$ if and only if 
\begin{equation} \label{eq:g-dot}
d_{\dR} (g_\tha \cdot \om) = g_\tha \cdot \ka,
\end{equation}
where $g_\tha \cdot -$ denotes the usual action of $\wti{\GL(A)}$-valued functions on $A$-valued differential forms. 

Now, suppose
\[
(\theta, \mu, \nu) \in (\LE{n} \tensor \Om(S))_{-1} = \End(A) \tensor \Om^1(S) \dsum A[n] \tensor \Om^{n+1}(S) \dsum A[n+1] \tensor \Om^{n+2}(S)  
\]
is a degree $-1$ element of the $L_\infty$-algebra $\LE{n}$. Then by combining Eq.\ \ref{eq:g-dot} with the formulas for the brackets \eqref{eq:Ome-bracket}, a direct calculation shows that $(\theta, \mu, \nu)$ is  Maurer-Cartan if and only if $\theta \in \MC(\End(A) \tensor \Omega^1(S))$ and the following equalities hold in the complex $\Omega(S;A)$: 
\begin{equation} \label{eq:MC-deRham}
d_{\dR}( g_\tha \cdot \nu) =0, \quad (-1)^{n+1}d_{\dR} (g_\tha \cdot \mu) = g_\tha \cdot \nu.
\end{equation}
Similarly, $(\theta, \nu) \in (\LB{n} \tensor \Om(S))_{-1} = \End(A) \tensor \Om^1(S) \dsum A[n+1] \tensor \Om^{n+2}(S)$ is a Maurer-Cartan element of $\LB{n}$ if and only if $\theta$ is Maurer-Cartan and $d_{\dR}( g_\tha \cdot \nu) =0$ in the complex $\Omega(S;A)$. As result of these identifications, we see that, for each $m$, the isomorphism \eqref{e:sullbased} extends to isomorphisms:
\begin{equation}\label{eq:EA-iso}
\begin{split}
  (\sint \LE{n})_m\to^\cong\{(g,\mu,\nu)\in &\Map_\ast(\Delta^m,\widetilde{\GL}(A))\times \Omega^{n+1}(\Delta^m;A)\times \Omega^{n+2}(\Delta^m;A)\\
  &\text{s.t.}~ (-1)^{n+1}d_{\dR} (g\cdot\mu)=-g\cdot\nu,~d_{\dR}(g\cdot\nu)=0.\}\\
(\theta,\mu,\nu)&\mapsto (g_\theta,\mu,\nu).
\end{split}
\end{equation}
and
\begin{equation}\label{eq:BA-iso}
\begin{split}
  (\sint \LB{n})_m& \to^\cong \{(g,\nu)\in \Map_\ast(\Delta^m,\widetilde{\GL}(A)) \times \Omega^{n+2} (\Delta^m;A)~|~d_{\dR}(g\cdot\nu)=0\}\\
(\theta,\nu)&\mapsto (g_\theta,\nu).
\end{split}
\end{equation}
With respect to these identifications, the universal fibration $p_{A(n)} \maps \LE{n} \to \LB{n}$ defined in \eqref{eq:uni-fib} is simply the projection
\[
  (g,\mu,\nu) \mapsto (g,\nu).
\] 
Furthermore, by comparing \eqref{eq:EA-iso} and \eqref{eq:BA-iso} with the homotopy quotients $(\sint EA[n])//\wti{\GL(A)}$ and $(\sint A[n+1])//\wti{\GL(A)}$ constructed in Example \ref{ex:A-int}, we see that the assignments
\begin{align*}
  (\theta,\mu,\nu)&\mapsto (g_\theta\cdot \mu,g_\theta(1),\ldots,g_\theta(m))\\
  (\theta,\nu)&\mapsto (g_\theta\cdot \nu,g_\theta(1),\ldots,g_\theta(m))\intertext{define simplicial maps}
  \sint \LE{n}&\to^\psi (\sint EA[n])//\widetilde{\GL(A)}\\
  \sint \LB{n}&\to^\varphi (\sint A[n+1])//\widetilde{\GL(A)}.
\end{align*}   
By combining these with the commuting diagram \eqref{eq:diag-hoq} induced by the signed de Rham differential \eqref{eq:sign-dR}, we obtain the pullback squares 
\[
  \xymatrix{
    \sint \LE{n} \ar[r]^-\psi \ar[d] & (\sint EA[n])//\widetilde{\GL(A)} \ar[d] \\
    \sint \LB{n} \ar[r]^-\varphi \ar[d] & (\sint A[n+1])//\widetilde{\GL(A)} \ar[d] \\
    \sint \End(A) \ar[r]^{\tau_{\le 1}} & N\widetilde{\GL(A)}
  }
\]
where we take inhomogeneous coordinates on $N\widetilde{\GL(A)}$ and the bottom map is the 1-truncation of $\int\End(A)$ as in  Example~\ref{ex:intnerve}.
\end{example}

\subsubsection{Integrating fibrations} \label{sec:int-fib}
The compatibility of the functor $\sint \maps \lnaft \to \LnG{\infty}$ with the CFO structure on $\lnaft$ and the iCFO structure on Lie $\infty$-groupoids was characterized in \cite[Sec.\ 9]{RZ}.
We recall the key results. Firstly, a fibration $ f \maps (L,\el) \to (L',\el')$, in the sense of Sec.\ \ref{sec:LnA}, is called \df{quasi-split} if the induced map in homology $H(f_1) \maps H(L) \to H(L')$ is surjective in all degrees, and $H_0(L) \cong \ker H_0(f_1) \oplus H_0(L')$ as Lie algebras. 
\begin{example}[Quasi-split fibrations]\label{ex:q-split}
Every acyclic fibration is quasi-split, and so are the fibrations $q_{\leq m} \maps \tleq{m} L \fib \tlt{m}L$ appearing in the Postnikov tower \eqref{diag:Ptowermap} of any Lie $n$-algebra. Finally, as in Sec.\ \ref{sec:uni-fib}, for $n \geq 1$ and a finite-dimensional vector space $A$, the strict fibrations $\pi_E \maps \LE{n} \to \End(A)$ and $\pi_B \maps \LB{n} \to \End(A)$ are both quasi-split.
In contrast, the universal fibration \eqref{eq:uni-fib} $p_{A(n)/L'} \maps \LEO{n}{L'} \xto{} \LBO{n}{L'}$ over any base $L'$ is not quasi-split. (However, see Prop.\ \ref{prop:uni-fib-int} below.) 
\end{example}
As shown in \cite{RZ}, the functor $\sint(-)$ behaves like an exact functor between CFOs when restricted to quasi-split fibrations. More precisely:

\begin{proposition}[Thm.\ 9.16 \cite{RZ}] \label{prop:int-exact}
The integration functor $\sint \maps \LnA{n}^{\ft} \to \LnG{\infty}$
preserves finite products, weak equivalences, acyclic fibrations, and
sends quasi-split fibrations to Kan fibrations. Furthermore, $\sint(-)$ preserves pullback squares in $\lnaft$ of the form
\[
\begin{tikzpicture}[descr/.style={fill=white,inner sep=2.5pt},baseline=(current  bounding  box.center)]
\matrix (m) [matrix of math nodes, row sep=2em,column sep=3em,
  ampersand replacement=\&]
  {  
(\ti{L},\ti{\el}) \& (L,\el) \\
(L',\el') \& (L'',\el'') \\
}
; 
  \path[->,font=\scriptsize] 
   (m-1-1) edge node[auto] {$$} (m-1-2)
   (m-1-1) edge node[auto,swap] {$$} (m-2-1)
   (m-1-2) edge node[auto] {$f$} (m-2-2)
   (m-2-1) edge node[auto] {$$} (m-2-2)
  ;

  \begin{scope}[shift=($(m-1-1)!.4!(m-2-2)$)]
  \draw +(-0.25,0) -- +(0,0)  -- +(0,0.25);
  \end{scope}
\end{tikzpicture}
\]
whenever $f$ is a fibration.
\end{proposition}      

The class of fibrations mapped to Kan fibrations by the integration functor turns out to be strictly larger than the quasi-split ones. Indeed, as we will now show, every universal fibration $p_{A(m)/L'} \maps \LEO{n}{L'} \to \LBO{n}{L'}$ integrates to a Kan fibration.
\begin{proposition}\label{prop:uni-fib-int}
For all $n \geq 1$ and each finite-dimensional vector space $A$, the universal fibration \eqref{eq:uni-fib} over the base $L'$ integrates to a Kan fibration 
\[
\sint p_{A(n)/L'} \maps \sint \LEO{n}{L'} \to \sint \LBO{n}{L'}
\]
between Lie $\infty$-groups over $\sint L'$.
\end{proposition}
The proof of Prop.\ \ref{prop:uni-fib-int} will require a trivially modified version of a result from \cite{Hen}. In what follows, we adopt the notation and results presented in Example \ref{ex:LEA-int}.

\begin{lemma}[Lem.\ 5.11 \cite{Hen}]\label{lem:prim}
Let $m \geq 1$, $k=0,\ldots,m$, and let $\jm \maps \Lam^m_k \to \Del^m$ denote the inclusion.  
Let $P \in \Mfd$ denote the pullback of the diagram
\[
\Om^{n+1}(\Lam^m_k;A) \xto{(-1)^{n+1}d_{\dR}} \Om^{n+2}_{\cl}(\Lam^{m}_k;A) \xleftarrow{\jm^\ast} \Om^{n+2}_{\cl}(\Del^m;A).
\]
Then the natural map $(\jm^\ast, d_{\dR}) \maps \Om^{n+1}(\Del^m,A) \to P$ is a surjective submersion.
\end{lemma}  

\begin{proof}[Proof of Prop.\ \ref{prop:uni-fib-int}]
It suffices to prove the statement for the universal fibration over the trivial base  
$p_{A(n)} \maps \LE{n}\to \LB{n}$, since $p_{A(n)/L'} = p_{A(n)} \times \id_{L'}$ and integration preserves products. For brevity, we denote by $\ph \maps X \to Y$ and $\pi \maps Y \to Z$ the  
morphisms of Lie $\infty$-groups $\sint p_{A(n)} \maps \sint \LE{n} \to \sint \LB{n}$ and  
$\sint \pi_B \maps \sint \LB{n} \to \sint \End(A)$, respectively. Since $\pi_B$ is a quasi-split fibration of Lie $(n+2)$-algebras,
Prop.\ \ref{prop:int-exact} implies that $\pi$ is a Kan fibration. 

Let $m \geq 1$, $k=0,\ldots,m$, and let $\jm \maps \Lam^m_k \to \Del^m$ denote the inclusion.
We will show that the map
\begin{equation} \label{eq:ufi1}
X_m \xto{(\ph, \jm^\ast)} Y_m \times_{\hrn{Y}} \hrn{X}
\end{equation}
is a surjective submersion. Let $P_1$ denote the pullback of the diagram
\[
Y_m \xto{(\pi,\jmm)}Z_m \times_{\hrn{Z}} \hrn{Y} \xleftarrow{\pi \times \ph \vert_{\horn}} 
Y_m \times_{\hrn{Y}} \hrn{X}.
\]
Since $\pi$ is a Kan fibration, the above map $Y_m \to Z_m \times_{\hrn{Z}} \hrn{Y}$ is a surjective submersion. Hence, $P_1$ is a manifold and the induced map
\[
\pi^{P_1} \maps P_1 \to Y_m \times_{\hrn{Y}} \hrn{X}
\]
is a surjective submersion. It follows from the identifications \eqref{eq:EA-iso} and \eqref{eq:BA-iso}, along with Eq.\ \ref{eq:MC-deRham}, that an element $( (\Ga,\ro), \ro', (\Ups ,\nu,\eta)) \in P_1$ consists of functions $\Ga \maps \Del^m \to \wti{\GL(A)}$ and $\Ups \maps \horn \to \wti{\GL(A)}$, along with $A$-valued differential forms
\[
\ro,\ro' \in \Om^{n+2}(\Del^m;A), \quad \lam \in \Om^{n+1}(\horn;A), \quad  \eta \in \Om^{n+2}(\horn;A)
\]
which satisfy the following equalities:
\begin{align} \label{eq:ufi2a}
& \Ga \vert_{\horn} = \Ups &  & \eta = \jmm \ro = \jmm \ro'\\ \label{eq:eq:ufi2b}
&d_{\dR}(\Ga \cdot \ro) = d_{\dR}(\Ga \cdot \ro')=0 &  & d_{\dR}(\Ups \cdot \eta) =0
\\ \label{eq:ufi2c}
&  (-1)^{n+1}d_{\dR}(\Ups \cdot \lam) = \Ups \cdot \eta.  & & 
\end{align}    
Next, let $P$ denote the pullback of the diagram in Lemma \ref{lem:prim}. Then Eqs.\ \ref{eq:ufi2a} -- \ref{eq:ufi2c} imply that the function
$\psi \maps P_1 \to P$
\[
\psi((\Ga,\ro), \ro', (\Ups ,\lam,\eta)):=( \Ups \cdot \lam, \Ga \cdot \ro)
\]  
is a well-defined map between manifolds. Let $P_2$ denote the pullback of the diagram
\[
P_1 \xto{\psi} P \leftarrow  \Om^n(\Del^m;A) 
\]
where the map on the right-hand side is the surjective submersion $(\jm^\ast, (-1)^{n+1}d_{\dR})$ from Lem.\ \ref{lem:prim}. Hence, $P_2$ is a manifold, and the induced map
\[
\pi^{P_2} \maps P_2 \to P
\] 
is a surjective submersion. An element in $P_2$ consists of an element
$( (\Ga,\ro), \ro', (\Ups ,\lam,\eta))$ of  $P_1$ along with a $A$-valued $n+1$-form $\ti{\al} \in \Om^{n+1}(\Del^m;A)$ satisfying
\[
\jmm \ti{\al} = \Ups \cdot \lam, \quad (-1)^{n+1}d_\dR \ti{\al} =  \Ga \cdot \ro.
\] 
Finally, recall that, via the isomorphism \eqref{eq:EA-iso}, an element $(\Ga, \al,\be) \in X_m$ consists of
a function $\Ga \maps \Del^m \to \wti{\GL(A)}$ as above, and $A$-valued differential forms
$\al \in \Om^{n+1}(\Del^m;A)$, $\be \in \Om^{n+2}(\Del^m;A)$ satisfying
\[
(-1)^{n+1}d_{\dR}(\Ga \cdot \al) =  \Ga \cdot \be, \quad d_{\dR}(\Ga \cdot \be) = 0. 
\]
Hence, the function $\chi \maps P_2 \to X_m$
\[
\chi\bigl( (\Ga,\ro), \ro', (\Ups ,\lam,\eta), \ti{\al} \bigr):= (\Ga, \Ga^{-1} \cdot \ti{\al},\ro)
\] 
is a well-defined map between manifolds. It is easy to see that $\chi$ is onto, and
its composition with the map  $(\ph, \jm^\ast) \maps X_m  \to  Y_m \times_{\hrn{Y}} \hrn{X}$ 
from Eq.\ \ref{eq:ufi1} is
\[
(\ph, \jm^\ast) \cc \chi = \pi^{P_2} \cc \pi^{P_1}.
\]
Since the right-hand side of the above equality is a surjective submersion, the surjectivity of $\chi$ implies (e.g.\ \cite[Lem.\ 6.8]{RZ}) that $(\ph, \jm^\ast)$ is a surjective submersion as well. 
\end{proof}

We conclude this section with a lemma which will play a critical role in the induction step of our proof of Thm.\ \ref{t:main3}.

\begin{lemma} \label{l:postsplit}
  Let $L \in \Linfft$. Let $n>0$. Consider the fibration of Lie $(n+1)$-algebras \eqref{eq:q-fib}  $q_{\leq n}\colon \tau_{\le n} L \to \tau_{<n} L$.  Then the canonical morphism \ref{eq:dusk-diag}
  \[
    q \maps \sint\tau_{\le n} L \to \tau_{\le n+1}\left(\sint\tau_{\le n}L ,\sint  q_{\leq n}\right)
  \]
    admits a section. 
\end{lemma}
\begin{proof}
    Let $A=H_nL$. By Prop.\ \ref{prop:classify} and Example \ref{ex:Ptower}, the map $\tau_{\le n}L \to \tau_{<n}L$ fits into in a  pullback square
    \[
        \xymatrix{
         \tau_{\le n}L \ar[d]_{q_{\leq n}} \ar[r] & \LEO{n}{\tlt{n}L} \ar[d]^{p_{A(n)/ \tlt{n}L}} \\
         \tau_{<n}L \ar[r] & \LBO{n}{\tlt{n}L}
        }
    \]
    By Prop. \ref{prop:uni-fib-int} and Prop.\ \ref{prop:int-exact}, integrating this square gives a pullback square
    \[
        \xymatrix{
         \sint \tau_{\le n} L  \ar[d]_{\sint p} \ar[r] & \sint \LEO{n}{\tlt{n}L} \ar[d]^{\sint p_{A(n)/ \tlt{n}L} } \\
         \sint \tau_{<n} L  \ar[r] & \sint \LBO{n}{\tlt{n}L}
        }
    \]
    in which both vertical maps are Kan fibrations. Since integration preserves products, we have 
$\sint p_{A(n)/ \tlt{n}L} = \sint (p_{A(n)} \times \id_{\tlt{n}L }) = (\sint p_{A(n)}) \times \id_{\int \tlt{n}L}$. Therefore, to show the lemma, it suffices to show it in the universal case over the trivial base, i.e. the map
    \[
        \sint \LE{n} \to \tau_{\le n+1}(\sint \LE{n},\sint p_{A(n)})
    \]
    splits. Indeed, as we saw in Example \ref{ex:LEA-int}, this map fits into a pullback square
   \[
        \xymatrix{
         \sint \LE{n} \ar[r]^-\psi \ar[d] & (\sint EA[n])//\widetilde{\GL(A)} \ar[d]^f \\
         \sint \LB{n} \ar[r]^-\varphi & (\sint A[n+1])//\widetilde{\GL(A)} 
        }
   \]
   Therefore, it suffices to show that 
   \[
        (\sint EA[n])//\widetilde{\GL(A)}\to\tau_{\le n+1}((\sint EA[n])//\widetilde{\GL(A)},f)
   \]
   splits.  But this follows from Example \ref{ex:ab-int}.  Indeed, $EA[n]\to A[n+1]$ is a map of abelian $L_\infty$-algebras, and the inclusion \eqref{eq:DK-inc}  $K_\bl (-[1]) \xto{\iota} \sint(-)$ is natural. Specializing to the present case, this determines an inclusion 
   \[
        \xymatrix{
            WK(A,n+1)//\widetilde{\GL(A)} \ar[rr]^\sigma \ar[d] && (\sint EA[n])//\widetilde{\GL(A)} \ar[d]^f \\
            K(A,n+2)//\widetilde{\GL(A)} \ar[rr] && (\sint A[n+1])//\widetilde{\GL(A)}
        }
   \]
    By inspection, the relative $(n+2)$-truncation of the map
    \[ 
        (\sint A[n+1])//\widetilde{\GL(A)}\to N\widetilde{\GL(A)}
    \]
    gives a retract, over $N\widetilde{\GL(A)}$ of the bottom map.  Taking the relative $(n+1)$-truncation of $f$, we obtain an isomorphism
    \begin{align*}
        \tau_{\le n+1}&((\sint EA[n])//\widetilde{\GL(A)},f)\\
        &\cong (\sint A[n+1])//\widetilde{\GL(A)}\times_{K(A,n+2)//\widetilde{\GL(A)}}WK(A,n+1)//\widetilde{\GL(A)}. 
    \end{align*}
    Letting $\pi_2$ denote the projection of the fiber product onto the second factor, we see that the composition 
    \begin{align*}
        \tau_{\le n+1}&((\sint EA[n])//\widetilde{\GL(A)},f)\\
        &\cong (\sint A[n+1])//\widetilde{\GL(A)}\times_{K(A,n+2)//\widetilde{\GL(A)}}WK(A,n+1)//\widetilde{\GL(A)}\\
        &\to^{\pi_2} WK(A,n+1)//\widetilde{\GL(A)}\to^\sigma (\sint EA[n])//\widetilde{\GL(A)}
    \end{align*}
    gives the desired section.
\end{proof}

\begin{remark}
    Lemma~\ref{l:postsplit} is false if we allow $n=0$.  We can see this for $\g$ a Lie algebra, so that $\tau_{\le 0}\g=\g$ and $\tau_{<0}\g=0$, and we are considering the map
    \[
        \sint\g\to \tau_{\le 1}\sint\g=NG.
    \]
    This map does not admit a splitting in general. To see this, let $\Top$ denote the category of (weak Hausdorff compactly generated) topological spaces with the usual model structure.  Consider the geometric realization
    \begin{align*}
        |-|\colon \sSh(\sm)&\to\Top\\
                F&\mapsto \colim_{\Delta^k\times U\to F} |\Delta^k|\times U
    \end{align*}
    Unwinding the proof of \cite[Proposition 2.3]{Dugger}, we see that this is a left Quillen functor from the (global) presheaf model structure on simplicial presheaves to spaces.  In particular, it preserves global weak equivalences.  Therefore, as discussed in \cite[Sections 6,8]{Hen}, the map
    \[
        |\sint\g|\to |\tau_{\le 1}\sint\g|
    \]
    is weakly equivalent to the map
    \[
        \ast\to BG.
    \]
    We see that the non-contractibility of $G$ is an obstruction to the lemma in the $n=0$ case.
\end{remark}

\section{Lie $\infty$-groups}\label{s:liegroup}
In this section and for the rest of the paper, we specialize to the case $\C=\sm$ of smooth Banach manifolds with surjective submersions as covers as in \cite{Hen,RZ}.  Our terminology and notation largely follows \cite{Hen,W,BG,RZ}. Note that Banach manifolds with surjective submersions as covers satisfy all the axioms considered in \cite{W} for a category with covers, including \cite[Axiom 5]{W} Godemont's Theorem (see \cite[5.9.5]{VAR}). We assemble some basic results we need here.

By a {\em simplicial manifold}, we mean a diagram
\[
    X\colon \Delta^{\op}\to\sm.
\]
By a {\em finite dimensional} simplicial manifold, we mean that the manifold of $n$-simplices $X_n$ is finite dimensional for all $n$.

\begin{definition}\label{d:good}\mbox{}
    \begin{enumerate}
       \item A map $U\to X$ is {\em reduced} if $U_0=X_0=\ast$.
        \item A simplicial manifold $X\in\ssm$ is {\em good} if for all $n$, $X_n\simeq \pi_0 X_n$, i.e.  $X_n$ is a disjoint union of contractible pieces for all $n$.  It is {\em good up to level $n$}, if for all $k\le n$, $X_k\simeq \pi_0 X_k$.  Similarly, a hypercover $U\to X$ is {\em good} (resp. {\em good up to level $n$}) if the total space $U$ is a good simplicial manifold (resp. is good up to level $n$).
    \end{enumerate}
\end{definition}

\begin{example}
	Let $X$ be a manifold, viewed as a constant simplicial manifold.  Let $U\to X$ be a cover. Then $U$ is a {\em good cover} in the ordinary sense if and only if its nerve $NU\to X$ (as in~\eqref{e:nerve}) is a good hypercover.
\end{example}

\begin{definition}
    A {\em Lie} $\infty$-groupoid is a smooth $\infty$-groupoid in $\sm$.  Similarly, by a {\em Lie} $\infty$-group, $n$-group/oid, etc. we mean a smooth $\infty$-group, $n$-group/oid, etc. in $\sm$.
\end{definition}

\subsection{Descent categories}
We will also want to consider the analogue of hypercovers, $n$-stacks, etc. with surjective local diffeomorphisms in lieu of surjective submersions. Unfortunately, the category of Banach manifolds with surjective local diffeomorphisms fails to satisfy Axiom~\ref{a:term} for a category with covers: for $\dim X>0$, $X\to\ast$ is not a local diffeomorphism!  We can repair this by working in a closely related formalism introduced by Behrend and Getzler in \cite{BG}.  

Recall that a (sub-canonical) {\em descent category} consists of a category $\D$ with a subcategory of covers satisfying the following axioms:
\begin{enumerate}
    \item\label{a:dfinlim} $\D$ has all finite limits,
    \item\label{a:dterm} The map $\ast\to\ast$ is a cover, where $\ast\in\D$ denotes the terminal object.\footnote{This condition was unintentionally omitted in \cite{BG}.}
    \item\label{a:dbasechange} pullbacks of covers are covers,
    \item\label{a:drightcancellation} if $f$ and $gf$ are covers, then so is $g$,
    \item\label{a:dsubcan} covers are effective epimorphisms.
\end{enumerate}
Descent categories provide a setting for working with possibly ``singular'' objects in geometric contexts, with the covers being ``smooth'' morphisms.  Relevant examples include:
\begin{enumerate}
    \item The category of schemes with surjective \'etale, smooth or flat morphisms as the covers,
    \item The category of analytic spaces with surjective submersions as covers,
    \item The category of Banach analytic spaces with surjective submersions as covers,
    \item $\D=\Sh(\C)$, with $\C$ a category with covers, and with the local epimorphisms as covers.
\end{enumerate}

\begin{remark}
    Note that the two key differences between a category with covers and a descent category are that a) in a category with covers $\C$, the map $X\to\ast$ is a cover for all $X\in \C$, and b) a descent category $\D$ has all finite limits. Therefore, while smooth manifolds and surjective submersions form a category with covers, the failure of limits to exist prevents them from forming a descent category.  If we work instead with $C^\infty$-schemes and surjective submersions, then these form a descent category, but for singular $C^\infty$-schemes $X$, the map $X\to\ast$ will not be a cover.
\end{remark}

\begin{remark}\label{r:desc<=>catcov}
    The theory of (sub-canonical) descent categories closely parallels that of categories with covers, and vice versa. One can see this, for example, by comparing \cite{BG} and \cite{W,W2}. More precisely, \cite{W2} is explicitly written to apply to both categories with covers and descent categories, while a careful reading of \cite{W} makes clear that all of the results of \cite[Sections 2, 3, 4, and 5]{W} rely only on the axioms shared in common by categories with covers and descent categories; therefore, they apply, by the same proofs, to simplicial objects in a descent category $\D$. As a result, with the sole exception of Lemma~\ref{l:prods}, all of the results in Section~\ref{s:prelim} and Sections~\ref{s:post1} and~\ref{s:min} carry over unchanged, via the same proofs, to the setting of simplicial objects in descent categories. 
    
    Further, the results of \cite[Section 6]{W}, and their proofs, carry over unchanged to descent categories, so long as we explicitly require that for a simplicial group $G\in\sD$, the maps $G_n\to\ast$ must be covers for all $n$.\footnote{N.b. this assumption is only relevant for \cite[Proposition 6.6 and Theorem 6.7(2)]{W}.} With this assumption added, the results of Section~\ref{sec:smooth-kinv} carry over as well.
\end{remark}

We record the following elementary consequence of the definitions.
\begin{lemma}\label{l:desccore}
    Given a descent category $\D$, let $\C\subset\D$ be the full subcategory consisting of objects $X$ for which the map $X\to \ast$ is a cover, and define a map in $\C$ to be a cover if it is a cover in $\D$. Then $\C$ is a category with covers. 
\end{lemma}
\begin{proof}
    Axiom~\ref{a:dterm} for a descent category implies that $\C$ contains the terminal object of $\D$. The definition of $\C$ ensures that $\ast\in\C$ is terminal and that Axiom~\ref{a:term} for a category with covers. That Axioms~\ref{a:rightcancellation} and~\ref{a:subcan} for categories with covers hold for $\C$ follows immediately from the corresponding axioms for descent categories. For Axiom~\ref{a:basechange}, it suffices to observe that if $X\to Y$ is a cover in $\C$, and $Z\to Y$ is any map in $\C$, then $Z\times_Y X\to Z\to \ast$ is a composition of covers in $\D$ (by Axioms~\ref{a:dfinlim} and ~\ref{a:dbasechange}).  Therefore $Z\times_Y X\in\C$, Axiom~\ref{a:basechange} for categories with covers holds in $\C$.
\end{proof}

We now introduce the main descent category we work in.
\begin{definition}
    Let  $f\colon F\to G$ be a map in $\Sh(\sm)$. We say that $f$ a {\em local diffeomorphism} (submersion) if for all sections $U\to G$ where $U\in\sm$, the pullback 
\[
    F\times_G U\to U
\]
is a local diffeomorphism (submersion) of smooth manifolds. We say it is a surjective local diffeomorphism (surjective submersion) if it is a surjective map of sheaves. 
\end{definition} 

\begin{remark}
    As an immediate consequence of the definition, note that if $F\to U$ is a local diffeomorphism or submersion of sheaves with $U$ representable, then $F\in\sm$ as well.
\end{remark}

\begin{remark}
    The notion of surjective submersion of sheaves on smooth manifolds here is strictly stronger than that appearing in \cite[Definition 7.4]{Hen}. For example, let $G$ be a finite dimensional Lie group, $\Gamma$ a finitely generated discrete group, and $\rho\colon\Gamma\to G$ a homomorphism with nondiscrete image.  Then $G/\Gamma\to \ast$ satisfies the conditions of \cite[Definition 7.4]{Hen}, but fails to be a submersion in our sense, e.g. because $G/\Gamma$ is not a manifold.  This example also shows that surjective submersions of sheaves, in the above sense, fail to have the right cancellation property (Axiom~\ref{a:drightcancellation}): the map $G\to G/\Gamma$ is a surjective submersion of sheaves, as is the composite $G\to G/\Gamma\to\ast$, but $G/\Gamma\to\ast$ is not. 
\end{remark}

\begin{lemma}\label{l:desccat}
    The category $\Sh(\sm)$ with surjective local diffeomorphisms as covers forms a descent category.
\end{lemma}
\begin{proof}
    The only nontrivial statement to check is the right cancellation property (Axiom~\ref{a:drightcancellation}).  Let $F\to^f G\to^g H$ be a pair of maps of sheaves such that $f$ and $gf$ are surjective local diffeomorphisms. We must prove that $g$ is one as well.  Let $U\to H$ be any map with $U\in\sm$. It suffices to prove that $G\times_H U$ is a smooth manifold.  Indeed, consider the following diagram of pullback squares
    \[
        \xymatrix{
            F\times_H U\ar[r] \ar[d] & G\times_H U\ar[r] \ar[d] & U\ar[d]\\
            F\ar[r]^f & G\ar[r]^g & H
        }.
    \]
    Because $gf$ is a surjective local diffeomorphism, $F\times_H U$ is a manifold and $F\times_H U\to U$ is a surjective local diffeomorphism. Because local diffeomorphisms of smooth manifolds satisfy right cancellation, if $G\times_H U$ is representable, it immediately follows that $G\times_HU\to U$ is a surjective local diffeomorphism.

    Recall that sheaves are colimits of representables, i.e. we have an isomorphism of sheaves
    \[
        \varinjlim_{V\to G\times_H U} V\cong G\times_H U.
    \]
    where the colimit is indexed by the category of elements of $G\times_H U$. Because $F\times_HU\to U$ is a surjective local diffeomorphism of smooth manifolds, $F\times_H U$ admits an atlas $\mathcal{A}$ consisting of open balls which map diffeomorphically onto their image in $U$. Because $F\times_HU\to G\times_HU$ is a surjective local diffeomorphism of sheaves, the collection of maps
    \[
        V\to G\times_H U,
    \]
    where $V\subset F\times_H U$ is an open ball in the atlas $\mathcal{A}$, forms a cofinal subcategory of the category of elements. Moreover, because $F\times_HU\to G\times_HU$ is a local diffeomorphism of sheaves, and because every such $V$ maps diffeomorphically onto its image in $U$, we see that  every such $V$ has the property that the map $V\to G\times_HU$ is an open embedding, in the sense that the pullback along along any map $W\to G\times_HU$, with $W$ a smooth manifold, is an open embedding. We conclude that the collection of all such $V\to G\times_HU$ gives a smooth atlas of $G\times_HU$, and thus that $G\times_HU$ is representable as claimed.
\end{proof}

\begin{notation}\label{n:etloc}
    Below, we will at times want to treat $\Sh(\sm)$ with surjective local diffeomorphisms as a descent category, and at others to treat $\Sh(\sm)$ with local epimorphisms as a category with covers. To attempt to reduce confusion, we denote the descent category by $\Sh(\sm)_{\et}$, and the category with covers by $\Sh(\sm)_{\loc}$.
\end{notation}

\begin{remark}\label{r:discgroupkinvar}
    Let $G\in\sSet$ be a discrete simplicial group. Viewing sets as 0-dimensional manifolds, we see that for all $n$, the map $G_n\to\ast$ is a cover in the descent category $\Sh(\sm)_{\et}$. As a result, as observed in Remark~\ref{r:desc<=>catcov}, if $f\colon X\to Y$ is a map of simplicial sheaves satisfying the assumptions of Theorem~\ref{t:kinvar} and $\pi_n F$ is a discrete group, then the conclusions of Theorem~\ref{t:kinvar} apply to $f\colon X\to Y$, by the same proof.  Similar comments apply to the other results in Section~\ref{sec:smooth-kinv}.
\end{remark}

\subsection{\'Etale hypercovers and fibrations.}
A primary motivation for introducing the descent category $\Sh(\sm)_{\et}$ is to treat the following classes of morphisms.
\begin{definition} \label{def:etale-maps} 
\mbox{} 
    \begin{enumerate}
        \item  An {\em \'etale hypercover} is a hypercover in $\Sh(\sm)_{\et}$, i.e. a map $f\colon U\to X$ of simplicial sheaves such that, for all $n\ge 0$, the boundary-filling maps
	       \[
		          \mu_n(f)\colon U_n\to M_n(f)
	       \]
	       are surjective local diffeomorphisms.
        \item An {\em \'etale Kan fibration} is a Kan fibration in $\Sh(\sm)_{\et}$, i.e. a map $f\colon X\to Y$ of simplicial sheaves such that for all $n>0$ and all $i$, the horn-filling maps 
        \[
            \lambda^n_i(f)\colon X_n\to \Lambda^n_i(f)
        \]
        are surjective local diffeomorphisms. We say $f$ is an {\em \'etale covering Kan fibration} if in addition, $f_0\colon X_0\to Y_0$ is a surjective local diffeomorphism.
    \end{enumerate}
\end{definition}

We have the following analogues of Lemmas~\ref{l:covkan} and \ref{l:moorehyp}. As observed in Remark~\ref{r:desc<=>catcov}, these statements hold via the same proofs.
\begin{lemma}\label{l:et}\mbox{}
    \begin{enumerate}
        \item  Let $f\colon U\to X$ be an \'etale covering Kan fibration. Then for all $n$, the map
	       \[
		          f_n\colon U_n\to X_n
	       \]
	       is a surjective local diffeomorphism. In particular, if $X\in\ssm$, then $U\in\ssm$ as well.
        \item Let $f\colon X\to Y$ be an \'etale hypercover. Then $\tau_{\le n}(X,f)\cong \tau_{<n}(X,f)\to X$ is an \'etale $n$-hypercover and $p\colon X\to \tau_{\le n}(X,f)\cong \tau_{<n}(X,f)$ is an \'etale hypercover.
    \end{enumerate}
\end{lemma}

We have a mild strengthening of Lemma~\ref{l:kanngpd} in this setting.
\begin{lemma}\label{l:etkan}
    Suppose we have a diagram in $\ssm$ 
    \[
        \xymatrix{
         X \ar[rr]^f \ar[dr]_g && Y \ar[dl]^h\\
         & Z
        }
    \]
    in which $g$ is an \'etale $n$-stack. 
    \begin{enumerate}
        \item If $f$ is a Kan fibration, then $f$ is an \'etale $n$-stack. 
        \item If $f$ is a covering Kan fibration and $d_i\colon Z_k\to Z_0$ is a cover for all $k$ and $i$, then $h$ is also an \'etale $n$-stack. 
    \end{enumerate}
\end{lemma}
\begin{proof}
    For the first statement, by Lemma~\ref{l:kanngpd} (applied in the category with covers $\ssm$), we only need to show that $f$ is an \'etale Kan fibration. For this, consider the commuting triangle
    \[
        \xymatrix{
            X_k \ar[r]^{\lambda^k_i(f)} \ar[dr]_{\lambda^k_i(g)} & \Lambda^k_i(f) \ar[d] \\
            & \Lambda^k_i(g)
        }
    \]
    By assumption, $\lambda^k_i(g)$ is a surjective local diffeomorphism and $\lambda^k_i(f)$ is a surjective submersion. We conclude that $\lambda^k_i(f)$ must be a local diffeomorphism (by the inverse function theorem). 

    For the second statement, the argument above shows that $f$ is an \'etale Kan fibration. This statement is now the special case of the statement in Lemma~\ref{l:kanngpd} interpreted in the descent category $\Sh(\sm)_{\et}$.
\end{proof}

\subsection{Locally minimal Kan fibrations.}
\begin{definition}\label{d:locmin}
    Let $f\colon X\to Y$ be a Kan fibration in $\ssm$.  We say $f$ is {\em locally minimal} if the map $X_0\to \pi_0(Y_0\times_Y X)$ is a local diffeomorphism of sheaves and, for all $n>0$, the natural map
    \[
        X_n\to \pi_0(P^{\ge n}(f))
    \]
    is a surjective local diffeomorphism of sheaves. A {\em locally minimal model} for $f$ is a commuting triangle in $\ssm$
    \[
        \xymatrix{
            \tilde{X} \ar[dr]_{\tilde{f}} \ar[rr]^\sim && X \ar[dl]^f \\
            & Y
        }
    \]
    where $\tilde{f}$ is a locally minimal Kan fibration and the horizontal map is a local weak equivalence.
\end{definition}

We have counterparts of Lemmas~\ref{l:minpost} and~\ref{l:duskmooremin} above.  Recall from Notation~\ref{n:etloc} that $\Sh(\sm)_{\loc}$ denotes the category with covers given by sheaves with local epimorphisms as covers, while $\Sh(\sm)_{\et}$ denotes the descent category given by sheaves with surjective local diffeomorphisms as covers.
 
\begin{lemma}\label{l:locmin}\mbox{} 
    \begin{enumerate}
        \item\label{l:locminhyp} \'Etale hypercovers of simplicial manifolds are locally minimal Kan fibrations.
        \item\label{l:locmincrit}  A Kan fibration $f\colon X\to Y$ is locally minimal if and only if the map
            \[
                p_n\colon \tau_{<n+1}(X,f)\to \tau_{\le n}(X,f)
            \]
            is an \'etale hypercover of simplicial sheaves for all $n\ge 0$.
        \item\label{l:locmindusk} Let $f\colon X\to Y$ be a locally minimal Kan fibration. Suppose that $\pi_0 P^{\ge n}(f)$ is a smooth manifold. Then $\tau_{\le n}(X,f)\to Y$ is a locally minimal Kan fibration, and $X\to\tau_{\le n}(X,f)$ is a Kan fibration in $\ssm$.
        \item\label{l:locminmoore} Let $f\colon X\to Y$ be a locally minimal Kan fibration. Suppose that $\Image(\mu_n(f))$ is a smooth manifold. Then $\tau_{<n}(X,f)\to Y$ is locally minimal Kan fibration in $\ssm$.
    \end{enumerate}
\end{lemma}
\begin{proof}
    The first statement follows immediately from the definitions and \cite[Theorem 3.7]{W}, noting that the proof carries over verbatim if we work in a descent category, rather than a category with covers.

    For the second statement, by Lemma~\ref{l:mooredusk} applied in $\sSh(\sm)_{\loc}$, the map $p_n$ is an local trivial $(n+1)$-fibration of simplicial sheaves. As the map $p_n$ is the identity on $(n-1)$-skeleta, it is an \'etale hypercover if and only if the map 
     \[
        \mu_n(p_n)\colon X_n\to M_n(p_n)\cong \pi_0 P^{\ge n}(f)
     \]
     is a surjective local diffeomorphism of sheaves. But this is a surjective local diffeomorphism for all $n$ if and only if $f$ is locally minimal. 

    For the third statement, consider the map of smooth manifolds
    \[
        X_n\to\pi_0P^{\ge n}(f).
    \]
    By our assumption that $f$ is locally minimal, this is a surjective local diffeomorphism, and in particular a cover.  The statement now follows from Lemma~\ref{l:kan} appied in $\ssm$, and from Lemma~\ref{l:trtr}.

    For the last statement, the same reasoning as in the proof of \cite[Proposition 7.9]{Hen}, we conclude that $\tau_{<n}(X,f)\to Y$ is a Kan fibration in $\ssm$. Its local minimality now follows from the second statement and Lemma~\ref{l:trtr}.
\end{proof}

Recall, following \cite[p. 1034]{Hen}, that a {\em finite dimensional diffeological group} is a group object in $\Sh(\sm)$ of the form $G/A$ where $G$ is a finite dimensional Lie group and $A$ is a finitely generated, but possibly non-discrete subgroup of $G$. 
\begin{lemma}\label{l:locminfd}
    Let $X$ be a Lie $\infty$-group such that for all $i$, $\pi_iX$ is a finite dimensional diffeological group.  Then:
    \begin{enumerate}
        \item a locally minimal model $X'\to^\sim X$ is necessarily a finite dimensional Lie $\infty$-group.
        \item the dimensions $\{\dim X'_n\}_{n\ge 0}$ are invariants of $X$: for any weakly equivalent Lie $\infty$-group $Y\simeq X$ and any locally minimal model $Y'\to^\sim Y$, $\dim X'_n=\dim Y'_n$ for all $n$.
    \end{enumerate}
\end{lemma}
\begin{proof}
    For the first statement, by \cite[Proposition 7.9]{Hen}, $\tau_{<n}X'$ is a Lie $\infty$-group for all $n$. Using that $X'=\varprojlim_n \tau_{<n}X'$ and that the map $X'\to\tau_{<n}X'$ is an isomorphism on $k$ skeleta for $k<n$, it suffices to show that $\tau_{<n}X'$ is finite dimensional for all $n$. We do this by induction.  For the base case, $\tau_{<1}X'=\Delta^0$ (this is true by definition for any smooth $\infty$-group in any category with covers). Now suppose that $\tau_{<n}X'$ is finite dimensional. The map of simplicial manifolds
    \[
        \tau_{<n+1}X'\to\tau_{<n}X'
    \]
    factors as 
    \[
        \tau_{<n+1}X'\to\tau_{\le n}X'\to\tau_{<n}X'.
    \]
    By Lemma~\ref{l:locmin}(\ref{l:locmincrit}), 
    \[
        \tau_{<n+1}X'\to\tau_{\le n}X'
    \]
    is an \'etale hypercover. By Proposition~\ref{p:postbund} applied in $\sSh(\sm)_{\loc}$, 
    \[
        \tau_{\le n}X'\to\tau_{<n}X
    \]
    is a principal $K(\pi_nX,n)//\pi_1X\to N\pi_1X$ bundle. By assumption, $\pi_n X$ is a finite dimensional diffeological group. Let $m=\dim\pi_n X$ (equivalently, $m$ is the dimension of the universal cover of $\pi_nX$ which is a finite dimensional Lie group). By Proposition~\ref{p:postbund}, applied in $\sSh(\sm)_{\loc}$, for all $k$, the maps
    \[
        (\tau_{\le n}X')_k\to(\tau_{<n}X')_k
    \]
    are fiber bundles with fibers $m\binom{k}{n}$-dimensional non-Hausdorff manifolds.  By Lemma~\ref{l:et}, the maps 
    \[
        (\tau_{<n+1}X')_k\to (\tau_{\le n}X')_k
    \]
    are surjective \'etale maps, in particular the fibers over any section are 0-dimensional smooth manifolds. We conclude that, for all $k$, the maps of smooth manifolds
    \[
        (\tau_{<n+1}X')_k\to (\tau_{<n}X')_k
    \]
    are surjective submersions, and the fibers at each point are $m\binom{k}{n}$-dimensional smooth manifolds. Because $(\tau_{<n}X')_k$ is finite dimensional for all $k$, by the inductive hypothesis, we conclude that $(\tau_{<n+1}X')_k$ is finite dimensional for all $k$ as well. This concludes the inductive step and thus the proof of the first statement.

    The second statement follows from the argument above via a simple induction.  Indeed, let $X'$ and $Y'$ be locally minimal models as above. We induct on $n$ to show that $\dim(\tau_{<n} X')_k=\dim(\tau_{<n}Y')_k$ for all $k$ and $n$. The base is trivial, as $X'_0=Y'_0=\ast$, and this implies that $\tau_{<1} X'=\tau_{<1}Y'=\Delta^0$. Suppose we have shown that, for $\ell\le n$, $\dim (\tau_{<\ell}X')_k=\dim (\tau_{<\ell} Y')_k$ for all $k$. Then, as noted above, the maps of smooth manifolds
    \begin{align*}
        (\tau_{<n+1} X')_k&\to (\tau_{<n}X')_k\\
        (\tau_{<n+1} Y')_k&\to (\tau_{<n}Y')_k
    \end{align*}
    are surjective submersions and the fibers at each point are $m\binom{k}{n}$-dimensional smooth manifolds.  Note here that we are using that $X\simeq Y$ implies that $\pi_n X\cong \pi_n Y$ and thus $m=\dim \pi_n X=\dim \pi_n Y$.  By the inductive hypothesis, we conclude that $\dim (\tau_{<n+1}X')_k=\dim (\tau_{<n+1}Y')_k$ for all $k$. This concludes the induction step and thus the proof.
\end{proof}

\subsection{Good \'etale hypercovers.}

Our goal in this subsection is to state and prove Proposition~\ref{p:hypexist}. We begin by recalling some facts about Reedy diagrams and simplicial spaces. We follow the notation of Hirschhorn \cite{Hirschhorn}.  Let $\Delta$ denote the ordinal category, let $\overrightarrow{\Delta^{\op}}$ denote the subcategory of surjective maps (in $\Delta$), and for $[n]\in\Delta$, recall that the {\em latching category} $\partial(\overrightarrow{\Delta^{\op}}\downarrow [n])$ is the full subcategory of $(\overrightarrow{\Delta^{\op}}\downarrow [n])$ on all objects except the identity map. This provides one of the canonical examples of a Reedy category, and in particular, it carries a degree function 
\[
\delta([n+1]\to[k]):=k
\]  
Given a simplicial manifold
\[
X\colon\Delta^{\op}\to \sm,
\]
by a standard exercise (cf. \cite[Proposition 15.8.6]{Hirschhorn}), for all $n\ge 0$, we can identify the degeneracies 
\begin{align}\label{e:latch}
	s(X_n)&:=\bigcup_{0\le i\le n} s_i(X_n)\subset X_{n+1}\nonumber\\
	&\cong L_{n+1} X 
\end{align}
where the latching object $L_{n+1}X$ is the colimit
\[
L_{n+1}X:=\colim_{\partial(\overrightarrow{\Delta^{\op}}\downarrow [n+1])} X.
\]

\begin{lemma}\label{l:degen}
	Let $X$ be a finite dimensional simplicial manifold. Let $n\ge 0$. The map
		\[
		\colim_{\partial(\overrightarrow{\Delta^{\op}}\downarrow[n+1])^{\op}} X\to s(X_n)
		\]
		is a homotopy colimit (of topological spaces). In particular, if $X\to Y$ is a map of simplicial manifolds which restricts to a homotopy equivalence $X_k\simeq Y_k$ for $k\le n$, then $s(X_n)$ is homotopy equivalent to $s(Y_n)$.
\end{lemma}

\begin{proof}
	The claim follows immediately from the isomorphism \eqref{e:latch} and the theory of Reedy categories.  First, observe that because every finite dimensional smooth manifold is triangularizable \cite{Cairns,Whitney}, a finite dimensional simplicial manifold $X$ is Reedy cofibrant as a simplicial space  \cite[Corollary 15.8.8]{Hirschhorn}.  Therefore 
	\[
	X\colon \partial(\overrightarrow{\Delta^{\op}}\downarrow[n+1])\to \sm
	\]
	is a Reedy cofibrant diagram too \cite[Lemma 15.3.7]{Hirschhorn}. The statement now follows from \cite[Theorem 19.9.1(1)]{Hirschhorn} because $\partial(\overrightarrow{\Delta^{\op}}\downarrow[n+1])$ is a Reedy category with fibrant constants (e.g. it satisfies the criterion of \cite[Proposition 15.10.2]{Hirschhorn}).
\end{proof}

\begin{proposition}\label{p:hypexist}
	Let $X$ be a finite dimensional simplicial manifold.  Assume that there exists:
	\begin{enumerate}
		\item a pointed manifold $y\in Y$, 
		\item a monotone function $f\colon \Nb\to\Nb$, and 
		\item local diffeomorphisms $\varphi_n\colon X_n\to^\cong Y^{f(n)}$ for each $n$ which conjugate the degeneracy maps into inclusions of factors, i.e. for all $0\le i\le n-1$, there exists an ordered inclusion $I\colon \{1,\ldots,f(n-1)\}\into \{1,\ldots,f(n)\}$ such that 
		\[
			\varphi_n\circ s_i\circ\varphi_{n-1}^{-1}=\iota_I
		\]
		where $\iota_I\colon Y^{f(n-1)}\to Y^{f(n)}$ is the map which inserts copies of the basepoint $y$ in all coordinates indexed by the complement $I(\{1,\ldots,f(n-1)\})$. 		
	\end{enumerate}
	Then any cover of $X_0$ by contractible opens extends to a good \'etale hypercover $U\to X$. 
	
	Similarly, if $X_n=\ast$ for some $n$ (and hence $X_m=\ast$ for all $m\le n$), then any cover of $X_{n+1}$ by contractible opens extends to a good \'etale hypercover $U\to X$.
\end{proposition}

\begin{remark}
		By Lemma~\ref{l:et}, any \'etale hypercover of a finite dimensional simplicial manifold is again finite dimensional.
\end{remark}

\begin{remark}
Statement \eqref{t:main1.1} of Theorem \ref{t:main1} is an immediate corollary of this proposition. Indeed, for a Lie group $G$, the nerve $NG$ satisfies the assumptions of the proposition, with $y\in Y$ given by $e\in G$ and $f(n)=n$.  Similarly, let $A$ be an abelian Lie group. Then for any $m$, the simplicial manifold $K(A,m)$ satisfies the assumptions of the proposition with $y\in Y$ given by $0\in A$ and $f(n)=\binom{n}{m}$. 
\end{remark}

\begin{remark}
	The construction is a modification of the classical construction of ``split hypercovers'' (cf. Artin--Mazur \cite[Chapter 8]{AM}).  The main difference is that split simplicial manifolds cannot in general satisfy the Kan condition, so we need to work more carefully to handle the degeneracies in the induction step of the construction. Lemma~\ref{l:degen} and the theory of tubular neighborhoods for stratified spaces  (cf. Quinn \cite{Q}) provide the key tools we need to accomplish this.  Indeed, the proposition's assumptions on $X$ are intended to guarantee that $X_n$, stratified by the order of degeneracy, is a manifold homotopy stratified space. It would be interesting to know if this is true in general, or if $X$ is assumed Kan.
\end{remark}

\begin{proof}[Proof of Proposition~\ref{p:hypexist}]
	Let $U_0\to X_0$ be a cover of $X_0$ by contractible opens. We prove, by induction on $n$, that there exists an $n$-truncated \'etale hypercover
	\[
	f_{\le n}\colon U_{\le n}\to X_{\le n}	
	\]
	such that for all $i\le n$, the map $U_i\to \pi_0(U_i)$ is a homotopy equivalence, and such that $(f_{\le n})|_{(U_{\le n})_{\le n-1}}=f_{\le n-1}$.
	
	For the base of the induction, $n=0$, we can just take $U_0\to X_0$ to be the given cover.
	
	For the inductive step, suppose we have an $n$-truncated good \'etale hypercover
	\[
	f_{\le n}\colon U_{\le n}\to X_{\le n}
	\]
	extending $U_0\to X_0$. Taking the relative $n$-coskeleton, we obtain a hypercover
	\[
	U^n:=\csk_n U_{\le n}\times_{\Csk_n X}X\to X
	\]
	By \cite[Lemma 3.3]{BG} applied in $\sSh(\sm)_{\et}$, $U^n\to X$ is an \'etale hypercover.  Because $X$ is a finite dimensional simplicial manifold, we conclude by Lemma~\ref{l:et} that $U^n$ is as well.
	
	By our assumption on $X$, the collection of closed submanifolds $\{\theta^*(X_k)\}_{\theta\colon [n+1]\onto [k]}$ of $X_{n+1}$ is {\em locally flat}, i.e. each point $x\in X_{n+1}$ admits a neighborhood $U$ and a diffeomorphism $U\to \Rb^m$ such that each $U\cap \theta^*(X_k)$ is taken into a linear subspace under this diffeomorphism.  In particular, by \cite[Proposition 1.5]{Q}, the stratification of $X_{n+1}$ by the order of degeneracy makes $X_{n+1}$ into a manifold homotopy stratified space.  Similarly, by Lemma~\ref{l:et}, the map $f\colon U^n\to X$ gives a local diffeomorphism $f_k\colon U^n_k\to X_k$ on $k$-simplices for all $k$.  Therefore, the collection of closed submanifolds $\{\theta^*(U^n_k)\}_{\theta\colon [n+1]\onto [k]}$ of $U^n_{n+1}$ is also locally flat, and the stratification of $U^n_{n+1}=M_{n+1}U_{\le n}\times_{M_{n+1}X}X_{n+1}$ by the order of degeneracy makes $U^n_{n+1}$ into a manifold homotopy stratified space as well. 
	
	By \cite[Proposition 3.2]{Q}, there exists an open neighborhood $V\subset U^n_{n+1}$ of $s(U_n)$ such that $V$ deformation retracts onto $s(U_n)$. Further, because smooth maps between finite dimensional manifolds are dense in the space of continuous maps (e.g. \cite[Ch. 2.2, Theorem 2.6]{Hirsch}) and the space of maps is locally path connected, we can assume without loss of generality that the deformation retraction is smooth. Now let $U_{n+1}$ be the disjoint union of $V$ and some locally finite cover by open balls of $U^n_{n+1}\setminus s(U_n)$. By Lemma~\ref{l:degen} and our assumption on $U_{\le n}$, $V\simeq s(U_n)\simeq \pi_0(s(U_n))$.  We conclude that $U_{n+1}$ is a disjoint union of contractible pieces.  Further, the inclusion $s(U_n)\subset V$ and the map 
	\[
	U_{n+1}\to U^n_{n+1}=M_{n+1}U_{\le n}\times_{M_{n+1}X}X_{n+1}
	\]
	define an $n+1$-truncated \'etale hypercover
	\[
	f_{\le n+1}:U_{\le n+1}\to X_{\le n+1}
	\]
	(where the $n$-truncation of $f_{\le n+1}$ equals $f_{\le n}$ and where $f_{n+1}\colon (U_{\le n+1})_{n+1}\to X_{n+1}$ is the map $U_{n+1}\to M_{n+1}U_{\le n}\times_{M_{n+1}X}X_{n+1}\to X_{n+1}$). By construction $U_i\simeq \pi_0(U_i)$ for all $i\le n+1$. This completes the inductive step and thus the proof.

	If $X_n=\ast$, then $X_{\le n}=\Delta^0_{\le n}$ (since $X_m$ is a retract of $X_n$ for all $m\le n$).  We can then take $U_{\le n}=X_{\le n}=\Delta^0_{\le n}$, and observe that in the $(n+1)^{st}$ stage of the induction, we have
	\[
		\csk_n U_{\le n}\times_{\Csk_n X}X=X.
	\]
	We can thus take $U_{n+1}$ to be the given cover of $X_{n+1}$ by contractible opens, and proceed with the induction as above.
\end{proof}

\section{A locally minimal model for the Sullivan integral}\label{s:base}
Let $\g$ be a finite dimensional Lie algebra and let $G$ be the simply connected Lie group integrating $\g$. The goal of this section is to prove that there exists a finite dimensional model of the Sullivan integration $\sint\g$. 

We begin with some differential topology. Let $\Map_\ast(\Delta^n,G)$ denote the Banach manifold of $C^r$-maps $f\colon \Delta^n\to G$ such that $f(0)=e$ as in Example~\ref{ex:intnerve} above.

\begin{lemma}\label{l:key}
	For all $n$, the map
	\[
	\pi\colon \Map_\ast(\Delta^n,G)\to \Map_\ast(\Delta^n,G)/\sim_{\text{rel }\partial}
	\]
	is a locally trivial principal bundle for the group $\Map_\ast^0((\Delta^n,\partial\Delta^n),(G,e))$, i.e. the identity component of the based $n$-fold loop group of $G$.
\end{lemma}
\begin{proof}
	It suffices to prove that there exists an isomorphism over $\Map_\ast(\Delta^n,G)$
	\begin{align*}
	Map_\ast(\Delta^n,G)\times_{\Map_{\ast}(\Delta^n,G)/\sim_{\text{rel }\partial}}\Map_\ast(\Delta^n,G)&\to^\cong \Map_\ast(\Delta^n,G)\times \Map_\ast^0((\Delta^n,\partial\Delta^n),(G,e))\\
 (x_1,x_2)&\mapsto (x_1,\varphi(x_1,x_2))
	\end{align*}
	such that $\varphi$ satisifes the cocycle condition 
	\[
	\varphi(x_1,x_2)\cdot\varphi(x_2,x_3))=\varphi(x_1,x_3).
	\]
    Indeed, granting this, we see that the source of this isomorphism, viewed as a subspace of $\Map_\ast(\Delta^n,G)\times\Map_\ast(\Delta^n,G)$, is the equivalence relation $\sim_{\text{rel }\partial}$. This subspace is closed, by inspection, and by the isomorphism above, is an embedded submanifold whose projection onto the first factor is a surjective submersion.  By Godemont's Theorem \cite[5.9.5]{VAR}, we conclude that the map
    \[
        \Map_\ast(\Delta^n,G)\to \Map_\ast(\Delta^n,G)/\sim_{\text{rel }\partial}
    \]
    is a surjective submersion.  Using the implicit function theorem and the isomorphism above, we obtain a local trivialization as desired.
 
	To construct $\varphi$, we use the Lie group structure on $G$.  Indeed, the space
	\[
	Map_\ast(\Delta^n,G)\times_{\Map_{\ast}(\Delta^n,G)/\sim_{\text{rel }\partial}}\Map_\ast(\Delta^n,G)
	\]
	consists, by definition of pairs of maps $x_1,x_2\in \Map_\ast(\Delta^n,G)$ which are homotopic rel boundary.  In particular, their boundaries are identical, so taking $\varphi$ to be the pointwise product
	\[
	\varphi(x_1,x_2):=x_1^{-1}\cdot x_2
	\]
	we obtain the desired isomorphism satisfying the cocycle condition.
\end{proof}


We can now state and prove the main result of this section.
\begin{theorem}\label{t:splitexist}
	Let $\g$ be a finite dimensional Lie algebra and let $G$ be the associated simply connected Lie group.  Let $n\in \N_{\ge 1}\cup\{\infty\}$, and let $h\colon U\to NG$ be a reduced, \'etale hypercover which is good up to level $n-1$ (Definition~\ref{d:good}). There exists a commuting square
	\begin{equation}\label{e:finsull}
	\xymatrix{
		 \widetilde{U}^n \ar[r]^-\sim \ar[d]_{q_n}  & \tau_{\le n}\sint\g \ar[d]^{\tau_{\le 1}} \\
		 U \ar[r]^-\sim & NG
	}
	\end{equation}
    such that the left vertical map is a minimal Kan fibration, and $\widetilde{U}^n$ is finite dimensional, reduced and good up to level $n-1$.  Further, the map 
    \[
        \widetilde{U}^n\to^\sim\tau_{\le n}\sint\g
    \]
    is a locally minimal model for $\tau_{\le n}\sint\g$.
\end{theorem}

Taking $n=\infty$ as in Remark~\ref{r:infpost} and writing $q$ and $\cG^0$ in lieu of $q_\infty$ and $\widetilde{U}^\infty$, we obtain statement \eqref{t:main1.2} of Theorem \ref{t:main1} and Theorem ~\ref{t:main4} in the case when $L$ is a Lie algebra:
\begin{corollary}\label{c:splitexist}
    Let $\g$ be a finite dimensional Lie algebra and let $G$ be the associated simply connected Lie group.  For any reduced, good \'etale hypercover $U\to NG$, there exists a commuting square
	\begin{equation}\label{e:finsullinf}
	\xymatrix{
		 \cG^0 \ar[r]^-\sim \ar[d]_{q}  & \sint\g \ar[d]^{\tau_{\le 1}} \\
	   U \ar[r]^-\sim & NG
	}
	\end{equation}
    such that the left vertical map is a minimal Kan fibration, and $\cG^0$ is finite dimensional, reduced and good.  Further, the map 
    \[
        \cG^0\to^\sim\sint\g
    \]
    is a locally minimal model for $\sint\g$.
\end{corollary}
    
\begin{proof}[Proof of Theorem~\ref{t:splitexist}]
Note that by \cite[Theorem 7.5]{Hen}, $\tau_{\le n}\sint\g$ is a Lie $n$-group for all $n$. Fix $n\ge 1$ and let $U\to NG$ be a reduced, \'etale hypercover which is good up to level $(n-1)$ as in the statement of the theorem. We prove by induction on $1\le m\le n$ that there exists a commuting square
\begin{equation}\label{e:finsullind}
\xymatrix{
	\widetilde{U}^m  \ar[d]_{q_m} \ar[r]^-\sim & \tau_{\le m}\sint\g \ar[d] \\
	U \ar[r]^-\sim & NG
}
\end{equation}
in which the left vertical map is a minimal Kan fibration, the horizontal maps are local weak equivalences, $\tau_{\le m-1}(\widetilde{U}^m,q_m)\cong \widetilde{U}^{m-1}$ over $U$, and both $\widetilde{U}^m$ and $\widetilde{U}^{m-1}$ are good up to level $(n-1)$.  In the case where $n=\infty$, we complete the transfinite induction by taking the inverse limit in $n$ of the squares~\eqref{e:finsull}. Because $\tau_{\le m-1}(\widetilde{U},q_m)^m\cong \widetilde{U}^{m-1}$ over $U$, we obtain a well-defined finite dimensional reduced, good simplicial manifold $\widetilde{U}^\infty$ fitting into the square \eqref{e:finsull} as desired. By construction, for all $m$, the map
\[
    \tau_{\le m}(\tilde{U}^\infty,q)\cong \widetilde{U}^m\to \widetilde{U}^{m-1}\cong\tau_{\le m-1}(\tilde{U}^\infty,q)
\]
is a minimal Kan fibration. By Lemma~\ref{l:postlim}, we conclude that $q_\infty$ is a minimal Kan fibration as well. 

For the base of the induction, $m=1$, $\tau_{\le 1}\sint\g=NG$, by \cite[Examples 5.5 and 7.2]{Hen}, so we can just take $\widetilde{U}^1:=U$ as above.

Now suppose we have a commuting square as above for $1\le m<n$. Consider the maps
\[
\tau_{\le m+1}\sint\g\to \tau_{<m+1}\sint \g\to \tau_{\le m}\sint \g.	
\]	
We claim that there exists a lift
\[
\xymatrix{
	& \tau_{<m+1}\sint \g \ar[d]\\
	\widetilde{U}^m \ar[r] \ar@{..>}[ur]^{\exists} & \tau_{\le m}\sint\g
}
\]
To prove this claim, recall that given a diagram of simplicial objects in a category $\C$
\[
	\xymatrix{
	& Z \ar[d]^f\\
	X \ar[r] & Y
	}
\]
and a lift of $(k-1)$-truncated objects
\[
\xymatrix{
	& Z_{\le k-1} \ar[d]\\
	X_{\le k-1} \ar[r] \ar[ur] & Y_{\le k-1}
},
\]
to give an extension of this to a lift of $k$-truncated objects is equivalent to giving a solution of the lifting problem 
\[
	\xymatrix{
	s(X_{k-1}) \ar[r] \ar[d] & Z_k \ar[d]\\
	X_k \ar[r] \ar@{..>}[ur]^{\exists} & M_k(f)
	},
\]
where, to avoid representability issues, we view the above diagram as taking place in the category of presheaves on $\C$ and invoke the Yoneda lemma.

We will apply this repeatedly here.  First, note that, by Lemma~\ref{l:mooredusk}, the map 
\[
	\tau_{<m+1}\sint\g\to \tau_{\le m}\sint \g
\]
is an $(m+1)$-hypercover. Moreover, it restricts to an equality on $m-1$-skeleta.   Taking $k=m$, we obtain a lift of $(m-1)$-truncated objects
\[
\xymatrix{
	& (\tau_{<m+1}\sint \g)_{\le m-1} \ar@{=}[d]\\
	(\widetilde{U}^m) _{\le m-1}\ar[r] \ar[ur] & (\tau_{\le m}\sint\g)_{\le m-1}
}.
\]
To extend this to a lift of $m$-truncated objects, we need to solve the lifting problem 
\[
\xymatrix{
	s(\widetilde{U}^m_{m-1}) \ar[r] \ar[d] & (\tau_{<m+1}\sint\g)_m=(\sint\g)_m \ar[d]\\
	\widetilde{U}^m_m \ar[r] \ar@{..>}[ur] & (\tau_{\le m}\sint\g)_m
}.
\]
Note that we can identify $M_m(\tau_{<m+1}\sint\g\to \tau_{\le m}\sint \g)$ with $(\tau_{\le m}\sint\g)_m$ precisely because the map is an equality on $m-1$-skeleta.

As observed in the discussion preceding Lemma~\ref{l:key}, the right vertical map in the square above is isomorphic to the map
\[
\pi\colon \Map_\ast(\Delta^m,G)\to \Map_\ast(\Delta^m,G)/\sim_{\text{rel }\partial}.
\]
By Lemma~\ref{l:key}, this is a locally trivial principal bundle, in particular it is a smooth Serre fibration. By our inductive assumption that $m<n$ and $\widetilde{U}^m$ is good up to level $n-1$, $\widetilde{U}^m_m$ is a disjoint union of contractible pieces. In other words, if we denote by $\pi \widetilde{U}^m$ the simplicial set with $k$-simplices $\pi_0(\widetilde{U}^m_k)$, then the canonical map
\[
(\widetilde{U}^m)_{\le n-1}\to (\pi \widetilde{U}^m)_{\le n-1}
\]
is a levelwise homotopy equivalence of finite dimensional $(n-1)$-truncated simplicial manifolds. Because $m<n$, Lemma~\ref{l:degen} implies that
\[
s(\widetilde{U}^m_{m-1})\simeq s((\pi \widetilde{U}^m)_{m-1}).
\]
In particular, the cofibration $s(\widetilde{U}^m_{m-1})\to \widetilde{U}^m_m$ induces an injection
\[
\pi_0(s(\widetilde{U}^m_{m-1}))\to \pi_0(\widetilde{U}^m_m).
\]
We conclude that $s(\widetilde{U}^m_{m-1})\to \widetilde{U}^m_m$ is a trivial cofibration onto the connected components in its image.  Therefore, by the lifting property of Serre fibrations, on each such component $\widetilde{U}^m_{m,a}$, there exists a lift
\[
\xymatrix{
	s(\widetilde{U}^m_{m-1})_a \ar[d] \ar[r] & (\tau_{<m+1}\sint\g)_m \ar[d]\\
	\widetilde{U}^m_{m,a} \ar[ur] \ar[r] & (\tau_{\le m}\sint \g)_m
}
\]
Moreover, by our inductive assumption that $\widetilde{U}^m$ is good up to level $n-1$, every connected component of $\widetilde{U}^m_m$ is contractible.  Therefore, for all the connected components $\widetilde{U}^m_{m,b}$ not containing $s(\widetilde{U}^m_{m-1})$, the pullback bundle
\[
(\tau_{<m+1}\sint\g)|_{\widetilde{U}^m_{m,b}}\to \widetilde{U}^m_{m,b}
\]
is trivializeable.  Choosing any section for each such component $\widetilde{U}^m_{m,b}$, we obtain the desired lift
\[
\widetilde{U}^m_m\to (\tau_{<m+1}\sint\g)_m.
\]
As observed above, this defines a commuting diagram of $m$-truncated simplicial objects
\[
\xymatrix{
	& (\tau_{<m+1}\sint\g)_{\le m}\ar[d] \\
	\widetilde{U}^m_{\le m} \ar[r] \ar[ur] & (\tau_{\le m}\sint\g)_{\le m}
}
\]
We now need to extend this to a map of $(m+1)$-truncated objects. As above, such an extension is equivalent to a solution of the lifting problem
\[
\xymatrix{
	s(\widetilde{U}^m_m) \ar[r] \ar[d] & (\tau_{<m+1}\sint\g)_{m+1} \ar[d]\\
	\widetilde{U}^m_{m+1} \ar[r] \ar@{..>}[ur] & M_{m+1}\left(\tau_{<m+1}\sint\g\to\tau_{\le m}\sint\g\right)
}.
\]
But as observed above, $\tau_{<m+1}\int\g\to \tau_{\le m}\int\g$ is an $(m+1)$-hypercover by Lemma~\ref{l:mooredusk}; in particular, the right vertical map in the square above is an isomorphism.
Therefore, the lift 
\[
	\widetilde{U}^m_{\le m}\to (\tau_{<m+1}\sint\g)_{\le m}
\]
extends uniquely to a map of $m+1$-truncated objects
\[
	\widetilde{U}^m_{\le m+1}\to (\tau_{<m+1}\sint\g)_{\le m+1}
\]
and, in fact, to the desired lift of simplicial objects
\[
	\xymatrix{
		& \tau_{<m+1}\sint \g \ar[d]\\
		\widetilde{U}^m \ar[r] \ar[ur]& \tau_{\le m}\sint\g
	},
\]
by the same argument (with $m+k$ in lieu of $m+1$ in the lifting problem above).

We now have a commuting diagram of simplicial manifolds
\[
	\xymatrix{
	& \tau_{\le m+1}\sint\g \ar[d] \\
	& \tau_{<m+1}\sint\g \ar[d] \\
	\widetilde{U}^m \ar[r] \ar[ur] & \tau_{\le m}\sint\g
	}.
\]
The top vertical map is a minimal Kan fibration by Lemma~\ref{l:duskmoore}. By Proposition~\ref{p:postbund}, it is also a locally trivial principal $K(\pi_{m+1} G,m+1)$-bundle; note that, here we are using that $\pi_1\sint\g$ is the simply connected group $G$, and $\pi_m\sint\g=\pi_mG$ for $m\ge 2$ \cite[Example 6.1]{Hen} to conclude that the action of $\pi_1\sint\g$ on $\pi_m\sint\g$ is trivial. 

Define $\widetilde{U}^{m+1}$ to be the fiber product
\[
	\widetilde{U}^{m+1}:=\widetilde{U}^m\times_{\tau_{<m+1}\sint\g}\tau_{\le m+1}\sint\g.
\]
Then $\widetilde{U}^{m+1}\to \widetilde{U}^m$ is a minimal Kan fibration, and a principal $K(\pi_{m+1} G,m+1)$-bundle over the finite dimensional simplicial manifold $\widetilde{U}^m$. Let $q_{m+1}$ denote the composite
\[
    \widetilde{U}^{m+1}\to\widetilde{U}^m\to^{q_m}U
\]
Because $K(\pi_{m+1} G,m+1)$ is a discrete simplicial group with $K(\pi_{m+1}G,m+1)_{\le m}=(\Delta^0)_{\le m}$, we immediately see that  $\widetilde{U}^{m+1}$ is a finite dimensional reduced Kan simplicial manifold, and that
\[
    \tau_{\le m}(\widetilde{U}^{m+1},q_{m+1})\cong\widetilde{U}^m
\]
over $U$. Further, for all $k$, we have a locally trivial principal $K(\pi_m G,m)_k$-bundle
\[
		\widetilde{U}^{m+1}_k\to\widetilde{U}^m_k
\]
By our inductive assumption, $\widetilde{U}^m$ is good up to level $n-1$. Therefore, $\widetilde{U}^{m+1}$ is also good up to level $n-1$, as for $k\le n-1$, $\widetilde{U}^{m+1}_k$ is a disjoint union of trivializeable fiber bundles with discrete fiber over contractible bases.  By construction, it sits in the pullback square
\begin{equation}\label{e:lie3ind}
 \xymatrix{
 	\widetilde{U}^{m+1}  \ar[d] \ar[r]^-\sim & \tau_{\le m+1}\sint\g \ar[d] \\
 	\widetilde{U}^{m} \ar[r]^-\sim & \tau_{<m+1}\sint\g
 }
\end{equation}
The left and right vertical maps are minimal Kan fibrations, the bottom horizontal map is a local weak equivalence (by the 2 of 3 property and the construction above).  Therefore, we conclude from Lemma~\ref{l:rightproper} that the top vertical map is also a local weak equivalence.  Postcomposing this square on the bottom right by the weak equivalence $\tau_{<m+1}\sint\g\to^\sim \tau_{\le m}\sint\g$, we obtain the commuting square~\eqref{e:finsullind}. This completes the induction step.

It remains to show that $\widetilde{U}^m$ is locally minimal. We induct on $m$, using that, by Lemmas~\ref{l:locmin}(\ref{l:locmincrit}), it suffices to show that for all $k$, the map
\[
        \tau_{<k+1}\widetilde{U}^m\to \tau_{\le k}\widetilde{U}^m=\widetilde{U}^k      
\]
is an \'etale hypercover. For the base case, $m=1$, note that $h\colon U\to NG$ is a reduced \'etale hypercover. In particular,
\[
    \tau_{<1}U=\tau_{\le 0}U=\Delta^0,
\]
and for $k\ge 1$, because $NG$ is a Lie 1-group, 
\begin{align*}
    \tau_{<k+1}(U,h)&\cong \tau_{<k+1}U\\
    \tau_{\le k}(U,h)&\cong \tau_{\le k}U.
\end{align*}
That $U$ is locally minimal now follows immediately from Lemma~\ref{l:locmin}(\ref{l:locminhyp} and \ref{l:locmincrit}).

Now suppose we have shown that $\widetilde{U}^m$ is locally minimal. For all $k\ge 0$, we have a commuting square
\[
    \xymatrix{
        \tau_{<k+1}\widetilde{U}^{m+1} \ar[r] \ar[d] & \tau_{\le k}\widetilde{U}^{m+1} \ar[d] \\
        \tau_{<k+1}\widetilde{U}^m \ar[r] & \tau_{\le k}\widetilde{U}^m
    }
\]
By the induction and Lemma~\ref{l:locmin}(\ref{l:locmincrit}), the bottom horizontal map is an \'etale hypercover. By \cite[Proposition 7.9]{Hen}, the left vertical map is in $\ssm$. By Lemma~\ref{l:tauexact} and our construction above, it is an \'etale Kan fibration.  If $k\le m$, the right vertical map is an isomorphism. If $k>m$, it is a principal $K(\pi_{m+1}G,m+1)$-bundle. In either case, it is an \'etale Kan fibration.  We conclude by Lemma~\ref{l:et} that the upper horizontal map is an \'etale Kan fibration. But by Lemma~\ref{l:mooredusk}, it is a hypercover, and is thus an \'etale hypercover as claimed. This completes the induction and thus the proof.
\end{proof}

\section{Proof of main theorem}\label{s:ind}
We are now in a position to prove Theorems \ref{t:main3} and \ref{t:main4}, thereby completing the proof of the main result Theorem \ref{t:lie3}. 

\subsection{Proof of Theorem \ref{t:main3}} \label{sec:main3-pf}
 Let $L$ be a finite type Lie $\infty$-algebra, and set $\g:=H_0(L)$. By \cite[Theorem 5.10]{Hen}, the map
    \[
        \sint L\to \sint \g
    \]
    is a Kan fibration which factors as a tower of Kan fibrations
    \[
        \sint L\to \cdots\to\sint \tau_{\le n+1}L\to \sint \tau_{\le n}L\to\cdots\to\sint \tau_{\le 0}L=\sint \g.
    \]
    We claim that for each $n> 0$, the Kan fibration 
    \[
        \sint \tau_{\le n}L\to \sint \tau_{\le n-1}L
    \]
    admits a minimal model
    \begin{equation}\label{e:minstage}
        \xymatrix{
        \mathcal{L}_{n,n-1} \ar[r]^\sim \ar[dr]_{q_{n,n-1}} & \sint \tau_{\le n}L\ar[d] \\
        & \sint \tau_{\le n-1}L
        }.
    \end{equation}
    Granting this claim, by Proposition~\ref{p:postbund} and \cite[Theorem 6.4]{Hen}, the map $q_{n,n-1}$ is a principal $K(H_nL,n+1)//G\to NG$ bundle. Set $\mathcal{L}_0=\int \g$. For all finite $n\ge 1$, we define 
    \[
        \mathcal{L}_n=\mathcal{L}_{n,n-1}\times_{\sint \tau_{\le n-1}L}\cdots\times_{\sint \tau_{\le 1}L}\mathcal{L}_{1,0}.
    \]
    Observe that, for all finite $n\ge 1$, the map 
    \[
        \mathcal{L}_n\to\mathcal{L}_{n-1}
    \]
    is a principal $K(H_nL,n+1)//G\to NG$-bundle.  Moreover, the above map is isomorphism on $n$-skeleta. We can thus define
    \[
        \mathcal{L}:=\varprojlim_n\mathcal{L}_n
    \]
    to obtain a well-defined simplicial manifold. By construction, for all $\infty\ge n\ge 0$, the simplicial manifold $\mathcal{L}_n$ sits in a commuting triangle
    \[
         \xymatrix{
            \mathcal{L}_n \ar[r]^\sim \ar[dr] & \sint \tau_{\le n} L\ar[d] \\
        & \sint \g
        }.
    \]
    By Lemma~\ref{l:minpost}, for all $\infty\ge n\ge 0$ we conclude that this triangle is a minimal model for the Kan fibration $\sint\tau_{\le n}L\to \sint\g$.
    
    It remains to prove the claim.  We begin by noting that the map
    \[
        \tau_{\le n}L\to \tau_{\le n-1}L
    \]
    factors as
    \[
        \tau_{\le n}L\to^{q_{\le n}} \tau_{<n}L\to^{q_{<n}}\tau_{\le n-1}L.
    \]
    As noted in Example~\ref{ex:q-split}, the map $q_{\le n}$ is a quasi-split fibration, with fiber the minimal abelian $L_\infty$-algebra $H_nL[n]$. By \cite[Theorem 9.16]{RZ}, 
    \[
        \sint q_{\le n}\colon \sint \tau_{\le n} L \to \sint \tau_{<n}L
    \]
    is a Kan fibration, with fiber $\sint H_nL[n]$. By \cite[Theorem 6.4]{Hen}, $\pi_i\sint H_nL[n]=0$ for $i\neq n+1$ and $\pi_{n+1}\sint H_nL[n]=H_nL$. Further, it follows from the explicit construction of $\int$ that $\tau_{<n+1}\int q_{\le n}$ is an isomorphism. We are therefore in the situation of Proposition~\ref{p:postbund}.  By Proposition~\ref{p:postbund}, we can apply the relative $(n+1)$-truncation functor to obtain a pair of Kan fibrations
    \begin{equation}\label{e:minsplit2}
        \sint\tau_{\le n}L \to^\sim \tau_{\le n+1}(\sint\tau_{\le n}L ,\sint q_{\le n})\to\sint\tau_{<n} L.
    \end{equation}
    in which the second map is a principal $K(H_nL,n+1)//G\to NG$ bundle.  

    As in Eq.\ \ref{eq:qlt-split} from Example \ref{ex:Ptower}, the acyclic fibration $q_{<n}$ splits. Picking a section and integrating, we obtain a map
    \begin{equation}\label{e:minsplit1}
        \sint\tau_{\le n-1}L \to^\sim \sint\tau_{<n}L
    \end{equation}
    which is a weak equivalence by \cite[Theorem 9.15]{RZ}.
    
    Now define $\mathcal{L}_{n,n-1}$ as the pullback
    \[
        \mathcal{L}_{n,n-1}:=\sint\tau_{\le n-1}L\times_{\sint\tau_{<n}L}\tau_{\le n+1}(\sint\tau_{\le n}L ,\sint q_{\le n})
    \]
    By Lemma~\ref{l:minkan}, the projection onto the first factor is a minimal fibration and a principal $K(H_nL,n+1)//G\to NG$ bundle. 

    By Lemma~\ref{l:rightproper}, the projection onto the second factor is a weak equivalence. Note that, by the discussion above, the first map in~\eqref{e:minsplit2} is a weak equivalence, and thus a hypercover by \cite[Proposition 6.7]{RZ}. By Lemma~\ref{l:postsplit}, the hypercover
    \[
        \sint\tau_{\le n} L \to^\sim \tau_{\le n+1}(\sint\tau_{\le n}L ,\sint q_{\le n})
    \]
    admits a section $\sigma$. Post-composing $\sigma$ with the projection of $\mathcal{L}_{n,n-1}$ onto the second factor, we obtain a weak equivalence 
    \[
        \mathcal{L}_{n,n-1}\to^\sim \sint\tau_{\le n} L
    \]
    and thus the claimed minimal model as in~\eqref{e:minstage}. \hfill \qed
\subsection{Proof of Theorem \ref{t:main4}} \label{sec:main4-pf}    
Let $\cG^0\to^\sim \sint\g$ be a good, reduced, locally minimal model for $\sint\g$ as in statement \eqref{t:main1.2} of Theorem \ref{t:main1}. Let $\mathcal{L}\to^\sim\sint L$ be a minimal model for the Kan fibration $\sint L\to \sint \g$, and let 
     \[
        \cG:=\cG^0\times_{\sint\g}\mathcal{L}
     \]
    We begin by showing that $\cG$ is good, reduced, and finite dimensional. We will then show that $\cG$ is locally minimal. 

    For the former, first observe by \cite[Theorem 6.4]{Hen}, Lemma~\ref{l:minkan} and Corollary~\ref{c:minbun}, the map
    \[
        \cG\to\cG^0
    \]
    factors as a composition 
    \[
        \cG\to\cdots\cG^n\to\cG^{n-1}\cdots\to\cG^0.
    \]
    where each $\cG^n\to\cG^{n-1}$ is a principal $K(H_nL,n+1)//G\to NG$-bundle. We show by induction on $n$ that each $\cG^n$ is good, reduced, and finite dimensional. Because $\cG_{\le k}=\cG^N_{\le k}$ for $N>>0$, this suffices to show that $\cG$ is good, reduced, and finite dimensional.
    
    For the base case, this is just Theorem~\ref{t:splitexist}. Now suppose we have shown this for $\cG^{n-1}$. The map
    \[
        \cG^n\to \cG^{n-1}
    \]
    is a principal $K(H_nL,n+1)//G\to NG$ bundle over a good, reduced, finite dimensional Lie $\infty$-group. Since $n>0$, we conclude that $\cG^n$ is reduced. Further, for each $k$, 
    \[
        \cG^n_k\to \cG^{n-1}_k
    \]
    is a principal $(H_nL)^{\binom{k}{n+1}}$-bundle (by Proposition~\ref{p:postbund}) and $\cG^{n-1}_k$ is a disjoint union of contractible pieces. We conclude that this bundle is trivializeable.  As it has contractible, finite dimensional fibers ($L$ is finite type!), $\cG^n_k$ is finite dimensional and a disjoint union of contractible pieces.  We conclude that $\cG^n$ is good, reduced, and finite dimensional, as claimed.  This completes the inductive step and thus the proof that $\cG$ is good, reduced, and finite dimensional. 

    It remains to prove that $\cG$ is locally minimal.  We again induct on $n$ to prove that $\cG^n$ is locally minimal for all $n$. The base case $n=0$ is Corollary~\ref{c:splitexist}. 
    
    For the inductive step, suppose $\cG^{n-1}$ is locally minimal. We will work in $\sSh(\sm)$, and following Notation~\ref{n:etloc}, we will write $\sSh(\sm)_{\et}$ when we work in the descent category given by sheaves and surjective local diffeomorphisms, and we will write $\sSh(\sm)_{\loc}$ when we work in the category with covers given by sheaves and local epimorphisms. We will use the characterization of locally minimal given in Lemma~\ref{l:locmin}(\ref{l:locmincrit}). 

    For all $k\ge 0$, we have a commuting square of simplicial sheaves
   \begin{equation}\label{e:locminproof}
    \xymatrix{
        \tau_{<k+1}\cG^{n} \ar[r] \ar[d] & \tau_{\le k}\cG^{n}\ar[d] \\
        \tau_{<k+1}\cG^{n-1} \ar[r] & \tau_{\le k}\cG^{n-1}
        }
   \end{equation}
   By Lemma~\ref{l:locmin}(\ref{l:locmincrit}) and our inductive hypothesis, the bottom horizontal map is an \'etale hypercover of simplicial sheaves. By Lemma~\ref{l:mooredusk} applied in $\sSh(\sm)_{\loc}$, the top horizontal map is a local trivial fibration of simplicial sheaves. As noted above, the map 
    \[
        \cG^n\to \cG^{n-1}
    \]
    is a minimal Kan fibration and a principal $K(H_nL,n+1)//G\to NG$ bundle. If $k\le n$, the vertical maps in the square~\eqref{e:locminproof} are isomorphisms, and thus the top horizontal map is an \'etale hypercover of simplicial sheaves. If $k>n+1$, the vertical maps remain $K(H_nL,n+1)//G\to NG$ bundles, and in particular minimal Kan fibrations.  In this case, we conclude by Lemma~\ref{l:minkan}, applied in $\sSh(\sm)_{\loc}$, that the square is a pullback. Because \'etale hypercovers of simplicial sheaves are preserved under pullback (by Lemma~\ref{l:desccat} and \cite[Lemma 3.16]{BG}), we conclude that the top horizontal map is an \'etale hypercover, as claimed.  

   It remains to consider the case $k=n+1$, i.e. to show that the map $p\colon \tau_{<n+2}\cG^n\to \tau_{\le n+1}\cG^n$ is an \'etale hypercover of simplicial sheaves. For this, consider the commuting diagram
    \begin{equation}\label{e:intlocmin1}
        \xymatrix{
        \tau_{<n+2}\cG^n \ar[r] \ar[dr] & P_1 \ar[r] \ar[d] & \tau_{\le n+1}\cG^n \ar[d] \\
            & \tau_{<n+2}\cG^{n-1} \ar[r]  & \tau_{\le n+1}\cG^{n-1}
        }
    \end{equation}
    in which the square is a pullback, and $p$ is the composition of the top horizontal maps. By the inductive hypothesis, the bottom horizontal map is an \'etale hypercover, so the top horizontal map in the square is too (by Lemma~\ref{l:desccat} and \cite[Lemma 3.17]{BG}). By Lemma~\ref{l:kanngpd} applied in $\sSh(\sm)_{\loc}$, the map
    \[
        \tau_{<n+2}\cG^n\to P_1
    \]
    is a local Kan fibration. Because \'etale Kan fibrations of simplicial sheaves compose (by Lemma~\ref{l:desccat} and \cite[Lemma 3.15]{BG}), it suffices to show that it is in fact an \'etale Kan fibration. For this, consider the commuting diagram    
    \begin{equation}\label{e:intlocmin2}
        \xymatrix{
        \tau_{\le n+2}\cG^n \ar[r] \ar[d] & P_2 \ar[r] \ar[d] & \tau_{\le n+2}\cG^{n-1} \ar[d] \\
        \tau_{<n+2}\cG^n \ar[r] & P_1 \ar[r] & \tau_{<n+2}\cG^{n-1}
        }
    \end{equation}
    in which the right square is a pullback.  By \cite[Theorem 6.4]{Hen} and \cite[Proposition 7.9]{Hen}, the lower left and lower right corners are Lie $\infty$-groups. By Lemma~\ref{l:duskmoore} and Proposition~\ref{p:postbund}, applied in $\ssm$, the right vertical map is a principal $K(\pi_{n+2}G,n+2)$ bundle in $\ssm$, while the left vertical map is a principal $K(\ker \partial_{n+1},n+2)$ bundle, where 
    \[
        \partial_{n+1}\colon\pi_{n+2}G\to H_nL
    \]
    is the boundary map of Henriques' long exact sequence in homotopy groups \cite[Theorem 6.4]{Hen}. Because principal bundles are are preserved under pullback, the middle vertical map is also a principal $K(\pi_{n+2}G,n+2)$-bundle, now in $\sSh(\sm)$. 
    
    By Lemma~\ref{l:duskmin} applied in $\sSh(\sm)_{\loc}$, the right vertical map in the square~\eqref{e:intlocmin1} is a minimal local Kan fibration of simplicial sheaves. By \cite[Theorem 6.4]{Hen} and Proposition~\ref{p:postbund} applied in $\sSh(\sm)_{\loc}$, it is a principal $K(H_nL/\partial_{n+1}(\pi_{n+2}G),n)//G\to NG$-bundle. Therefore, so are the top and bottom right horizontal maps in the diagram~\eqref{e:intlocmin2} above. By Lemma~\ref{l:duskmin} applied in $\ssm$, the composition of the top horizontal maps is a minimal Kan fibration, and thus by \cite[Theorem 6.4]{Hen} and Proposition~\ref{p:postbund}, it is a principal $K(H_nL,n+1)//G\to NG$ bundle. We conclude that the upper left horizontal map
    \[
        \tau_{\le n+2}\cG^n\to P_2
    \]
    is a principal $K(\partial_{n+1}(\pi_{n+2}G),n+1)$-bundle in $\sSh(\sm)$. Therefore, as observed in Remark~\ref{r:discgroupkinvar}, Theorem~\ref{t:kinvar} and \cite[Theorem 6.7]{W} together imply that principal bundles for discrete simplicial groups in $\sSh(\sm)_{\et}$ are \'etale Kan fibrations. Because \'etale Kan fibrations compose (Lemma~\ref{l:desccat} and \cite[Lemma 3.17]{BG}), we conclude by Lemma~\ref{l:etkan} that the map
    \[
        \tau_{<n+2}\cG^n\to P_1
    \]
    is an \'etale Kan fibration, and thus that $\tau_{<n+2}\cG^n\to \tau_{\le n+1}\cG^n$ is an \'etale hypercover of simplicial sheaves, as claimed. \hfill \qed

\section{Integrating Lie $n$-algebras}\label{s:lie3n}
As observed in \cite[Section 7]{Hen}, one defect of Sullivan's integration functor 
\[
L \mapsto \int L
\]
is that it takes values in Lie $\infty$-groups even when $L$ is an ordinary Lie algebra.  Naively, we might hope that given a finite type Lie $n$-algebra, there exists an integration functor taking values in Lie $n$-groups. However, by \cite[Theorem 7.5, Example 7.10]{Hen}, there exist nontrivial obstructions to the existence of such a functor. In this final section, we return to this and record what can salvaged.  

\begin{definition}
    A {\em Lie $n^\ast$-group/oid} is a Lie $(n+1)$-group/oid $X$ such that $X=\tau_{<n+1}X$ and for all $0\le i\le n+1$, the map
    \[
        \lambda^{n+1}_i\colon X_{n+1}\to \Lambda^{n+1}_iX
    \]
    is a local diffeomorphism. 
    
    More generally, a map $f\colon X\to Y$ is an $n^\ast$-stack if $f=\tau_{<n+1}(f)$ and for all $0\le i\le n+1$, the map
    \[
        \lambda^{n+1}_i(f)\colon X_{n+1}\to \Lambda^{n+1}_i(f)
    \]
    is a local diffeomorphism.
\end{definition}

\begin{remark}
        Lie $1^\ast$-groupoids are essentially {\em effective} Weinstein groupoids $\mathcal{G}/M$ in the sense of Tseng and Zhu \cite[Definition 1.1]{TZ}, i.e. Weinstein groupoids such that the differentiable stack $\mathcal{G}$ is equivalent to a sheaf. Concretely, by Zhu \cite[Theorem 1.5]{Z}, the nerve of an effective Weinstein groupoid is a Lie $1^\ast$-groupoid, and \cite[Section 4.2]{Z} gives a reverse construction (from Lie $1^\ast$-groupoids to weakly truncated Weinstein groupoids) such that the two constructions are mutually inverse up to equivalence. See Zhu \cite{Z} for more details.
\end{remark}

\begin{remark}
    Let $X$ be a Lie $\infty$-groupoid.  Recall that, from the perspective of higher compositions of higher morphisms, the spans
    \[
        \xymatrix{
            \Lambda^k_i X & \ar[l]_-{\lambda^k_i} X_k \ar[r]^-{d_i} & X_{k-1}
        }  
    \]
    encode the spaces of possible compositions of $(k-1)$-morphisms (and their inverses if $i=0,n+1$). In general, these compositions are uncountably, continuously multi-valued (i.e. if $\dim X_k>\dim \Lambda^k_iX$). 
    
    From this perspective, a Lie $n^\ast$-groupoid is a Lie $(n+1)$-groupoid such that a) the map of simplicial sheaves $X\to^\sim \tau_{\le n}X$ is a local weak equivalence, and b) the composition of $n$-morphisms in $X$ is {\em discretely} multi-valued.
\end{remark}

We will need the following.
\begin{lemma}\label{l:n*comp}\mbox{}
    \begin{enumerate}
        \item Let $f\colon X\to Y$ be an \'etale hypercover. Then $\tau_{<n+1}f\colon \tau_{<n+1}(X,f)\to Y$ is an $n^\ast$-stack.
        \item The composition of two $n^\ast$-stacks $X\to^f Y\to^g Z$ is an $n^\ast$-stack.
        \item Let $f\colon X\to Y$ be an $n^\ast$-stack, $g\colon Z\to Y$ a map of simplicial manifolds and suppose the pullback $X_0\times_{Y_0}Z_0$ exists in $\sm$. Then $X\times_Z Y\to Z$ is an $n^\ast$-stack in $\ssm$.
    \end{enumerate}
\end{lemma}
\begin{proof}
    The first item is an immediate consequence of Lemma~\ref{l:et} and the fact that \'etale hypercovers are \'etale Kan fibrations, e.g. by \cite[Lemma 3.12]{BG} applied to the descent category $\sSh(\sm)_{\et}$. 
    
    For the second, consider the diagram of pullback squares
    \[
        \xymatrix{
        \Csk_{n,k}(f) \ar[r] \ar[d] & Y_k \ar[d] \\
        \Csk_{n,k}(gf) \ar[r] \ar[d] & \Csk_{n,k}(g) \ar[d] \ar[r] & Z_k\ar[d] \\
        \Csk_{n,k}X \ar[r] & \Csk_{n,k}Y \ar[r] & \Csk_{n,k}Z
        }
    \]
    If $Y_k=\Image(\sigma_n^k(g))$, then $\Image(\sigma_n^k(f))=\Image(\sigma_n^k(gf))$. Therefore, $gf=\tau_{<n+1}gf$ (because $X_k=\Image(\sigma_n^k(f))$ by our assumption on $f$).  The condition on $(n+1)$-horns now follows from the diagram~\eqref{e:filler}, equivalently just as in the proof of \cite[Theorem 2.17(3)]{W}.  
    
    For the third statement, by \cite[Theorem 2.17(4)]{W}, the assumption on $X_0\times_{Y_0}Z_0$ guarantees that $\pi_Z\colon X\times_Y Z\to Z$ exists in $\ssm$. By Lemma~\ref{l:duskmoorepull}, that $f=\tau_{<n+1}f$ implies that $\pi_Z=\tau_{<n+1}\pi_Z$.  The condition on $(n+1)$-simplices now follows from the observation that the map $\lambda^{n+1}_i(\pi_Z)$ is isomorphic to the pullback of the map $\lambda^{n+1}_i(f)$ along $Z_{n+1}\to Y_{n+1}$ and the fact that local diffeomorphisms are preserved under pullback.
\end{proof}

We can now prove Lie's Third Theorem for Lie $n$-algebras.
\begin{proof}[Proof of Corollary~\ref{c:lien*}]
    Finite type Lie algebras have unobstructed integration (either by Lie's Third Theorem, or by \cite[Theorem 7.5]{Hen} and Bott's Theorem that $\pi_2G=0$ \cite[Theorem A]{Bott} for any finite dimensional Lie group $G$). Therefore, assume $n>1$, and let $L$ be a finite type Lie $n$-algebra.

    Let $\cG\xto{\sim} \sint L$ be a finite dimensional locally minimal model as in Theorem \ref{t:main4}. By \cite[Theorem 6.4]{Hen}, for all $i\ge 1$, $\pi_i \cG\cong \pi_i\sint L$ is a finite dimensional diffeological group, i.e. a quotient of a Lie group by a finitely generated subgroup. In particular, for each $i\ge 1$, we have an exact sequence for each $i\ge 1$
    \[
        0\to H_{i-1}L/\partial_i(\pi_{i+1}G)\to \pi_i \cG\to N\to 0
    \]
    where $G$ is the simply connected Lie group integrating $\g=H_0L$, and $N\subset\pi_iG$.  By \cite[Proposition 7.9]{Hen}, for all $m$, $\tau_{<m}\cG$ is a Lie $\infty$-group. 
    
    Now assume that the integration is unobstructed, i.e. that $\partial_n(\pi_{n+1}G)\subset H_{n-1}L$ is a discrete subgroup.  By the short exact sequence above, this implies that $\pi_n\cG$ is a finite dimensional Lie group. By Proposition~\ref{p:postbund}, the map
    \[
        \tau_{\le n}\int L\to\tau_{<n}\int L
    \]
    is a principal $K(\pi_n\int L,n)//G\to NG$ bundle in $\ssm$, and by Lemma~\ref{l:duskmoore} it is a minimal Kan fibration. Now consider the commuting square (in $\sSh(\sm)$ for the moment)
    \[
        \xymatrix{
            \tau_{\le n}\cG \ar[r]^\sim \ar[d] & \tau_{\le n} \int L\ar[d] \\
            \tau_{<n}\cG \ar[r]^\sim & \tau_{<n}\int L
        }.
    \]
    By Lemma~\ref{l:duskmoore} the vertical maps are minimal Kan fibrations, while the horizontal maps are a weak equivalences by Theorem~\ref{t:main4} and Lemma~\ref{l:tauexact}. We conclude by Lemma~\ref{l:minkan} that the square is a pullback. But, by \cite[Theorem 2.17(4)]{W}, we conclude that $\tau_{\le n}\cG$ is a Lie $n$-group. By Theorem \ref{t:main4} and Lemma~\ref{l:locmin}(\ref{l:locmindusk}), we conclude that $\tau_{\le n}\cG$ is a locally minimal model, as claimed; its finite dimensionality follows by inspection from the above, or alternately from Lemma~\ref{l:locminfd} and \cite[Theorem 6.4]{Hen}.  
    
Now suppose instead that the integration is obstructed, i.e. $\partial_n(\pi_{n+1}G)\subset H_{n-1}L$ is not a discrete subgroup. By Theorem~\ref{t:main4}, $\cG$ is locally minimal. By Lemma~\ref{l:locmin}(\ref{l:locminmoore}), we conclude that $\tau_{<n+1}\cG$ is a locally minimal Lie $\infty$-group.  It remains to show that it is a Lie $n^\ast$-group, i.e. that for all $i$, the map of smooth manifolds
    \begin{equation}\label{e:obstructed}
        \lambda^{n+1}_i\colon (\tau_{<n+1}\cG)_{n+1}\to\Lambda^{n+1}_i\tau_{<n+1}\cG
    \end{equation}
    is a local diffeomorphism. 

    For this, consider the pair of maps in $\sSh(\sm)$
    \[
        \tau_{<n+1}\cG\to^p\tau_{\le n}\cG\to \Delta^0.
    \]
    As in the proof of Lemma~\ref{l:kanngpd}, for all $i$, we have a commuting diagram of sheaves
    \begin{equation*}
        \xymatrix{
            (\tau_{<n+1}\cG)_{n+1} \ar[r]^{\lambda^{n+1}_i(p)} \ar[dr]_{\lambda^{n+1}_i} & \Lambda^{n+1}_i(p) \ar[d] \ar[r] & (\tau_{\le n}\cG)_{n+1} \ar[d]^{\lambda^{n+1}_i}_\cong \\
            & \Lambda^{n+1}_i\tau_{<n+1}\cG \ar[r]^{p_\ast} & \Lambda^{n+1}_i\tau_{\le n}\cG
        }
    \end{equation*}
    where the square is a pullback. Note that the right vertical map is an isomorphism (by definition), and thus the left vertical map is too, i.e. the map~\eqref{e:obstructed} is isomorphic to the map $\lambda^{n+1}_i(p)$. Further, because a map of representables is a surjective local diffeomorphism of sheaves if and only if it is represented by a surjective local diffeomorphism in the usual sense, it suffices to show that $\lambda^{n+1}_i(p)$ is a local diffeomorphism. By Lemma~\ref{l:locmin}(\ref{l:locmincrit}), the map
    \[
        p\colon \tau_{<n+1}\cG\to \tau_{\le n}\cG
    \]
    is an \'etale $(n+1)$-hypercover of simplicial sheaves. By \cite[Lemma 3.12]{BG}, this implies that $\lambda^{n+1}_i(p)$ is a surjective local diffeomorphism. We conclude that $\tau_{<n+1}\cG$ is a Lie $n^\ast$-group as claimed.
\end{proof}

\begin{example}\label{ex:s2g}
As recall in the introduction, the integration of the so-called string Lie 2-algebra was the motivating example for Henriques' work \cite{Hen}, and the construction of finite dimensional models was a motivating example for \cite{Z,SP,W}. We give a new, streamlined construction of such a model here, thereby proving Corollary \ref{cor:FinString}.  For the sake of exposition, we consider the level $k=1$ situation. A more detailed analysis of this construction will appear in forthcoming work.

Let $\g$ be a simple Lie algebra of compact type, with integrating Lie group $G$.  Let $\langle -,-\rangle$ denote the Killing form on $\g$, and let $\str(\g)$ denote the minimal Lie 2-algebra given by
    \begin{align*}
        \str(\g)_0&=\g\\
        \str(\g)_1&=\R
    \end{align*}
    with brackets
    \begin{align*}
        \ell_2((A,r),(B,s))&=([A,B]_{\g},0)\\
        \ell_3((A,r),(B,s),(C,t))&=(0,\langle [A,B]_{\g},C\rangle).
    \end{align*}
   The map $\str(\g)\to\g$ is a quasi-split fibration, so by \cite[Theorem 9.16]{RZ}, 
   \[
        \sint\str(\g)\to\sint\g
   \]
   is a Kan fibration with fiber $\sint \R[1]$.  By \cite[Theorem 7.5, Lemma 8.2]{Hen}, $\tau_{\le 2}\sint\str(\g)$ is a Lie 2-group, and as observed in \cite[Equation (49), p. 1041]{Hen}, the map 
   \[
    \tau_{\le 2}\sint\str(\g)\to \tau_{\le 2}\int\g
    \]
    is a Kan fibration with fiber $K(S^1,2)$.  By inspection, it is actually a minimal Kan fibration, and by Proposition~\ref{p:postbund}, it is a principal $K(S^1,2)$-bundle.  

   Now let $f_1\colon U_1\to G$ be a cover of $G$ by contractible opens, and fix a preimage $e\in U_1$ of the identity in $G$. We can view this is a map of 1-truncated reduced simplicial objects $U_1\to NG_{\le 1}$.  Let $U\to NG$ be the reduced, \'etale hypercover extending $U$, i.e.
   \[
        U_n=\csk_{1,n}U\times_{\csk_{1,n}NG}NG_n=U_1^{\binom{n+1}{2}}\times_{G^{\binom{n+1}{2}}}^{f,\sk_1}G^n
   \]
   By \cite[Theorem 6.4]{Hen}, Bott's Theorem that $\pi_2G=0$ \cite[Theorem A]{Bott}, and Lemma~\ref{l:duskmoore}, we have the isomorphism $\tau_{\le 2}\sint\g\cong\tau_{<2}\sint\g$. Therefore, by Theorem~\ref{t:splitexist}, the map $U$ sits in a commuting diagram
   \[
        \xymatrix{
            U \ar@{=}[d] \ar[r]^{\tilde{f}} & \tau_{\le 2}\sint\g\ar[d]^\sim\\
            U \ar[r]_\sim^f & NG
        }.
   \]
    Indeed, unpacking the proof of Theorem~\ref{t:splitexist}, we see that the map $\tilde{f}\colon U\to \tau_{\le 2}\sint\g$ in the square above is specified by a choice of lift
    \[
        \xymatrix{
       \{e\} \ar[d] \ar[r]^-{c_e} & \Map_\ast(\Delta^1,G) \ar[d]^{\ev_1}\\
        U_1 \ar[r]^{f_1} \ar@{..>}[ur]^{\tilde{f}_1} & G
        }
    \]
    where $c_e$ denotes the constant path at the identity.
    
    We now define 
    \[
        \cString^{\ft}(G,U_1,\tilde{f}):=U\times_{\tau_{\le 2}\sint\g}\tau_{\le 2}\sint\str(\g).
    \]
    The projection onto the second factor exhibits this is a finite dimensional, locally minimal model for the string 2-group of $G$, which depends only on the choice of open cover $f\colon U_1\to G$ and the lift $\tilde{f}$ above.
\end{example}

We close with the example of Henriques \cite[Example 7.10]{Hen} which motivated Corollary~\ref{c:lien*}.
\begin{example}\label{ex:h2g}
    Let $\g$ be a simple Lie algebra of compact type and let $\str(\g)$ be its string Lie 2-algebra as in the previous example. Let $(p,q)\in\R$ be a pair of $\Q$-linearly independent real numbers, and let $L_{\Hen}$ be the quotient Lie 2-algebra sitting in the pushout square
    \[
        \xymatrix{
            \R[1] \ar[r]^-{(p,q)} \ar[d] & \str(\g)\oplus\str(\g) \ar[d] \\
            0 \ar[r] & L_{\Hen} 
        }.
    \]
    Then, as shown in \cite[Example 7.10]{Hen},
    \begin{align*}
        \pi_1\int L_{\Hen}&=G\times G\\
        \pi_2\int L_{\Hen}&=(\R^2/(p,q)\R)/\Z^2\intertext{where the inclusion $\Z^2\into \R^2/(p,q)\R$ is given by the composite $\Z^2\into \R^2\to \R^2/(p,q)\R$, and for $i>2$}
        \pi_i\int L_{\Hen}&=0
    \end{align*}
    In particular, the integration is obstructed and $\tau_{\le 2}\int L_{\Hen}$ is not a simplicial manifold! 
    
    We construct a integrating finite dimensional Lie $2^\ast$-group here. As in the previous example, let $f_1\colon U_1\to G$ be a cover of $G$ by contractible opens, and let $f\colon U\to NG$ be the reduced, \'etale hypercover of $NG$ extending $f_1$. As in the proof of Proposition~\ref{p:hypexist}, let $h_2\colon V_2\to M_2(f)\cong U_2$ be a cover of $U_2$ by contractible opens such that the inclusion $s(U_1)\to U_2$ factors through $V_2$. We obtain a map of 2-truncated simplicial manifolds 
    \[
        h_{\le 2}\colon V_{\le 2}\to NG_{\le 2}
    \]
    and let $h\colon V\to NG$ be the reduced, \'etale hypercover extending $h_{\le 2}$.  By construction, $h$ is good up to level 2. By Theorem~\ref{t:splitexist} and Bott's Theorem as in the previous example, we have a commuting diagram
     \[
        \xymatrix{
             && \tau_{<3}\sint(\g\times\g)\ar[d]^\sim\\
            V\times V \ar[rr]_\sim^{h\times h} \ar[urr]^{\tilde{h}} && NG\times NG
        }.
   \]
   Unpacking the proof of Theorem~\ref{t:splitexist} as in the previous example, we see that the map $\tilde{h}\colon V\times V\to \tau_{<3}\sint(\g\times\g)$ is specified by a choice of lift $\tilde{f}_1$
    \[
        \xymatrix{
        \{e\} \ar[d] \ar[r]^-{c_e}        & \Map_\ast(\Delta^1,G\times G) \ar[d]^{\ev_1}\\
        U_1\times U_1 \ar[r]^f \ar@{..>}[ur]^{\tilde{f}_1} & G\times G
        }
    \]
    as in the previous example, and a choice of lift
    \[
        \xymatrix{
            s(U_1\times U_1) \ar[d] \ar[rr]    && \Map_\ast(\Delta^2,G\times G) \ar[d]^{\partial^\ast}\\
            V_2\times V_2 \ar[rr]_-{\tilde{f}_1^{\times 3}\circ\partial^\ast\circ h} \ar@{..>}[urr]^-{\tilde{h}_2} &&\Map_\ast^0(\partial\Delta^2,G\times G)
        }.
    \]
    Now let $\tilde{V}$ by the pullback
    \[
        \xymatrix{
        \widetilde{V} \ar[r] \ar[d] & \tau_{\le 3}\int(\g\times\g)\ar[d] \\
        V\times V \ar[r] & \tau_{<3}\int(\g\times \g)
        }.
    \]
    As observed in the proof of Theorem~\ref{t:splitexist} (i.e. by \cite[Example 5.5]{Hen}, Lemma~\ref{l:duskmoore} and Proposition~\ref{p:postbund}), the map
    \[
        \widetilde{V}\to V\times V
    \]
    is a minimal Kan fibration and a principal $K(\Z^2,3)$-bundle. 
    
    Now denote by $p\colon L_{\Hen}\to\g\times\g$ the quotient of $L_\infty$-algebras.  This is a strict, quasi-split fibration, so by integrating, we obtain by \cite[Theorem 9.16]{RZ}, a Kan fibration 
    \[
        \sint p\colon \sint L_{\Hen}\to\sint(\g\times\g)
    \]
    As in the proof of Theorem~\ref{t:main3}, by Proposition~\ref{p:postbund}, 
    \[
        \tau_{\le 2}(\sint L_{\Hen},\sint p)\to \sint(\g\times\g)
    \]
    is a minimal Kan fibration and a principal $K(\R,2)$-bundle. Moreover, by Lemma~\ref{l:postsplit}, the hypercover
    \[
        \sint L_{\Hen}\to \tau_{\le 2}(\sint L_{\Hen},\sint p)
    \]
    admits a section $\sigma$.

    Now, viewed as a Lie 3-algebra, the integration of $L_{\Hen}$ is unobstructed, i.e. by \cite[Theorem 7.5]{Hen} and Lemma~\ref{l:tauexact}, we can apply $\tau_{\le 3}$ to obtain a Kan fibration in $\ssm$
    \[
        \tau_{\le 3}\sint p\colon \tau_{\le 3}\sint L_{\Hen}\to \tau_{\le 3}\sint(\g\times\g).
    \]
    By \cite[Theorem 6.4]{Hen} and Lemma~\ref{l:trtr}, we obtain a canonical isomorphism over $\tau_{\le 3}\sint\g$
    \[
        \tau_{\le 3}(\tau_{\le 2}(\sint L_{\Hen},\sint p))\cong \tau_{\le 2}(\tau_{\le 3}\sint L_{\Hen},\tau_{\le 3}\sint p).
    \]
    In particular, $\tau_{\le 3}\sint p$ factors as
    \[
        \tau_{\le 3}\sint L\to \tau_{\le 2}(\tau_{\le 3}\sint L,\tau_{\le 3}\sint p)\to \tau_{\le 3}\sint(\g\times\g)
    \]
    where the first map is a split hypercover and the second is a principal $K(\R,2)$-bundle.

    Let $\widetilde{V}'$ be the pullback
    \[
        \widetilde{V}':=\widetilde{V}\times_{\tau_{\le 3}\sint(\g\times\g)}\tau_{\le 2}(\tau_{\le 3}\sint L,\tau_{\le 3}\sint p).
    \]
    Post-composing with $\tau_{\le 3}\sigma$, we obtain a weak equivalence
    \[
        \tau_{\le 3}\sigma\colon \widetilde{V}'\to^\sim \tau_{\le 3}\sint L_{\Hen}   
    \]
    realizing $\widetilde{V}'$ as a locally minimal model for $\tau_{\le 3}\sint L_{\Hen}$.  Moreover, since $\pi_3\sint L_{\Hen}=0$ (by \cite[Example 7.10]{Hen}), by Lemma~\ref{l:duskmoore}, the map 
    \[
        \tau_{\le 3}\sint L_{\Hen}\to\tau_{<3}\sint L_{\Hen}
    \]
    is a minimal fibration and a weak equivalence.  By Lemma~\ref{l:minkan}, it is an isomorphism, and by the same logic, so is the map
    \[
        \tau_{\le 3}\widetilde{V}'\to^\cong \tau_{<3}\widetilde{V}'.
    \]
    We see that
    \[
        \cG_{\Hen}(G,U_1,V_2,\tilde{f_1},\tilde{h}_2,\sigma):=\tau_{\le 3}\widetilde{V'}\to^{\tau_{\le 3}\sigma}_\sim \tau_{\le 3}\sint L_{\Hen}\cong \tau_{<3}\sint L_{\Hen}
    \]
    exhibits $\cG_{\Hen}(G,U_1,V_2,\tilde{f_1},\tilde{h}_2,\sigma)$ as a locally minimal finite dimensional Lie $2^\ast$-group integrating $L_{\Hen}$; it depends only on the choices $f_1$, $h_2$, and $\sigma$, and the lifts $\tilde{f}_1$ and $\tilde{h}_2$ above.
\end{example}

\appendix

\section{Formulas for $k$-invariants of Kan fibrations of Lie $\infty$-groups}\label{app:kinvar}
In this section, we give the proof of Theorem~\ref{t:kinvar}. We briefly recall the context of the theorem.  We work in a category with covers $\C$.  We are given a smooth $\infty$-group $Y\in\sC$ such that $\pi_1 Y$ is a group in $\C$, and a Kan fibration $f\colon X\to Y$ such that $\tau_{<n}f$ is an isomorphism for some $n>1$. We write $F:=X\times_Y \Delta^0$ and we assume that $\pi_n F$ is a group object in $\C$. We now give the proof as follows.
    
\begin{proof}[Proof of Theorem~\ref{t:kinvar}]
    {\bf Step 1: Construct the map $\varphi$.} First, observe that our assumptions on $f$ imply that $f$ is an isomorphism on $(n-1)$-skeleta, and that for all $k\ge n-1$, the maps $\sigma_{n-1}^k(f)$ induce isomorphisms
    \[
        \Image(\sigma_{n-1}^k(f))\to^\cong Y_k.
    \]
    Next, because $Y_0=\ast$ and $f$ is Kan, we have a commuting diagram
    \[
        \xymatrix{
            F\times Y_n \ar[rr]^-{\sim} \ar[d] && (\Delta^n\times Y_n)\times_Y X \ar[r] \ar[d] & X \ar[d]^f \\
            \Delta^0\times Y_n \ar[rr]_-{\sim}^-{(\delta^0)^n\times 1} && \Delta^n\times Y_n \ar[r]^-{\ev} & Y
        }
    \]
    in which both squares are pullbacks, all vertical maps are Kan fibrations, and the horizontal maps in the left square are both local weak equivalences (the bottom by inspection, the upper by Lemma~\ref{l:rightproper}).  Passing to homotopy sheaves over $Y_n$, we obtain a $\pi_1Y$-equivariant isomorphism of $Y_n$-groups
    \[
        \xymatrix{
            \pi_n F\times Y_n \ar[dr] \ar[rr]^-\cong && \pi_n((\Delta^n\times Y_n)\times_Y X) \ar[dl] \\
            & Y_n
        } 
    \]
    where the left diagonal arrow is just the projection, and the right diagonal arrow sends a based homotopy class of $n$-spheres to the basepoint.

    Denote by $\nu$ the composition
    \[
        \xymatrix{
            \hom_{Y_n}((\Delta^n\cup_{\partial\Delta^n}\Delta^n)\times Y_n,(\Delta^n\times Y_n)\times_Y X) \ar[rr]^-\nu \ar[d] && \pi_n F \\ \pi_n((\Delta^n\times Y_n)\times_Y X) \ar[rr]^-\cong && \pi_nF\times Y_n \ar[u]
        }
    \]
    We are now in a position to construct the map 
    \[
        \varphi\colon \LW(X/Y)\to K(\pi_nF,n+1)//\pi_1Y.
    \]
    As observed in Remark~\ref{r:hoquot}, \cite[Theorem 6.7]{W} and Lemma~\ref{l:stackcosk} imply that $K(\pi_nF,n+1)//\pi_1Y\to N\pi_1 Y$ is a smooth $(n+1)$-stack and is $(n+2)$-coskeletal. Therefore it suffices to define the map $\varphi$ on $(n+2)$-skeleta. 
    
    From the definition of $\LW$ and our assumption on $f$, 
    \begin{align*}
        \LW(X/Y)_{\le n}&\cong X_{\le n}\\ 
        \LW(X/Y)_{n+1}&\cong X_{n+1}\times_{\Image(\mu_n(f))}^{\mu_n(f)d_0,\mu_n(f)}X_n\\
        \LW(X/Y)_{n+2}&\cong X_{n+2}\times_{\Image(\sigma_{n-1}^{n+1}(f))}^{\sigma_{n-1}^{n+1}(f)d_0,\sigma_{n-1}^{n+1}(f)}X_{n+1}\times_{\Image(\mu_n(f))}^{\mu_n(f)d_0,\mu_n(f)}X_n
    \end{align*}
    with face and degeneracy maps as in \cite[Example 3.7]{W2}. 
    
    For the target, we consider $K(A,n)$ with ``inhomogeneous coordinates'' as in Definition~\ref{d:Ainhomog} and use the isomorphism of Lemma~\ref{l:relw=k[1]}
    \[
        K(\pi_nF,n+1)//\pi_1Y\to^\cong \LW((K(\pi_nF,n)//\pi_1Y)/N\pi_1Y)
    \]
    to obtain 
    \begin{align*}
        (K(\pi_nF,n+1)//\pi_1Y)_{\le n}&= (N\pi_1 Y)_{\le n}\\
        (K(\pi_nF,n+1)//\pi_1Y)_{n+1}&=(\pi_1Y)^{n+1}\times \pi_n F\\
        (K(\pi_nF,n+1)//\pi_1Y)_{n+1}&=(\pi_1Y)^{n+2}\times (\pi_nF)^{n+1}\times\pi_n F
    \end{align*}
    with face and degeneracy maps determined from those of $K(\pi_nF,n)//\pi_1Y$ as in \cite[Example 3.6]{W2}. We emphasize that this precise choice of coordinates is necessary for the formulas below to satisfy the simplicial identities, and thus define simplicial maps.
    
    We now define $\varphi$ on $k$-simplices for $k\le n$ by
    \begin{align*}
       \LW(X/Y)_k\cong X_{\le k}&\to^{\varphi_k} (\pi_1Y)^k=(K(\pi_nF,n+1)//\pi_1Y)_k\\
        x&\mapsto [\spine f(x)]
    \end{align*}
    where $[-]$ denotes based homotopy class.  Next, observe that the assignment
    \[
        (x_{n+1},x_n)\mapsto (d_0x_{n+1},x_n)
    \]
    defines a map over $Y_n$
    \[
       \LW(X/Y)_{n+1}\to \hom_{Y_n}((\Delta^n\cup_{\partial\Delta^n}\Delta^n)\times Y_n,(\Delta^n\times Y_n)\times_Y X).
    \]
    We define $\varphi$ on $(n+1)$-simplices by
    \begin{align*}
        \LW(X/Y)_{n+1}&\to^{\varphi_{n+1}} (\pi_1Y)^{n+1}\times\pi_n F\\
         (x_{n+1},x_n)&\mapsto ([\spine f(x_{n+1})],\nu(d_0x_{n+1},x_n))
    \end{align*}
    Similarly, observing that the assignments $(x_{n+2},x_{n+1},x_n)\mapsto (d_id_0x_{n+2},d_ix_{n+1})$ and $(x_{n+2},x_{n+1},x_n)\mapsto (d_0x_{n+1},x_n)$ define maps over $Y_n$
    \[
        \LW(X/Y)_{n+2}\to \hom_{Y_n}((\Delta^n\cup_{\partial\Delta^n}\Delta^n)\times Y_n,(\Delta^n\times Y_n)\times_Y X),
    \]
    we define $\varphi$ on $(n+2)$-simplices by
    \begin{align*}
         \LW(X/Y)_{n+2}&\to^{\varphi_{n+2}} (\pi_1Y)^{n+2}\times (\pi_nF)^{n+1}\times\pi_n F\\
         (x_{n+2},x_{n+1},x_n)&\mapsto ([\spine f(x_{n+2})],\{\nu(d_id_0 x_{n+2},d_i x_{n+1})\}_{i=0}^n,\nu(d_0x_{n+1},x_n)).
    \end{align*} 
    Using our choice of inhomogenous coordinates on $K(\pi_nF,n)$ as noted above, all of the simplicial identities now follow immediately from the formulas in \cite[Examples 3.6, 3.7]{W2}, except the identity
    \[
        d_{n+2}\circ\varphi_{n+2}=\varphi_{n+1}\circ d_{n+2}.
    \]
    Unpacking this and noting that $d_{n+2}\spine f(x_{n+2})=\spine f(d_{n+2}x_{n+2})$, we are reduced to verifying that
    \[
        [\spine f(x_{n+2})_{n+2}]\cdot \left((-1)^n\sum_{i=0}^{n}(-1)^i\nu(d_id_0x_{n+2},d_ix_{n+1})\right)=\nu(d_0d_{n+2}x_{n+2},d_{n+1}x_{n+1}).
    \]
    But, this follows from the corresponding formula for simplicial sets, i.e. the observation that for any reduced Kan complex $S$, a map
    \[
        (x,y)\colon \Delta^{n+1}\cup_{\sk_{n-1}\Delta^{n+1}}\Delta^{n+1}\to S
    \]
    gives rise to a relation
    \[
        [\spine(x)_{n+1}]\cdot\left((-1)^n\sum_{i=0}^n (-1)^i[(d_ix,d_iy)]\right)=[(d_{n+1}x,d_{n+1}y)]\in \pi_nS. 
    \]
    This completes the construction of the map $\varphi$.

    \smallskip
    \noindent
    \newline
    {\bf Step 2: Construct the map $\psi$.}
    We now construct the map 
    \[
        \psi\colon X\times_Y\LW(X/Y)\to WK(\pi_n F,n)//\pi_1Y.
    \]
    Calculating as above, we have
    \begin{align*}
        (X\times_Y\LW(X/Y))_{\le n-1}&\cong\LW(X/Y)_{\le n-1}\\
        (X\times_Y\LW(X/Y))_n&\cong X_n\times_{\Image(\mu_n(f))}X_n\\
        (X\times_Y\LW(X/Y))_{n+1}&\cong X_{n+1}\times_{\Image(\sigma_{n-1}^{n+1}(f))}X_{n+1}\times_{\Image(\mu_n(f))}^{\mu_n(f)d_0,\mu_n(f)}X_n\\
        (X\times_Y\LW(X/Y))_{n+2}&\cong X_{n+1}\times_{\Image(\sigma_{n-1}^{n+2}(f))}X_{n+2}\times_{\Image(\sigma_{n-1}^{n+1}(f))}^{\sigma_{n-1}^{n+1}(f)d_0,\sigma_{n-1}^{n+1}(f)}X_{n+1}\times_{\Image(\mu_n(f))}^{\mu_n(f)d_0,\mu_n(f)}X_n 
    \end{align*}
    with face and degeneracy maps as in \cite[Example 3.7]{W2}.  Similarly, using the isomorphism of Lemma~\ref{l:relw=k[1]}
    \[
        WK(\pi_nF,n)//\pi_1Y\to^\cong W((K(\pi_nF,n)//\pi_1Y)/N\pi_1Y)
    \]
    we have
    \begin{align*}
        (WK(\pi_nF,n)//\pi_1Y)_{\le n-1}&\cong (N\pi_1Y)_{\le n-1}\\
        (WK(\pi_nF,n)//\pi_1Y)_n&\cong (\pi_1Y)^n\times\pi_nF\\
        (WK(\pi_nF,n)//\pi_1Y)_{n+1}&\cong (\pi_1Y)^{n+1}\times(\pi_nF)^{n+1}\times\pi_nF\\
        (WK(\pi_nF,n)//\pi_1Y)_{n+2}&\cong (\pi_1Y)^{n+1}\times(\pi_n F)^{\binom{n+2}{2}}\times(\pi_nF)^{n+1}\times\pi_nF
    \end{align*}
    with face and degeneracy maps determined from those on $K(\pi_nF,n)//\pi_1Y$ as in \cite[Lemmas 3.3, 3.4]{W2}. We now define $\psi$ as follows:
    \begin{align*}
        X\times_Y\LW(X/Y)_{\le n-1}&\to^{\psi_{\le n-1}} (N\pi_1Y)_{\le n-1}\cong (WK(\pi_nF,n)//\pi_1Y)_{\le n-1}\\
        x&\mapsto [\spine f(x)]\\
        \\
        (X\times_Y\LW(X/Y))_n&\to^{\psi_n} (\pi_1Y)^n\times\pi_nF\\
        (x,y)&\mapsto ([\spine f(x)],\nu(x,y))\\
        \\
        (X\times_Y\LW(X/Y))_{n+1}&\to^{\psi_{n+1}} (\pi_1Y)^{n+1}\times(\pi_nF)^{n+1}\times\pi_nF\\
        (x_{n+1},y_{n+1},x_n)&\mapsto ([\spine f(x_{n+1})],\{\nu(d_ix_{n+1},d_iy_{n+1})\}_{i=0}^n,\nu(d_0y_{n+1},x_n))
    \end{align*}
    By the same arguments as above, this defines a map of $(n+1)$-truncated objects over $\varphi$, i.e. we have a commuting diagram
    \[
        \xymatrix{
            (X\times_Y \LW(X/Y))_{\le n+1} \ar[rr]^-{\psi_{\le n+1}} \ar[d] && (WK(\pi_nF,n)//\pi_1Y)_{\le n+1} \ar[d]\\
            \LW(X/Y) \ar[rr]^-{\varphi} && K(\pi_nF,n+1)//\pi_1Y
        }
    \]
    By \cite[Theorem 6.7]{W}, the right vertical map is an $n$-stack, and thus $(n+1)$-coskeletal (by Lemma~\ref{l:stackcosk}). We conclude that $\psi_{\le n+1}$ extends uniquely to a map 
    \[
        \psi\colon X\times_Y\LW(X/Y)\to WK(\pi_nF,n)//\pi_1Y
    \]
    over $K(\pi_nF,n+1)//\pi_1Y$ as claimed. 

    That the restriction of $\psi$ to the fibers induces an isomorphism on $\pi_1$ follows directly from the formulae above.  Indeed, for 
    \[
         (x,s_0^n\ast)\in\hom(\Delta^n/\partial\Delta^n,X\times_Y\LW(X/Y)\times_{\LW(X/Y)}\Delta^0)\cong\hom(\Delta^n/\partial\Delta^n,F)
    \]
    the formulas above give
    \[  
        \psi(x,s_0^n\ast)=\nu(x,s_0^n\ast)\in\pi_nF,
    \]
    i.e. on $\pi_n$, the restriction of $\psi$ to the fibers is isomorphic to the identity map
    \[
        \pi_nF\to^1_{\cong}\pi_nF\cong \pi_n K(\pi_nF,n)
    \]
    as claimed.

    \smallskip
    \noindent
    \newline
    {\bf Step 3: Construct the pullback square~\eqref{e:kinvar2}.} 
    By \cite[Theorem 6.7]{W}, the right vertical map in the square~\ref{e:kinvar1} is a $n$-stack, and by inspection, it is minimal. By Lemma~\ref{l:trtr}, the square~\eqref{e:kinvar1} induces a commuting square
    \[
        \xymatrix{
         \tau_{\le n}(X\times_Y \LW(X/Y),f|_{\LW(X/Y)}) \ar[r] \ar[d] & WK(\pi_nF,n)//\pi_1 Y \ar[d] \\
         \LW(X/Y) \ar[r]^-{\varphi} & K(\pi_nF,n+1)//\pi_1 Y
        }
    \]  
    By Lemmas~\ref{l:trtr} and~\ref{l:minpost}, our assumptions on $f$ guarantee that
    \[
        \tau_{\le n}(X\times_Y\LW(X/Y),f|_{\LW(X/Y)})\to\LW(X/Y)
    \]
    is a minimal Kan fibration. By Lemma~\ref{l:minkan} and the conclusion of Step 2 above, we obtain a weak equivalence, and thus an isomorphism of minimal Kan fibrations over $\LW(X/Y)$
    \[
        \tau_{\le n}(X\times_Y\LW(X/Y),f|_{\LW(X/Y)})\to^\cong \LW(X/Y)\times_{K(\pi_nF,n+1)//\pi_1Y}WK(\pi_nF,n)//\pi_1Y.
    \]
    By the same reasoning, the canonical map gives an isomorphism over $\LW(X/Y)$
    \[
        \tau_{\le n}(X\times_Y\LW(X/Y),f|_{\LW(X/Y)})\to^\cong \tau_{\le n}(X,f)\times_Y\LW(X/Y).
    \]
    Together, these isomorphisms determine an isomorphism over $\LW(X/Y)$
    \[
        \tau_{\le n}(X,f)\times_Y\LW(X/Y)\to^\cong \LW(X/Y)\times_{K(\pi_nF,n+1)//\pi_1Y}WK(\pi_nF,n)//\pi_1Y
    \]
    and thus the pullback square~\eqref{e:kinvar2} as claimed.
\end{proof}

\section{$L_\infty$-algebras and coalgebra bundles}\label{app:Linf}
\stoptocentries
\subsection{$L_\infty$-algebras as dg coalgebras} \label{sec:Linf-coalg}
In order to establish notation and conventions, we first summarize the dg coalgebraic characterization of $L_\infty$-algebras and weak $L_\infty$-morphisms. We follow the exposition in \cite[\Sec 2.4]{R} and work over a field $\kk$ of characteristic zero. For a graded vector space $V$, we denote by  $\S(V):=V \oplus S^{2}(V) \oplus S^{3}(V) \oplus \cdots$
the reduced conilpotent cocommutative coalgebra cofreely generated by $V$.  Let $\Phi \maps \S(V) \to \S(V')$ be a linear map. For $p,m \geq 1$ the  notation $\Phi^p_{n}$ is reserved for the restriction-projections $\Phi^p_{m}:= \pr_{\S^{p}(V')} \circ \Phi \vert_{\S^{m}(V)}$. Furthermore, we denote by $\Phi^1 \maps \S(V) \to V'$ the linear map $\Phi^{1} := \Phi^{1}_{1} + \Phi^{1}_{2} + \cdots$.

In particular, a degree zero linear map $F^{1} \maps \S(V) \to V'$ uniquely determines  a coalgebra morphism $F \maps \S(V) \to \S(V')$ via the following formula
\begin{equation} \label{eq:morph_formula1}
\begin{split}
F(v_{1}, \ldots, v_{m}) &=  F^{1}_{m}(v_{1}, \ldots,
 v_{m}) +  \sum_{p=2}^m F^p_m(v_{1}, \ldots,
 v_{m}), 
\end{split}
\end{equation}
where
\begin{equation} \label{eq:morph_formula2}
\begin{split}
F^{p}_{m}(v_{1}, \ldots, v_{m}) &=
\sum^{k_{1} + k_{2} +
  \cdots + k_{p}=m}_{k_{1},k_{2},\hdots,k_{p} \geq 1} \sum_{\sigma
  \in \Sh(k_{1},k_{2},\hdots,k_{p})} \frac{\epsilon(\sigma)}{p!}
F^{1}_{k_{1}}(v_{\sigma(1)}, \ldots, v_{\sigma(k_{1})}) \\
& \quad \odot F^{1}_{k_{2}}(v_{\sigma(k_{1} + 1)}, \ldots, v_{\sigma(k_{1}+k_{2})}) \odot \cdots 
\odot F^{1}_{k_{p}}(v_{\sigma(m-k_{p} + 1)}, \ldots, v_{\sigma(m)}).
\end{split}
\end{equation}
The composition $GF \maps \S(V) \to \S(V'')$ of coalgebra morphisms $F \maps \S(V) \to \S(V')$  and $G \maps \S(V') \to \S(V'')$ is the unique coalgebra morphism whose structure maps $(GF)^1_m \maps \S^m(V) \to V''$ are
\begin{equation}\label{eq:comp}
\begin{split}
(GF)^1_m(v_1,\ldots, v_m) & = \sum^m_{p=1} G^1_p F^p_m(v_1,\ldots, v_m)
\end{split}
\end{equation} 
Analogously, a degree $-1$ linear map $\delta^{1} \maps \S(V) \to V$ uniquely determines a degree $-1$ coderivation $\delta \maps \S(V) \to \S(V)$ via the following formula
\begin{equation} \label{eq:coder_formula1}
\delta_{m}(v_{1}, \ldots, v_{m}) = 
\delta^1_{m}(v_{1},\ldots,v_{m}) + \sum_{p=2}^m \delta ^p_m(v_{1}, \ldots, v_{m})
\end{equation}
where
\begin{equation} \label{eq:coder_formula2}
\begin{split}
\delta^{p}_{m}(v_{1}, \ldots, v_{m}) = \hspace{-.5cm}\sum_{\sigma \in \Sh(m-p+1,p-1)} \hspace{-.5cm}
\epsilon(\sigma) \delta^1_{m-p+1}(v_{\sigma(1)}, \ldots, v_{\sigma(m-p+1)}) \odot v_{\sigma(m-p+2)} \odot \cdots \odot v_{\sigma(m)}.
\end{split}
\end{equation}
Recall that an $L_\infty$-algebra structure $\el = \el_1,\el_2,\cdots$ on $L$ is  equivalent to equipping the coalgebra $\bar{S}(\bs L)$ with a codifferential $\delta$ of degree $-1$. The correspondence between the structure maps 
\[
\delta^1_m \maps \S^{m}(\bs L) \to \bs L
\]
and the brackets $\el$ is given by the formula
\begin{equation}\label{eq:struc_skew}
\delta^1_{m} = (-1)^{\frac{m(m-1)}{2}} \bs \circ
\el_{m} \circ {(\ds)}^{\tensor m}.
\end{equation}
The dg coalgebra $(\bar{S}(\bs L), \del)$ is often referred to as the (reduced) \df{Chevalley-Eilenberg coalgebra} of $(L,\el)$. 
Similarly, weak $L_\infty$-morphisms $f \maps (L,\el) \to (L',\el')$ are in one to one correspondence with degree 0 morphisms of reduced dg coalgebras $F \maps (\bar{S}(\bs L), \del) \to (\bar{S}(\bs L'),\del')$. The degree $0$ structure maps $F^1_k \maps \S^k(\bs L) \to \bs L$ are given by the formula
\begin{equation} \label{eq:morph_eq1}
 F^{1}_{k} = (-1)^{\frac{k(k-1)}{2}} \bs \circ
 f_{k} \circ {(\ds)}^{\tensor k}.
\end{equation}

\subsection{Classification of coalgebra bundles}
We denote by $\dgco$, the category of homologically $\Z$-graded, counital coaugmented conilpotent cocommutative dg coalgebras. Since we work over a field, such coalgebras are connected, in the sense of Quillen \cite[\Sec B.3]{Quillen}. Given $C \in \dgco$, we let $\ba{C}$ denote its coaugmentation ideal, so that $C = \kk \dsum \ba{C}$ as graded vector spaces. Note that $\ba{C}$ is a reduced conilpotent dg coalgebra, in agreement with our notation above for the reduced Chevalley-Eilenberg coalgebra $(\bar{S}(\bs L), \del)$ of an $L_\infty$-algebra. Indeed, the coderivation $\del$ extends to
a coderivation on the coaugmented counital cofree conilpotent graded coalgebra ${S}(\bs L) = \kk \dsum \ba{S}(\bs L)$ by setting $\del(1):=0$. Similarly, every reduced dg coalgebra morphism $F \maps (\bar{S}(\bs L), \del) \to (\bar{S}(\bs L'),\del')$ extends to a morphism $F \maps ({S}(\bs L), \del) \to ({S}(\bs L'),\del')$ in $\dgco$ by setting $F(1):=1$. Hence, from here on, we will tacitly identify the  category of Lie $\infty$-algebras as a full subcategory $\LnA{\infty} \sse \dgco$.

\begin{definition}[Def.\ 2.2.1 \cite{Prigge}] \label{def:coalg-bun}
A {\em trivialized coalgebra bundle} $(E,\pi,\vphi)$ over a base $B$ with fiber $C \in \dgco$ is a morphism $\pi \maps E \to B$ in $\dgco$ equipped with a ``local trivialization'', i.e.\ an isomorphism of graded coalgebras $\vph \maps B \tensor C \xto{\cong} E$ such that $\pr_B = \pi \cc \vphi$, and
such that $\vphi\vert_{1_B \tensor C} \maps C \to E$ is a morphism in $\dgco$.
\end{definition}
Throughout, we will frequently refer to trivialized coalgebra bundles simply as ``coalgebra bundles'' whenever it is convenient to do so. Note that the local trivialization is not necessarily compatible with the codifferentials.\footnote{This is parallel to how in twisted Cartesian products of simplicial objects, the trivializing isomorphisms are not required to respect the zeroth face maps.} Morphisms between coalgebra bundles $(E,\pi,\vphi)$ and $(E',\pi',\vphi')$ with base $B$ are defined to be those morphisms in the slice category $\dgco_{/B}$ which intertwine the local trivializations.   

\begin{example}\label{ex:cb1}
We will maintain a running example throughout this subsection, in order to set the stage for the proof of Prop.\ \ref{prop:classify}. Let $A$ be a non-graded vector space, and $m \geq 1$. Let $g=g_1 \maps (L,\el) \fib (L',\el')$  be a minimal fibration of Lie $\infty$-algebras which is also surjective in degree 0, and such that $\ker g = A[m]$ is $A$ concentrated in degree $m$. Then the dg coalgebra morphism $G \maps  
({S}(\bs L), \del) \to ({S}(\bs L'),\del')$ corresponding to $g$ is a coalgebra bundle over the base $(S(\bs L'),\del')$ with fiber $(S \bigl (\bs A[m] \bigr),0)$. Indeed, choose a section of $g$ and let $\phi \maps L' \dsum A[m] \xto{\cong} L$ denote the corresponding isomorphism of graded vector spaces. Then
$\phi$ induces a local trivialization
\begin{equation} \label{eq:coalgbun-split}
S(\bs L') \tensor S \bigl(\bs A[m]\bigr) \cong S \bigl(\bs L' \dsum \bs A[m] \bigr)  \xto{\vphi} S(\bs L).    
\end{equation}   
Explicitly, for $x_1,\ldots, x_k \in L'$ and $a_1,\ldots,a_k \in A$: 
\begin{equation} \label{eq:cb-split}
\begin{split}
\vphi\bigl( (\bs x_1,\bs^{m+1} a_1) \vee (\bs x_2,\bs^{m+1} a_2) & \vee \cdots \vee (\bs x_k,\bs^{m+1} a_k) \bigr) =\\
& \bs \phi(x_1,\bs^{m} a_1) \vee \bs \phi(x_2,\bs^{m} a_2) \vee \cdots \vee \bs \phi(x_k,\bs^{m} a_k).   
\end{split}
\end{equation}
\end{example}

Fix $B, C \in \dgco$. Denote by $\Bun_B(C)$ the set of isomorphism classes of trivialized coalgebra bundles over $B$ with fiber $C$. We now recall the ingredients needed for the classification of such bundles in terms of twisting functions, as developed by Prigge in \cite[\Sec 2.2]{Prigge}, which builds on Quillen's classification of principal dg coalgebra bundles in \cite[\Sec B.5]{Quillen}. 

\subsubsection{Coderivation dg Lie algebra}
First, let $\g_C:= \Coder(C)$ denote the subspace of the hom chain complex $\Hom_{\kk}(C,C)$ consisting of all codervations of $C$ as a graded coalgebra. The restriction of the commutator bracket and differential give $\g_C$ the structure of a $\Z$-graded dg Lie algebra.  

\begin{example}\label{ex:cb2}
Keeping the notation of Example \ref{ex:cb1}, consider the coderivation Lie algebra for $C= S(\bs A[m])$. Since $S(\bs A[m])$ is cofree, it is well known (e.g.\ \cite[\Sec 4.2]{DR}) that we may identify $\g_{S(\bs A[m])}$ with the hom complex $\Hom_{\kk}(S(\bs A[m]), \bs A[m])$ equipped with the ``convolution bracket''. The identification has an explicit description via a formula similar to the one for degree $-1$ coderivations given above in Eq.\ \ref{eq:coder_formula1}. Since $\bs A[m]$ is concentrated in degree $m+1$, the underlying graded vector space of $\g_{S(\bs A[m])}$ is concentrated in degrees 0 and $m+1$. Explicitly, we have
\[
\begin{split}
&{\g_{S(\bs A[m])}}_{0} = \Hom_{\kk}(\bs A[m], \bs A[m]) \cong \End(A), \\  
&{\g_{S(\bs A[m])}}_{m+1} = \Hom_{\kk}(\kk, \bs A[m]) = A[m+1].  
\end{split}
\]   
Using the above identifications, it is straightforward to show that they induce an isomorphism of dg Lie algebras $\g_{S(\bs A[m])} \cong A[m+1] // \End(A) = \LB{m}$.
\end{example}
\subsubsection{Universal coalgebra bundle}  
Let $(\bs \g_C \dsum \g_C, d, [\cdot,\cdot])$ denote the acyclic dg Lie algebra whose underlying graded vector space is the direct sum $\bs \g_C \dsum \g_C$, and whose dg Lie structure is uniquely determined by the equations: $d(\bs D)= D$, $[D, \bs D'] = \bs[D,D']_{\g_C}$, $[\bs D, \bs D']=0$ for all coderivations $D,D' \in \g_C$, along with the requirement that the inclusion $\g_C \sse \bs \g_C \dsum \g_C$ is a morphism of dg Lie algebras. This inclusion, in fact, gives the universal enveloping algebra $U(\bs \g_C \dsum \g_C)$ the structure of coalgebra in right $U(\g_C)$-modules.  
Tensoring with the counit $U(\g_C) \to \kk$ over $U(\g_C)$ gives a surjective morphism 
\[
\pi_{U(\g_C)} \maps U(\bs \g_C \dsum \g_C) \to U(\bs \g_C \dsum \g_C) \tensor_{U(\g_C)} \kk \cong S(\bs \g_{C})
\]
in $\dgco$. Here $S(\bs \g_{C})$ is the Chevalley-Eilenberg coalgebra of $\g_C$. As shown by Quillen in \cite[\Sec B.6]{Quillen}, the PBW theorem provides a local trivialization of $\pi_{U(\g_C)}$, and  
the corresponding triple represents the universal (trivialized) principal $U(\g_C)$ coalgebra bundle\footnote{In \cite{Quillen}, this coalgebra bundle is called the universal``principal $\g_C$-bundle''.}. 

In \cite[\Sec 2.2]{Prigge}, Prigge exhibited a universal coalgebra bundle with fiber $C$ by applying the following associated bundle construction to $U(\bs \g_C \dsum \g_C)$. Observe that, by definition, $\g_C$ acts naturally on $C$ via coderivations, giving $C$ the structure of a left $U(\g_C)$-module. The desired universal bundle is obtained, as before, by tensoring the counit $C \to \kk$ with
$U(\bs \g_C \dsum \g_C)$ over $U(\g_C)$:
\begin{equation} \label{eq:uni-cobun}
\pi_{C} \maps U(\bs \g_C \dsum \g_C) \tensor_{U(\g_C)} C \to S(\bs \g_{C}).
\end{equation}          
\newcommand{\fh}{\mathfrak{h}}
\begin{example}\label{ex:cb3}
Keeping the notation of the previous examples, let $C= S(\bs A[m])$. Recall that for any dg Lie algebra $\fh$, the PBW theorem \cite[Thm.\ B2.3]{Quillen} gives an isomorphism of dg coalgebras $S(\fh) \xto{\cong} U(\fh)$. Using this, along with the fact that the enveloping algebra of the abelian dg Lie algebra $\bs \g_{S(\bs A[m])}$ is simply $S(\bs \g_{S(\bs A[m])})$,
we obtain the following isomorphisms of graded coalgebras as in \cite[p.\ 291]{Quillen}:
\begin{align*}
U(\bs \g_{S(\bs A[m])}\dsum \g_{S(\bs A[m])}) &\tensor_{U(\g_{S(\bs A[m])})} S(\bs A[m]) \cong   \\
&\cong \Bigl( U(\bs \g_{S(\bs A[m])}) \tensor_\kk U(\g_{S(\bs A[m])}) \Bigr) \tensor_{U(\g_{S(\bs A[m])})} S(\bs A[m])\\
& \cong U(\bs \g_{S(\bs A[m])}) \tensor_\kk S(\bs A[m])\\ 
& \cong S(\bs \g_{S(\bs A[m])}) \tensor_\kk S(\bs A[m])\\
& \cong S(\bs \g_{S(\bs A[m])} \dsum \bs A[m])\\
& \cong S(\bs \LB{m} \dsum \bs A[m]) =  S(\bs \LE{m}).
\end{align*}
Note that the last isomorphism above follows from the isomorphism $\g_{S(\bs A[m])} \cong \LB{m}$ of dg Lie algebras exhibited in Example \ref{ex:cb2}, along with the definition of the dg Lie algebra $\LE{m}$. A direct calculation confirms that this sequence of coalgebra isomorphisms induces a commuting diagram in $\dgco$
\begin{equation} \label{diag:uni-cobun}
\begin{tikzdiag}{2}{4}
{
U(\bs \g_{S(\bs A[m])}\dsum \g_{S(\bs A[m])}) \tensor_{U(\g_{S(\bs A[m])})} S(\bs A[m])\& S(\bs \LE{m})  \\
S(\bs \g_{S(\bs A[m])})  \& S(\bs \LB{m})\\
};
\path[->,font=\scriptsize]
(m-1-1) edge node[auto] {$\cong$} (m-1-2)
(m-1-1) edge node[auto,swap] {$\pi_{S(\bs A[m])}$} (m-2-1)
(m-1-2) edge node[auto] {$P_{A(m)}$} (m-2-2)
(m-2-1) edge node[auto] {$\cong$} (m-2-2)
;
\end{tikzdiag}
\end{equation}
in which the horizontal maps are isomorphisms, and where $P_{A(m)}$ is the dg coalgebra morphism corresponding to the universal Lie $\infty$-algebra fibration $p_{A(m)} \maps \LE{m} \xto{} \LB{m}$, introduced in Sec.\ \ref{sec:uni-fib}, over the trivial base. 
\end{example}

\subsubsection{Twisting functions}
As in \cite[\Sec 5]{Quillen}, we denote by $\cT(B,\g_C)$ the set of twisting functions from the dg coalgebra $B$ to the coderivation dg Lie algebra $\g_C$. Recall that an element $\tha \in \cT(B,\g_C)$ is a degree $-1$ linear map $\tha \maps B \to \g_C$ which satisfies $\tha(1)=0$, along with the Maurer-Cartan equation
\[
d_{\g_C} \cc \tha + \tha \cc \del_B + \frac{1}{2} [-,-]_{\g_C} \cc (\tha \tensor \tha) \cc \Del_B =0. 
\] 
It follows from \cite[Prop.\ 6.2]{Quillen} that there is a bijection of sets
\begin{align} \label{eq:twist-iso}
\cT(B,\g_C) &\to^\cong \hom_{\dgco}(B, S(\bs \g_C))\\
\tha &\mapsto F_{\tha}.\nonumber
\end{align}
Using the notation Sec.\ \ref{sec:Linf-coalg}, $F_\tha$ may be characterized as the unique coalgebra morphism such that $F^1_\tha \vert_{\bar{B}} = \bs \tha$ as degree 0 linear maps in $\Hom_{\kk}(\bar{B},\bs \g_C)$. The compatibility of $F_\tha$ with the codifferentials on $B$ and the Chevalley-Eilenberg coalgebra of $\g_C$ is equivalent to $\tha$ satisfying the Maurer-Cartan equation.    

In \cite[Prop.\ 5.3]{Quillen}, Quillen established a bijection between $\cT(B,\g_C)$ and the set of isomorphism classes of principal $U(\g_C)$ coalgebra bundles over $B$ equipped with a fixed local trivialization. Prigge built on this result in \cite{Prigge} to obtain a classification theorem for coalgebra bundles. In what follows, we denote by $F^\ast \Bigl( U(\bs \g_C \dsum \g_C) \tensor_{U(\g_C)} C \Bigr)$ the pullback of the universal coalgebra bundle \eqref{eq:uni-cobun} along a morphism $F \maps B \to S(\bs \g_C)$ in $\dgco$. 
\begin{proposition}[Thm.\ 2.2.2 \cite{Prigge}] \label{prop:cobun-class}
Let $\tha \in \cT(B,\g_C)$ be a twisting function. The assignment
\[
\tha \mapsto F^\ast_{\tha} \Bigl( U(\bs \g_C \dsum \g_C) \tensor_{U(\g_C)} C \Bigr)
\]
induces a bijection of sets $\cT(B,\g_C) \xto{\cong} \Bun_B(C)$. 
\end{proposition}
The inverse to the isomorphism in Prop.\ \ref{prop:cobun-class} is explicitly given in
\cite[Prop.\ 2.2.4]{Prigge}. Let $(E,\pi)$ be a coalgebra bundle over $B$ with fiber $C$. Choose a local trivialization $\vphi \maps B \tensor C \xto{\cong} E$. There is a degree $-1$ linear map $\bar{\tha}_E \maps \bar{B} \to \g_C$ which sends $b \in \bar{B}$ to the coderivation $\bar{\tha}_E(b)$ on $C$ defined as 
\begin{equation} \label{eq:twist-form}
\bar{\tha}_E(b)(c):= \pr_C \cc \vphi^{-1} \cc \del_E \cc \vphi(b \tensor c).
\end{equation}  
Then $\bar{\tha}_E$ extends to a twisting function ${\tha}_E \maps {B} \to \g_C$ by setting
${\tha}_E(1)=0$. The proof that $(E,\pi, \vphi) \mapsto \tha_E$ indeed provides a well defined inverse to the pullback of the universal bundle is essentially identical to the proof given in \cite[Prop.\ 5.3]{Quillen} for the analogous statement concerning principal coalgebra bundles.    

\subsection{Proof of Proposition \ \ref{prop:classify}} \label{sec:class-proof}
Let us briefly recall the context behind the proposition. We are given a minimal fibration $f \maps (L,\el) \to (L',\el')$, and for $m\geq 1$, we consider the minimal fibration \eqref{eq:rPtowermap} 
$q^f_{\leq m} \maps \tleq{m}(L,f) \to \tleq{m-1}(L,f)$ appearing in the Postnikov tower for $f$.
By construction, $q^f_{\leq m}$ is surjective in degree 0, and its fiber is  $\ker q^f_{\leq m}=A[m]$, 
where $A=(\ker f)_m$. We choose a linear section \eqref{eq:eta} $\eta \maps \tleq{m-1}(L,f) \to \tleq{m}(L,f)$ of $q^f_{\leq m}$, which determines an isomorphism \eqref{eq:prA} $\phi \maps \tleq{m-1}(L,f) \dsum A[m] \xto{\cong} \tleq{m}(L,f)$ at the level of graded vector spaces. Hence, by taking $g = q^{f}_{\leq m}$ in Example \ref{ex:cb1}, we obtain a trivialized coalgebra bundle $\bigl(S (\bs \tleq{m-1}(L,f)), Q^{f}_{\leq m}, \vphi \bigr)$ with fiber $S(\bs A[m])$, where $Q^{f}_{\leq m}$ denotes the dg coalgebra morphism corresponding to $q^f_{\leq m}$. 

Hence, in order to prove both statements (1) and (2) of Prop.\ \ref{prop:classify}, 
it follows from Examples \ref{ex:cb2}, \ref{ex:cb3}, and Prop.\ \ref{prop:cobun-class}, that it is sufficient to show that the linear maps \eqref{eq:class1} $\psi_k \maps \Alt^k \tleq{m-1}(L,f) \to \End(A) \dsum A[m+1]$ agree with the formula for the twisting function $\bar{\tha}_{S (\bs \tleq{m}(L,f))}$ given in Eq.\ \ref{eq:twist-form}. We begin by using Eq.\ \ref{eq:morph_eq1}
to covert the $\psi_k$'s into a morphism between graded reduced coalgebras $\Psi \maps \S(\bs \tleq{m-1}(L,f)) \to \S (\bs \LB{m})$. We then apply the isomorphism \eqref{eq:twist-iso} between twisting functions and dg coalgebra morphisms to compare 
\[
\bs^{-1} \Psi^1 \maps  \S(\bs \tleq{m-1}(L,f)) \to \LB{m}            
\]
with $\bar{\tha}_{S (\bs \tleq{m-1}(L,f))}$. It follows from the definition of $\psi_k$, along with the formula \eqref{eq:struc_skew} for the coderivation $\del$ on $\S(\bs \tleq{m-1}(L,f))$ that we have
\[
\bs^{-1} \Psi^1_1( \bs x) = \bs^{-m-1} \cc \del^1_2\bigl(\bs x, \bs^{m+1}(-) \bigr)
\]
if $\deg{\bs x} =1$, and if $\sum_{i=1}^k \deg{\bs x_i} = m+2$:
\[
\bs^{-1} \Psi^1_k( \bs x_1, \ldots, \bs x_k) = \pr_{\bs A[m]} \cc \vphi^{-1} \cc \del^1_k\bigl(\bs x_1, \ldots, \bs x_k\bigr).
\] 
A direct calculation of  $\bar{\tha}_{S (\bs \tleq{m-1}(L,f))}(\bs x)$ and
$\bar{\tha}_{S (\bs \tleq{m-1}(L,f))}( \bs x_1, \ldots, \bs x_k)$ using formula \eqref{eq:twist-form}, along with the explicit description \eqref{eq:cb-split} for the local trivialization, gives the desired equality $ \bar{\tha}_{S (\bs \tleq{m}(L,f))} = \bs^{-1} \Psi^1$. \hfill \qed 

\starttocentries
\bibliographystyle{amsplain}
\bibliography{lie3}

\end{document}